\documentclass{article}
\usepackage{geometry}               
\usepackage{microtype}
\usepackage{graphicx,color}
\usepackage{subfigure}
\usepackage{booktabs}
\usepackage{amsmath}
\usepackage{amssymb}
\usepackage{mathtools}
\usepackage{natbib}
\usepackage{hyperref}
\usepackage{algorithm}
\usepackage{algpseudocode}
\usepackage{amsthm}

\newcommand{\A}{\bold{A}}
\newcommand{\ac}{\bold{a}}

\theoremstyle{plain}
\newtheorem{theorem}{Theorem}[section]
\newtheorem{proposition}[theorem]{Proposition}
\newtheorem{lemma}[theorem]{Lemma}
\newtheorem{corollary}[theorem]{Corollary}
\theoremstyle{definition}

\theoremstyle{remark}
\newtheorem{remark}[theorem]{Remark}

\begin{document}

\title{Generalized Continuous-Time Models for\\
Nesterov’s Accelerated Gradient Methods}

\author{
  Chanwoong Park\thanks{Department of Electrical and Computer Engineering, Seoul National University, Seoul 08826, South Korea.} \and
  Youngchae Cho\thanks{Automation and Systems Research Institute, Seoul National University, Seoul 08826, South Korea.} \and
  Insoon Yang\thanks{Department of Electrical and Computer Engineering, Automation and Systems Research Institute, Seoul National University, Seoul 08826, South Korea. Corresponding author. Email: insoonyang@snu.ac.kr}
}
\date{}




\maketitle

\begin{abstract}
Recent research has indicated a substantial rise in
interest in understanding Nesterov's accelerated gradient methods via their continuous-time models. However, most  existing studies focus on specific classes of Nesterov's methods, which hinders the attainment of an in-depth understanding and a unified perspective. 
To address this deficit, we present {\it generalized} continuous-time models that cover a broad range of Nesterov's methods, including those previously studied under existing continuous-time frameworks. Our key contributions are as follows. 
First, we identify the convergence rates of the generalized models, eliminating the need to determine the convergence rate for any specific continuous-time model derived from them. 
Second, we show that six existing continuous-time models are special cases of our generalized models, thereby positioning our framework as a unifying tool for analyzing and understanding these models. 
Third, we design a restart scheme for Nesterov's methods based on our generalized models and show that it ensures a monotonic decrease in objective function values. 
Owing to the broad applicability of our models, this scheme can be used to a broader class of Nesterov's methods compared to the original restart scheme.
Fourth, we uncover a connection between our generalized models and gradient flow in continuous time, showing that the accelerated convergence rates of our generalized models can be attributed to a time reparametrization in gradient flow.
Numerical experiment results are provided to support our theoretical analyses and results. 
\end{abstract}


\section{Introduction}
\label{sec:intro}

Nesterov's accelerated gradient methods are a class of optimization algorithms designed to speed up gradient-based methods \citep{nesterov1983method}. The key innovation in Nesterov's methods lies in the incorporation of a momentum term, which effectively anticipates the future trajectory of the iterates, thus enabling faster convergence. 
These methods have long been widely used as an effective tool for solving large-scale optimization problems across various domains, including machine learning~
\citep{dozat2016incorporating,Bottou2018,yang2022federated,li2023nesterov,xie2024adan}. 

Despite their widespread adoption, understanding why Nesterov's methods 
achieve accelerated convergence rates remains challenging, especially when relying on Nesterov's original proof, which 
 uses the estimate sequence technique. 
 To better understand Nesterov's methods, continuous-time models and analyses have recently emerged as an alternative~\citep{su2016differential,wibisono2016variational,betancourt2018symplectic,zhang2018direct,muehlebach2019dynamical,wilson2021lyapunov,shi2021understanding}. 
These models often provide valuable physical insights into specific forms of Nesterov's methods and have even led to the creation of new algorithms, such as a restart scheme. However, despite these advancements, to our knowledge, the literature lacks a systematic study that aims to unify the existing continuous-time models and offers a 
more fundamental perspective on Nesterov's methods as compared to individual research efforts.

\subsection{Contributions}
To bridge the existing gap in the literature and gain deeper insights into Nesterov's methods, we present novel ordinary differential equation (ODE) models that generalize a range of existing ODE models for Nesterov's methods. Our models are comprehensive in the sense that the class of Nesterov's methods they describe encompasses those covered by several existing ODE models from the literature, which we specify later in this paper.

To obtain our models, we first reformulate a three-sequence form of Nesterov's methods by slightly narrowing the class and rewriting the coefficients using an auxiliary sequence. 
Subsequently, we derive ODE models for the reformulated Nesterov's methods, explicitly accounting for the gradient at the future position. 
It is  shown that the flows generated by these ODE models converge to the iterates of the Nesterov methods of interest. Finally, generalizing the coefficients of these ODE models results in the proposed ODE models. 

We verify our assertion that the proposed models cover various Nesterov's methods previously studied in isolation by demonstrating that six ODE models from 
\citet{su2016differential,wilson2021lyapunov,muehlebach2019dynamical,shi2021understanding,chen2022revisiting}
are, in fact, special cases of our generalized models. 
This can be intuitively explained by the fact that our models inherit the flexibility of the Bregman Lagrangian flow (BLF) \citep{wilson2021lyapunov}, while incorporating the gradient correction similar to  
\citet{shi2021understanding,muehlebach2019dynamical,shi2021understanding,chen2022revisiting},
unlike the BLF. 
Thus, our models facilitate a unified analysis of 
these existing ODE models. For instance, we can recover the convergence rates of the six existing ODE models by identifying those of our generalized models. We conduct the convergence analysis for our models using a Lyapunov function.
 
Moreover, we develop a restart scheme for a  class of Nesterov's methods based on our generalized models.
Due to the broad applicability of our models, our restart scheme can be applied to a wider spectrum of Nesterov's methods compared to the restart scheme proposed by \citet{su2016differential}. 
We prove that the objective function values are monotonically decreasing along the iterates when using our proposed restart scheme. 

To further elucidate Nesterov's methods, we analyze our generalized models in the context of gradient flows. First, focusing on a certain class of our models, 
we rewrite them as an ODE of a single variable by introducing a time-dependent coefficient function. 
We then demonstrate that reparametrizing time through the coefficient function leads to a standard gradient flow model whose discrete-time counterpart corresponds to the gradient descent method. This indicates that 
we can control the convergence speed
by appropriately choosing the coefficient function, 
which partially reveals how our generalized models can achieve faster convergence rates compared to the standard gradient flow model in continuous time.

To support our theoretical analyses and results, we compare instances of our generalized models to existing ODE models through numerical experiments. The experiment results confirm that the added flexibility of our models allow them to capture the iterates of Nesterov’s methods more accurately than previous ODE models. 
We also provide empirical evidence demonstrating the effectiveness of our restart scheme and the connection between our models and gradient flow. 

Our main contributions are summarized as follows: 

\begin{itemize}
\item
We design novel ODE models for a range of Nesterov's methods that have only been studied separately using gradient correction and unified coefficients~(Section~\ref{sec:gen}). 

\item We identify the convergence rates of our generalized models, eliminating the need for individual convergence analyses for any ODE model  within our generalized framework (Section~\ref{sec:gen}). 

\item
We demonstrate that our generalized models cover six existing continuous-time models, thereby providing a unified framework to analyze each of these models (Section~\ref{sec:ex}). 

\item
We develop a restart scheme for Nesterov's methods  using our models. This is more general compared to the restart scheme proposed by \citet{su2016differential} and guarantees a monotonic decrease in objective function values (Section~\ref{sec:restart}).

\item 
Using a time reparametrization technique, we show that our models reduce to the standard gradient flow that describes gradient descent. This implies that \emph{adjusting the speed of time} induces acceleration in continuous time (Section~\ref{sec:time}). 

\end{itemize}

\subsection{Related work}

First proposed by \citet{nesterov1983method}, 
Nesterov's accelerated gradient methods continue to be widely used as optimal first-order methods \citep{nemirovskij1983problem}. \citet{nesterov2003introductory} introduced the estimate sequence technique to analyze these methods. Based on this discrete-time analysis, several variations and extensions of Nesterov's methods have been proposed and analyzed~\citep{nesterov2005smooth, nesterov2008accelerating, tseng2008accelerated, beck2009fast, baes2009estimate}. 

However, the estimate sequence is algebraically complicated, making it difficult to attain an intuitive understanding of why Nesterov' methods achieve acceleration.  Numerous approaches have been explored to comprehend these methods better. \citet{allen2014linear} analyzed Nesterov's methods through the lens of linear coupling between gradient descent and mirror descent. \citet{bubeck2015geometric} offered a geometric interpretation of Nesterov's methods. \citet{lessard2016analysis} and \citet{hu2017dissipativity} proposed control theoretic frameworks using dissipativity. 
Additionally, \citet{fazlyab2018analysis, sanz2021connections} developed a linear matrix inequality framework based on integral quadratic constraints  to analyze Nesterov's methods. 
Another closely related approach is using \emph{performance estimation problems} (PEP) to analyze the convergence rate of first-order methods~\citep{drori2014performance, taylor2017smooth}.
The PEP framework was generalized to analyze continuous-time models for first-order optimization methods as well~\citep{moucer2023systematic,kim2024convergence}.   

As mentioned earlier, there has been significant interest in describing and explaining Nesterov's methods through continuous-time models. 
\citet{su2016differential} introduced an ODE model for a popular form of Nesterov's methods and developed an efficient restart scheme using this model. 
This ODE model was further extended to the mirror descent setting by \citet{krichene2015accelerated} 
and was used to design fast optimization algorithms in \citet{attouch2018fast, attouch2019fast}.
\citet{wibisono2016variational} and \citet{wilson2021lyapunov} proposed and analyzed the BLF,  a general ODE model for various accelerated gradient methods. 
\citet{shi2019acceleration, shi2021understanding} modeled Nesterov's methods as high-resolution ODEs, including a Hessian term through an elaborate approximation. 
The models were extended by \citet{chen2022revisiting} to accommodate varying step sizes. 
\citet{chen2022gradient} revealed that Nesterov's methods use a gradient correction term that implicitly contains velocity, achieving convergence rates of $o(1/k^3)$ for gradient norms and $o(1/k^2)$ for the objective value minimization. 
\citet{diakonikolas2019approximate} introduced a technique using a duality gap to analyze first-order methods, linking continuous-time methods with discretization error handling. 
\citet{chen2021unified} proposed the concept of ``strong Lyapunov conditions" to provide a unified convergence analysis for first-order convex optimization.
Without relying on vanishing step sizes, \citet{berthier2021continuized} presented a stochastic ODE by taking gradient steps at random times. 
\citet{suh2022continuous} analyzed an ODE for accelerated gradient methods using conservation laws in a new coordinate system, simplifying existing analyses and achieving an $O(1/k^2)$ convergence rate. 
\citet{kim2023unifying} developed a novel Lagrangian, continuous-time models and discrete-time algorithms that handle both convex and strongly convex objective functions in a unified way. 
\citet{toyoda2024unified} introduced a unified Euler--Lagrange system and its corresponding Lyapunov function to analyze accelerated gradient methods, covering various existing convergence rates. 
\citet{maskan2024variational} proposed a variational perspective connecting high-resolution ODEs to the Euler--Lagrange equation.
\citet{adly2024accelerated} derived high-resolution dynamics from Nesterov and Ravine methods and linked them to fast convergence and Levenberg--Marquardt regularization. 
The continuous-time interpretations have been extended to accelerated first-order methods for Riemannian optimization as well~\citep{alimisis2020continuous, Kim2022}.

Notably, in continuous-time analysis, discretization methods that connect ODEs and optimization algorithms play an important role. \citet{betancourt2018symplectic} identified a rate-matching discretization for the BLF model \citep{wibisono2016variational,wilson2021lyapunov} by replacing the Bregman Lagrangian with a Hamiltonian.
Similarly, \citet{muehlebach2021optimization} and \citet{francca2021dissipative} achieved rate-matching discretization using symplectic integrators based on Lagrangian or Hamiltonian formulations. 
\citet{zhang2018direct} proposed a new optimization method by applying the Runge--Kutta integrator to ODE models, while 
\citet{scieur2017integration} derived accelerated algorithms by applying multi-step discretization to the basic gradient flow. 
Meanwhile, \citet{muehlebach2019dynamical} derived Nesterov's methods by applying the semi-implicit Euler method to their novel ODE models. 
\citet{ushiyama2023unified} introduced a unifying discretization framework using the so-called ``weak discrete gradient" to simplify the development and analysis of algorithms. 
\citet{muehlebach2023accelerated} designed a new accelerated algorithm for constrained optimization by discretizing accelerated gradient flow, demonstrating both convergence and efficiency with velocity-based constraints.

Many papers have studied restarting techniques for  Nesterov's methods as, which prevent the algorithms from stagnating in suboptimal convergence patterns and reduce oscillatory behaviors, leading to faster and more robust convergence. \citet{o2015adaptive} proposed a heuristic adaptive technique that resets the momentum whenever the objective value strictly increases. Extending this approach, \citet{su2016differential} developed a restart scheme framed within the context of ODEs to maintain a high velocity along the trajectory of the objective value. 
The authors further established a linear convergence rate for its continuous-time counterpart. Moreover, \citet{aujol2022fista} designed a restart scheme for the fast iterative shrinkage-thresholding algorithm (FISTA) presented by \citet{beck2009fast}. \citet{maulen2023speed} analyzed the effects of the restart scheme of \citet{su2016differential} on inertial dynamics with Hessian-driven damping proposed by \citet{attouch2016fast}, laying the groundwork for accelerating Hessian-driven inertial algorithms.

\section{Preliminaries}

\subsection{Convex optimization}

We consider the problem of minimizing a continuously differentiable convex function $f: \mathbb{R}^n \to \mathbb{R}$, that is,
\begin{equation}\label{eqn:opt}
\min_{x \in \mathbb{R}^n} f(x).
\end{equation}
Throughout the paper, we assume that $f$ is \emph{L-smooth}. In other words,
$\nabla f(x)$ is $L$-Lipschitz continuous, that is, 
\begin{equation*}
    \|\nabla f(x) - \nabla f(y)\| \leq L\|x-y\| \quad \forall x,y\in\mathbb{R}^n.
\end{equation*}
For convenience, we assume that the objective function has a unique minimum at $x^* = 0$ and the optimal value is $0$, that is $f(0) = 0$. 
For $\mu > 0$, the function $f$ is said to be \emph{$\mu$-strongly convex} if 
\[
f(y) \geq f(x) + \left \langle \nabla f(x), y-x \right \rangle + \frac{\mu}{2} \|y-x\|^2. 
\]
Note that $f$ is \emph{convex}  when the strong convexity parameter $\mu$ is zero.

\subsection{Bregman divergences}

The Bregman divergence is a popular tool to handle optimization methods in non-Euclidean settings. 
Given a continuously differentiable strictly convex function $g: \mathbb{R}^n \to \mathbb{R}$, 
the \emph{Bregman divergence} 
$D_g:\mathbb{R}^n \times \mathbb{R}^n \to [0, \infty)$ of $g$ is defined as
\[
D_g (y,x) = g(y) - g(x) - \left \langle \nabla g(x), y-x \right \rangle. 
\]
In particular, when $g(x) = \frac{1}{2}\|x\|^2$, 
we have $D_g (y,x) = \frac{1}{2} \| y-x\|^2$. 

The function $f$ is said to be \emph{$\mu$-uniformly convex with respect to $g$} provided that 
\[
D_f (x, y) \geq \mu D_g (x, y) \quad \forall x, y \in \mathbb{R}^n.
\]
When $g(x) = \frac{1}{2}\|x\|^2$, this condition reduces to that for the $\mu$-strong convexity of $f$. 

\subsection{Nesterov's method}

Nesterov's method for solving problem~\eqref{eqn:opt}, possibly with $\mu = 0$, can be written in the following three-sequence form:
\begin{equation}\label{eqn:NA}
\begin{split}
y_k&= x_k + \frac{\theta_k \gamma_k}{\gamma_k + \mu \theta_k} (z_k - x_k)\\
x_{k+1} &= y_k - s_k \nabla f(y_k)\\
z_{k+1} &= x_k + \frac{1}{\theta_k} (x_{k+1} - x_k),
\end{split}
\end{equation}
where $x_0 = y_0 = z_0$, $s_k \geq 0$ and $0<\theta_k<1$.
Also, $\theta_k$, $\gamma_k$ and $s_k$ satisfy
\begin{equation}\label{eqn:para}
 \gamma_{k+1} = (1-\theta_k) \gamma_k + \mu \theta_k = \frac{\theta_k^2}{s_k}.
\end{equation}
The scheme can be expressed in the following two-sequence form as well:
\begin{equation}\label{eqn:NAtwo}
    \begin{split}
            y_k&=x_k + b_k(x_k-x_{k-1})\\
            x_{k+1} &= y_k - s_k \nabla f(y_k),
    \end{split}
\end{equation}
where 
\begin{equation*}
        b_k=\frac{\theta_k \gamma_k}{\gamma_k + \mu \theta_k}\frac{1-\theta_{k-1}}{\theta_{k-1}}.
\end{equation*}

If $0 < s_k \leq \frac{1}{L}$ for all $k$, then Nesterov's method achieves $ f(x_k) = O \big (\prod_{i=0}^{k-1} (1-\theta_i) \big)$ ~\citep[Theorem 2.2.1]{nesterov2018lectures}. Furthermore, for the frequently used constant step size $s_k=\frac{1}{L}$, we have $b_k \approx \frac{k}{k+3}$ if $\mu=0$~\citep[Algorithm 2]{tseng2008accelerated} and $b_k = \frac{\sqrt{L}-\sqrt{\mu}}{\sqrt{L}+\sqrt{\mu}}$ if $\mu>0$. 
The detailed derivation of this method can be found in \citet{nesterov2018lectures}.

\section{Generalized Continuous-Time Models for Nesterov's Methods}\label{sec:gen}

We derive generalized continuous-time models for Nesterov-type methods that preserve the lookahead-gradient structure and admit a unified Lyapunov analysis.

\subsection{Discrete reformulation via the growth sequence $\{A_k\}$}\label{sec:gen:discrete}
Following the idea suggested by \citet{wilson2021lyapunov}, we introduce an auxiliary sequence $A_k\in\mathbb{R}$ and express the coefficients of \eqref{eqn:NA} in terms of $A_k$.
The sequence $\{A_k\}$ serves as a rate-encoding parameter for Nesterov-type schemes.

Specifically, we rewrite \eqref{eqn:NA} in the following three-sequence form:
\begin{equation}\label{eqn:NA2}
\begin{split}
y_k&=x_k + a_k(z_k-x_{k})\\
x_{k+1} &= y_k - s_k \nabla f(y_k)\\
z_{k+1} &= x_k + \frac{1}{\theta_k}(x_{k+1}-x_k),
\end{split}
\end{equation}
with $\theta_k=(A_{k+1}-A_k)A^{-1}_{k+1}$ and
\begin{equation}\label{eqn:s}
a_k=
\begin{cases}
\begin{aligned}
&\frac{A_{k+1}-A_k}{A_{k+1}}&&\text{if }\mu=0\\
&\frac{A_{k+1}-A_k}{2A_{k+1}-A_k}&&\text{if }\mu>0
\end{aligned}
\end{cases}
\quad\text{and}\quad
s_k=
\begin{cases}
\begin{aligned}
&\frac{(A_{k+1}-A_k)^2}{A_{k+1}}&&\text{if }\mu=0\\
&\frac{(A_{k+1}-A_k)^2}{\mu A_{k+1}^2}&&\text{if }\mu>0.
\end{aligned}
\end{cases}
\end{equation}
In Appendix \ref{app:reform}, we provide a detailed derivation of \eqref{eqn:NA2} and discuss its convergence rate.\footnote{While the introduction of $A_k$ is attributed to \citet{wilson2021lyapunov}, \citet{wilson2021lyapunov} used a Lyapunov-based analysis to prove convergence of \eqref{eqn:NA2}, which differs from our approach.}
We refer to \eqref{eqn:NA2} with $\mu=0$ and $\mu>0$ as (NAG-C) and (NAG-SC), respectively.

It is worth noting that \eqref{eqn:NA2} covers a specific range of \eqref{eqn:NA} because $\gamma_0$ is confined to the specific values $1/A_0$ and $\mu$ in its derivation.
Nonetheless, \eqref{eqn:NA2} is general enough to encompass Nesterov's methods with constant step sizes,
as modeled by \citet{su2016differential, shi2019acceleration, shi2021understanding, muehlebach2019dynamical, berthier2021continuized},
and their variants studied by \citet[Section 4.4]{d2021acceleration}.
This flexibility arises because we can choose $\{A_k\}$ more freely than in existing models, where $\{A_k\}$ is typically restricted to a fixed inductive form.

\subsection{Continuous-time limits: (ODE-C) and (ODE-SC)}\label{sec:gen:ode}
We now formulate ODE models corresponding to (NAG-C) and (NAG-SC), which will serve as concrete instances of our generalized framework.
Through basic algebraic manipulations, for a constant $h>0$, (NAG-C) and (NAG-SC) are respectively transformed into
\begin{equation}\label{eqn:C-h}
\begin{split}
y_k&=x_k + a_k(z_k-x_k)\\
\frac{z_{k+1}-z_k}{h} &= -\frac{A_{k+1}-A_k}{h}\nabla f(y_k) \\
\frac{x_{k+1}-x_k}{h} &= \frac{A_{k+1}-A_k}{h A_k} (z_{k+1}-x_{k+1}),
\end{split}
\end{equation}
and
\begin{equation}\label{eqn:SC-h}
\begin{split}
y_k &= x_k + a_k (z_k-x_k)\\
\frac{z_{k+1}-z_k}{h} &= -\frac{A_{k+1}-A_k}{h A_{k+1}}\bigg(z_k-y_k+\frac{1}{\mu}\nabla f(y_k)\bigg) \\
\frac{x_{k+1}-x_k}{h} &= \frac{A_{k+1}-A_k}{h A_k} (z_{k+1}-x_{k+1}).
\end{split}
\end{equation}
By introducing $t=hk$ and applying the Euler method, we obtain the coupled ODE systems
\begin{equation}\label{eqn:ODE-C}\tag{ODE-C}
\begin{split}
Y &= X+\ac (t)(Z-X)\\
\dot{Z} &= -\dot{\A}(t)\nabla f(Y) \\
\dot{X} &= \frac{\dot{\A}(t)}{\A(t)}(Z-X),
\end{split}
\end{equation}
when $\mu=0$, and
\begin{equation}\label{eqn:ODE-SC}\tag{ODE-SC}
\begin{split}
Y &= X + \ac(t) (Z-X)\\
\dot{Z} &= -\frac{\dot{\A}(t)}{\A(t)} \bigg(Z-Y+\frac{1}{\mu}\nabla f(Y)\bigg) \\
\dot{X} &= \frac{\dot{\A}(t)}{\A(t)}(Z-X),
\end{split}
\end{equation}
when $\mu>0$, where $X(t),Y(t),Z(t):[0,\infty)\to\mathbb{R}^n$, $X(0)=Y(0)=Z(0)=x_0$,
$\A:[0,\infty)\to(0,\infty)$ and $\ac:[0,\infty)\to[0,1]$.
We define $\A(hk)=A_k$ and $\ac(hk)=a_k$.

The ODE model proposed by \citet{wilson2021lyapunov} to describe (NAG-C) and (NAG-SC) shares a similar form with \eqref{eqn:ODE-C} and \eqref{eqn:ODE-SC}.
The key difference is that \citet{wilson2021lyapunov} uses $\nabla f(X)$, whereas \eqref{eqn:ODE-C} and \eqref{eqn:ODE-SC} employ $\nabla f(Y)$.
In this sense, \eqref{eqn:ODE-C} and \eqref{eqn:ODE-SC} are novel ODE models that preserve the lookahead-gradient structure in continuous time.

\subsection{Generalized ODE class (G-ODE) and unified Lyapunov analysis}\label{sec:gen:gode}
To enable consolidated analyses across a broader family than specific instances such as \eqref{eqn:ODE-C}--\eqref{eqn:ODE-SC},
we further generalize them as follows:
\begin{equation}\label{eqn:G-ODE-C}\tag{G-ODE-C}
    \begin{split}
        Y&=X+\ac(t)(Z-X) \\
        Z&=X+e^{-\alpha(t)}\dot{X} \\
        \frac{d}{dt}\nabla g (Z) &= -e^{\alpha(t)+\beta(t)}\nabla f(Y)
    \end{split}
\end{equation}
when $f$ is convex, and
\begin{equation}\label{eqn:G-ODE-UC}\tag{G-ODE-UC}
    \begin{split}
        Y&=X+\ac (t)(Z-X) \\
        Z&=X+e^{-\alpha(t)}\dot{X} \\
        \frac{d}{dt}\nabla g(Z)&=-\dot{\beta}(t)(\nabla g(Z)-\nabla g(Y)) - \frac{e^{\alpha(t)}}{\mu}\nabla f(Y)
    \end{split}
\end{equation}
when $f$ is $\mu$-uniformly convex with respect to $g$, where $X(0)=Y(0)=Z(0)=x_0$.
The coefficient functions $\alpha,\beta:[0,\infty)\to\mathbb{R}$ are assumed to satisfy
\begin{equation}\label{eqn:condition}
e^{\alpha (t)} \geq \dot{\beta}(t) >0, \quad \beta(t) \to \infty \mbox{ as $t \to \infty$}.
\end{equation}
These conditions are similar to the ideal scaling conditions used by \citet{wibisono2016variational, wilson2021lyapunov}.
The generalized ODE models are determined by four functions $\ac,g,\alpha$, and $\beta$.
In particular, $\ac(t)$ can be any measurable function with $0\leq \ac(t)\leq 1$, so that $Y$ is a convex combination of $X$ and $Z$.

As with \eqref{eqn:ODE-C} and \eqref{eqn:ODE-SC}, it is important to note that \eqref{eqn:G-ODE-C} and \eqref{eqn:G-ODE-UC} evaluate $\nabla f$ at $Y$, unlike \citet{wibisono2016variational, wilson2021lyapunov, berthier2021continuized}, which use $\nabla f(X)$.
Since $Y=X+\ac(t)e^{-\alpha(t)}\dot{X}$ can be interpreted as a future position of $X$, incorporating $\nabla f(Y)$ yields an ODE model more aligned with the lookahead mechanism of Nesterov's methods.

We now state our first main result: a unified Lyapunov analysis yielding an $O(e^{-\beta(t)})$ convergence rate.
Define $V:\mathbb{R}^n\times\mathbb{R}^n\times[0,\infty)\to\mathbb{R}$ by
\begin{equation*}
V(X,Z,t)=
\begin{cases}
D_g(0,Z)+e^{\beta(t)}f(X), & \text{if }\mu=0,\\
e^{\beta(t)}\big(\mu D_g(0,Z)+f(X)\big), & \text{if }\mu>0.
\end{cases}
\end{equation*}

\begin{theorem}\label{thm:conv}
Let $X(t), Y(t), Z(t): [0,\infty) \rightarrow \mathbb{R}^n$ be a solution to \eqref{eqn:G-ODE-C} when $f$ is convex, and a solution to \eqref{eqn:G-ODE-UC} when $f$ is $\mu$-uniformly convex with respect to $g$.
Then, the energy function $\mathcal{E}(t)=V(X(t),Z(t),t)$ is monotonically non-increasing, and
\[
f(X(t)) \leq e^{-\beta(t)} \mathcal{E}(0) \quad \forall t \geq 0.
\]
\end{theorem}
Its proof can be found in Appendix~\ref{app:conv}.
Theorem~\ref{thm:conv} implies that one does not need to perform convergence analyses on a case-by-case basis for an ODE model if it falls within \eqref{eqn:G-ODE-C} or \eqref{eqn:G-ODE-UC}.
In this sense, the primary advantage of our models is that they provide a unified tool for analyzing various ODE models for accelerated first-order methods.

\subsection{Why include the auxiliary variable $Y$? A discretization gap in BLF}\label{sec:blf-gap}
To further motivate why we keep the auxiliary variable $Y$ explicitly in continuous time, we revisit the BLF discretization route of \citet{wilson2021lyapunov} and identify a discrete-time auxiliary-sequence gap.

Throughout this subsection, we focus on the convex case ($\mu=0$), the Euclidean geometry $g(x)=\frac12\|x\|^2$,
and the ideal scaling $e^{\alpha(t)}=\dot\beta(t)$.
Under these conditions, BLF can be reduced to the following two-state system:
\begin{equation}\label{eqn:blf-reduced}
\begin{split}
Z &= X + \frac{e^{\beta}}{\frac{d}{dt}e^{\beta}}\,\dot X,\\
\dot Z &= -\Big(\frac{d}{dt}e^{\beta}\Big)\nabla f(X).
\end{split}
\end{equation}
Notably, the continuous-time BLF evolves only the pair $(X,Z)$.

Approximating $e^{\beta(hk)}$ by a positive increasing sequence $\{A_k\}$ and $\frac{d}{dt}e^{\beta}$ by
$\frac{A_{k+1}-A_k}{h}$, a symplectic discretization of \eqref{eqn:blf-reduced} yields the following updates,
referred to as the first family in \citet{wilson2021lyapunov}:
\begin{equation}\label{eqn:blf-first-family}
\begin{split}
x_{k+1} &= \frac{A_{k+1}-A_k}{A_{k+1}} z_k + \frac{A_k}{A_{k+1}} y_k,\\
z_{k+1} &= z_k - (A_{k+1}-A_k)\nabla f(x_{k+1}).
\end{split}
\end{equation}
A crucial observation is the appearance of the auxiliary sequence $\{y_k\}$,
which has no direct counterpart in the continuous BLF dynamics \eqref{eqn:blf-reduced}.

\begin{proposition}[Auxiliary-sequence gap]\label{prop:blf-gap}
The recursion \eqref{eqn:blf-first-family} is not closed on the state pair $(x_k,z_k)$:
given $(x_0,z_0)$ and $\{A_k\}$, the next iterate $(x_{k+1},z_{k+1})$ is not uniquely determined
unless an additional evolution rule for $\{y_k\}$ is specified.
Equivalently, \eqref{eqn:blf-first-family} defines a family of discrete-time methods parameterized by the choice of $\{y_k\}$.
\end{proposition}
\begin{proof}
The update for $x_{k+1}$ in \eqref{eqn:blf-first-family} depends explicitly on $y_k$.
To see non-uniqueness, take the convex function $f\equiv 0$ so that $\nabla f\equiv 0$.
Then \eqref{eqn:blf-first-family} reduces to
$z_{k+1}=z_k$ and
$x_{k+1}=\frac{A_{k+1}-A_k}{A_{k+1}}z_k+\frac{A_k}{A_{k+1}}y_k$,
so $z_k\equiv z_0$ for all $k$, while $x_{k+1}$ can be altered arbitrarily by choosing $y_k$ arbitrarily.
\end{proof}

To recover a specific accelerated method from the family \eqref{eqn:blf-first-family},
\citet{wilson2021lyapunov} impose an additional update rule for the auxiliary sequence $\{y_k\}$:
\begin{equation}\label{eqn:blf-closing-wilson}
y_{k+1}=x_k-s_k\nabla f(x_k).
\end{equation}
For stating an exact algebraic equivalence with our three-sequence form, it is convenient to adopt the following
index-shifted but equivalent closing rule:
\begin{equation}\label{eqn:blf-closing-shift}
y_{k+1}=x_{k+1}-s_k\nabla f(x_{k+1}).
\end{equation}

\begin{proposition}[Closing \eqref{eqn:blf-first-family} yields \eqref{eqn:NA2}]\label{prop:blf-to-na2}
Assume the convex case $\mu=0$ and the coefficient relations in \eqref{eqn:s}, so that
$\theta_k=\frac{A_{k+1}-A_k}{A_{k+1}}$, $a_k=\theta_k$, and $s_k=\frac{(A_{k+1}-A_k)^2}{A_{k+1}}$.
Consider \eqref{eqn:blf-first-family} together with the closing rule \eqref{eqn:blf-closing-shift},
and initialize $x_0=y_0=z_0$.
Define the relabeled sequences
\begin{equation}\label{eqn:blf-relabel}
\tilde x_k := y_k,\quad \tilde y_k := x_{k+1},\quad \tilde z_k := z_k.
\end{equation}
Then $(\tilde x_k,\tilde y_k,\tilde z_k)$ satisfies \eqref{eqn:NA2} (with $\mu=0$).
\end{proposition}
\begin{proof}
(1) Using \eqref{eqn:blf-first-family} and $a_k=\theta_k$, we have
\begin{equation}\label{eq:blf-proof-1}
\tilde y_k=x_{k+1}=\theta_k z_k+(1-\theta_k)y_k=\tilde x_k+a_k(\tilde z_k-\tilde x_k).
\end{equation}

(2) By \eqref{eqn:blf-closing-shift},
\begin{equation}\label{eq:blf-proof-2}
\tilde x_{k+1}=y_{k+1}=x_{k+1}-s_k\nabla f(x_{k+1})=\tilde y_k-s_k\nabla f(\tilde y_k).
\end{equation}

(3) From \eqref{eqn:blf-first-family},
\[
\tilde z_{k+1}-\tilde z_k=z_{k+1}-z_k=-(A_{k+1}-A_k)\nabla f(x_{k+1}).
\]
Multiplying by $\theta_k=(A_{k+1}-A_k)/A_{k+1}$ yields
\[
\theta_k(\tilde z_{k+1}-\tilde z_k)
=-\frac{(A_{k+1}-A_k)^2}{A_{k+1}}\nabla f(\tilde y_k)
=-s_k\nabla f(\tilde y_k).
\]
Combining this with \eqref{eq:blf-proof-1}--\eqref{eq:blf-proof-2} and rearranging yields
\begin{equation}\label{eq:blf-proof-3}
\tilde z_{k+1}=\tilde x_k+\frac{1}{\theta_k}\big(\tilde x_{k+1}-\tilde x_k\big).
\end{equation}
\end{proof}

\begin{remark}[Indexing convention for the auxiliary update]
The closing rule is often written as \eqref{eqn:blf-closing-wilson}.
We adopt \eqref{eqn:blf-closing-shift} because it yields a direct algebraic equivalence with our three-sequence form.
\end{remark}

In summary, BLF evolves only the two-state pair $(X,Z)$, whereas its discretization introduces an auxiliary sequence $\{y_k\}$ that must be specified externally to obtain a closed accelerated algorithm.
In contrast, \eqref{eqn:ODE-C}--\eqref{eqn:ODE-SC} and \eqref{eqn:G-ODE-C}--\eqref{eqn:G-ODE-UC} incorporate the auxiliary variable $Y(t)$ directly in continuous time, so the discrete three-sequence structure \eqref{eqn:NA2} arises naturally without auxiliary-variable injection.

\subsection{Consistency with Nesterov's accelerated methods}\label{sec:gen:consistency}
Meanwhile, we can derive \eqref{eqn:ODE-C} and \eqref{eqn:ODE-SC} from \eqref{eqn:G-ODE-C} and \eqref{eqn:G-ODE-UC}, respectively, by letting
\begin{equation*}
    e^{\alpha(t)}=\frac{\dot{\A}(t)}{\A(t)}, \quad \beta(t)=\ln \A(t), \quad g(x)=\frac{1}{2}\|x\|^2.
\end{equation*}
The following proposition proves that the continuous-time flows generated by \eqref{eqn:ODE-C} and \eqref{eqn:ODE-SC} converge to the iterates of (NAG-C) and (NAG-SC), respectively.
\begin{proposition}\label{thm:Nes-Flow}
Assume that $\A$ satisfies:
\begin{enumerate}
  \item $\A(t)$ is differentiable.
  \item $\displaystyle{\lim_{t\to\infty}} \A(t)=\infty$, $\dot{\A}(t)>0$, and $\A(hk)=A_k$.
  \item $\A(t)$ is independent of $h$.
\end{enumerate}
Furthermore, assume that $\ac$ is a continuous function such that $\ac(hk)=a_k$.
Then, \eqref{eqn:ODE-C} and \eqref{eqn:ODE-SC} respectively have a unique solution $(X(t),Y(t),Z(t))$ such that
\[
 \lim_{h\rightarrow 0} \max_{0\leq k \leq \left \lfloor \frac{T}{h}\right \rfloor} \|x_k-X(hk)\| = 0,
\]
where $\{x_k\}$ is generated by \eqref{eqn:NA2}.
\end{proposition}
For its proof, refer to Appendix~\ref{app:Nes-Flow}. Moreover, examples of valid $\A(t)$ and $\ac(t)$ can be found in Section~\ref{sec:ex}.
Together with Theorem~\ref{thm:conv}, Proposition~\ref{thm:Nes-Flow} confirms that our generalized ODE models can accommodate continuous-time counterparts of \eqref{eqn:NA2}.

In Section~\ref{sec:ex}, we show that \eqref{eqn:G-ODE-C} and \eqref{eqn:G-ODE-UC} cover not only \eqref{eqn:ODE-C} and \eqref{eqn:ODE-SC} but also various continuous-time models in the literature.

\section{Relationships with Existing ODE Models}\label{sec:ex}
\begin{figure}[t]
	\begin{center}\centerline{\includegraphics[width=0.5\columnwidth]{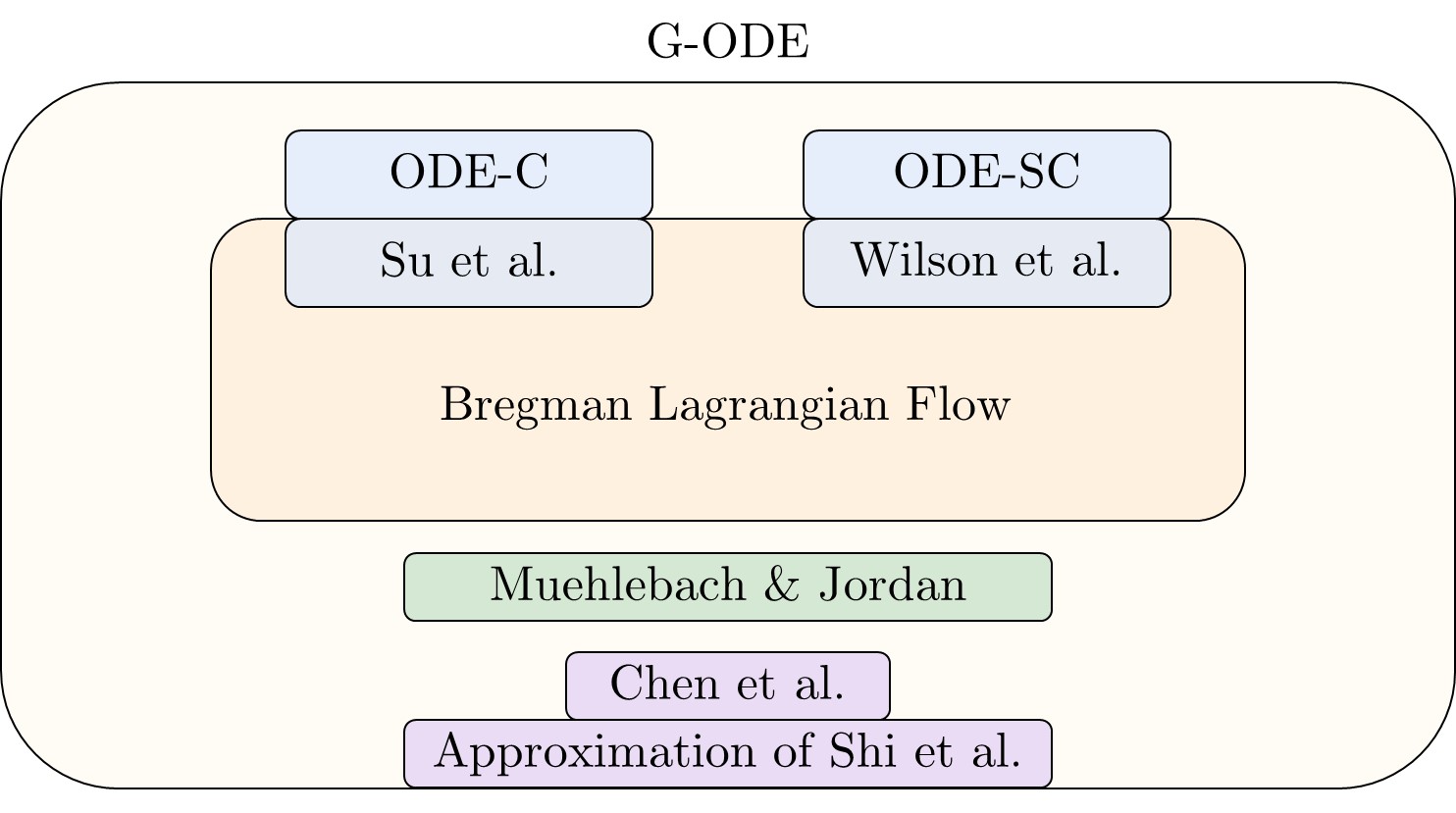}}
\caption{Relationships between our models and the six existing models. In the case of 
\citet{shi2021understanding},
our generalized models contain an approximation of it.}
		\label{fig:diagram}
	\end{center}
\end{figure}

In this section, as our second main contribution, we explain how existing continuous-time models for Nesterov's methods are encompassed within our framework. 
Specifically, we derive six existing models developed by \citet{su2016differential,wilson2021lyapunov,muehlebach2019dynamical,shi2021understanding,chen2022revisiting} from our generalized models. 
Some of these existing models arise directly from the special cases \eqref{eqn:ODE-C} and \eqref{eqn:ODE-SC}.
Moreover, the model proposed by  \citet{shi2021understanding} is shown to fit within our models  after applying approximations. Figure~\ref{fig:diagram} illustrates the relationships between our models and the six existing models. 

\subsection{ODE-C and the ODE model 
of \citet{su2016differential}
}\label{sec:comp:su}

Suppose first that the objective function $f$ is  convex. We consider the ODE model~\eqref{eqn:ODE-C} with 
 the coefficient function $\A(t)$ chosen as
\begin{equation*}
    \A(t)=\frac{(t+\varepsilon)^2}{4L}
\end{equation*}
for some constant $\varepsilon>0$. 
It satisfies the three conditions on $\A$ in Proposition \ref{thm:Nes-Flow} with 
\[
A_k=\frac{(hk+\varepsilon)^2}{4L}.
\]
According to \eqref{eqn:s}, the coefficients for (NAG-C) are given by
\begin{equation}\nonumber
\begin{split}
    a_k &= \frac{A_{k+1}-A_k}{A_{k+1}}=h\frac{2hk+2\varepsilon+h}{(hk+\varepsilon+h)^2}\\
    s_k&=\frac{(A_{k+1}-A_k)^2}{A_{k+1}}=\frac{h^2}{4L}\frac{(2hk+2\varepsilon+h)^2}{(hk+\varepsilon+h)^2}\leq\frac{h^2}{L}.
    \end{split}
\end{equation}
Thus, whenever $h \leq 1$, the stepsize condition $s_k \leq 1/L$ for convergence is satisfied. 
We also let
\[
\ac(t)= h\frac{2(t+\varepsilon)+h}{(t+\varepsilon+h)^2}
\] 
to ensure that $\ac(hk) = a_k$ as desired. 

With this choice of coefficient functions $\A(t)$ and $\ac (t)$, \eqref{eqn:ODE-C} can be expressed as
\begin{equation} \label{eqn:ODE-CV-h}
    \ddot{X} + \frac{3}{t+\varepsilon}\dot{X}+\frac{1}{L}\nabla f\bigg (X+h\frac{(t+\varepsilon+\frac{h}{2})(t+\varepsilon)}{(t+\varepsilon+h)^2}\dot{X}\bigg )=0.
\end{equation}
Letting both $h$ and $\varepsilon$ tend to zero yields the following second-order ODE:
\begin{equation}\label{eqn:sec}
    \ddot{X}(t) + \frac{3}{t}\dot{X}(t)+\frac{1}{L}\nabla f(X(t))=0.
\end{equation}
With the change of variables $\tilde{X}(t) = X (\sqrt{L} t )$, \eqref{eqn:sec} is equivalent to the ODE model proposed by 
\citet{su2016differential}; see Equation (3) in 
\citet{su2016differential}.
It is straightforward to verify that \eqref{eqn:sec} is included in \eqref{eqn:G-ODE-C}. 

At the other extreme, we consider the case where $h = 1$, departing from the discussion on the relationship between our models and the model  
of \citet{su2016differential}.
For a sufficiently large $k$, if $\varepsilon\to0$ and $t=k$, the coefficients $s_k$ and $b_k$ used in \eqref{eqn:NAtwo} can be approximated as 
\begin{equation}\label{eqn:Nes-CV-sb}
\begin{split}
 s_k &=\frac{1}{4L}\frac{(2k+1)^2}{(k+1)^2}\approx\frac{1}{L}\\
 b_k&=\frac{2k+1}{2k-1}\frac{(k-1)^2}{(k+1)^2}=1-\frac{6k^2-2}{2k^3+3k^2-1}\approx \frac{k}{k+3}.
\end{split}
\end{equation}
These are well-known coefficients for Nesterov's method \eqref{eqn:NAtwo} with a constant step size, which corresponds to \eqref{eqn:sec}. 
We discuss the differences between the trajectories of \eqref{eqn:ODE-C} for $h\to0$ and $h=1$ in Section \ref{sec:exp}. 

\subsection{ODE-SC and the ODE model 
of~\citet{wilson2021lyapunov}
}\label{sec:comp:wilson}

Now, suppose that $f$ is $\mu$-strongly convex. 
We consider \eqref{eqn:ODE-SC} with the following coefficient function:
\[
\A(t)=e^{\sqrt{\frac{\mu}{L}}t}.
\] 
If we choose
\[
A_k=e^{\sqrt{\frac{\mu}{L}}hk},
\]
then the three conditions on $\A$ in Proposition \ref{thm:Nes-Flow} are satisfied. As in the previous subsection, 
it follows from \eqref{eqn:s} that we can deduce the coefficients for (NAG-SC) as
\begin{equation}\label{eqn:coeff_sc}
\begin{split}
a_k&=\frac{A_{k+1}-A_k}{2A_{k+1}-A_k}=\frac{e^{\sqrt{\frac{\mu}{L}}h}-1}{2e^{\sqrt{\frac{\mu}{L}}h}-1}\\
s_k &= \frac{(A_{k+1}-A_k)^2}{\mu A_{k+1}^2} = \frac{1}{\mu}(1-e^{-\sqrt{\frac{\mu}{L}}h})^2 \leq \frac{h^2}{L}.
\end{split}
\end{equation}
The inequality holds because $1 - e^{-x} \leq x$ for $x \geq 0$. 
Again, when $h \leq 1$, we have $s_k \leq 1/L$ as desired. 
Setting
\[
\ac (t)= \frac{e^{\sqrt{\frac{\mu}{L}}h}-1}{2e^{\sqrt{\frac{\mu}{L}}h}-1}
\]
yields $\ac (hk) = a_k$.

The ODE model~\eqref{eqn:ODE-SC} with this choice of $\A (t)$ and $\ac (t)$ can be written as
\begin{equation}\label{eqn:ODE-SC-h}
    \ddot{X}+(2-\ac(t))\sqrt{\frac{\mu}{L}}\dot{X}+\frac{1}{L}\nabla f\bigg (X+\ac(t)\sqrt{\frac{L}{\mu}}\dot{X}\bigg )=0.
\end{equation}
Letting $h\to 0$, the second-order ODE is further simplified to 
\begin{equation}\label{eqn:wilson_ode}
    \ddot{X}+2\sqrt{\frac{\mu}{L}}\dot{X}+\frac{1}{L}\nabla f(X)=0.
\end{equation}
Using the change of variables $\tilde{X}(t) = X(\sqrt{L}t)$, it can be verified that this ODE model is equivalent to the one proposed by 
\citet{wilson2021lyapunov}; see Equation (7) in \citet{wilson2021lyapunov}.

Similar to the previous subsection, we consider the case when $h = 1$, setting aside the relationship between our models and the model 
of \citet{wilson2021lyapunov}.
When $h = 1$, for a sufficiently large  $\kappa = L /\mu$, the coefficients $s_k$ and $b_k$ used in \eqref{eqn:NAtwo} can be approximated as
\begin{equation}\label{eq:sbsb}
\begin{split}
    s_k&=\frac{1}{\mu}(1-e^{-\sqrt{\frac{\mu}{L}}})^2\approx\frac{1}{L}\\
    b_k  &= \frac{e^{-\sqrt{\frac{\mu}{L}}}}{2-e^{-\sqrt{\frac{\mu}{L}}}} \simeq \frac{\sqrt{\kappa}-1}{\sqrt{\kappa}+1}.
\end{split}
\end{equation}
Again, these are the well-known coefficients used in Nesterov's method for strongly convex objective functions~\citep{nesterov2003introductory}.

\subsection{G-ODE and the Bregman Lagrangian flow 
of~\citet{wilson2021lyapunov}}

When $\ac (t) = 0$, \eqref{eqn:G-ODE-C} and \eqref{eqn:G-ODE-UC} are given by 
\begin{equation}\label{eqn:BLF-C}
\begin{split}
Z&=X+e^{-\alpha}\dot{X} \\
        \frac{d}{dt}\nabla g (Z) &= -e^{\alpha+\beta}\nabla f(X)
\end{split}
\end{equation}
and
\begin{equation}\label{eqn:BLF-SC}
\begin{split}
Z&=X+e^{-\alpha}\dot{X} \\
        \frac{d}{dt}\nabla g(Z)&=-\dot{\beta}(\nabla g(Z)-\nabla g(X)) - \frac{e^{\alpha}}{\mu}\nabla f(X),
\end{split}
\end{equation}
respectively.
These two systems of ODEs coincide with Equations (4) and (6) in 
\citet{wilson2021lyapunov},
  indicating that our generalized ODE model covers
the BLF. 
It is important to note that our model is  derived independently and is not related to Bregman Lagrangians. 

\citet{wilson2021lyapunov} obtained the scheme~\eqref{eqn:NA2} using BLF. 
However, their discretization process suddenly introduces a new variable $y_k$, which does not exist in BLF; see Section~\ref{sec:blf-gap}. 
Moreover, their Lyapunov function is described for $f(y_k)$ instead of $f(x_k)$. By contrast, our G-ODE models yield a more natural rate-matching discretization because $Y(t)$ is already defined in both \eqref{eqn:G-ODE-C} and \eqref{eqn:G-ODE-UC}. Therefore, G-ODE describes Nesterov's method and its accelerated convergence rate naturally without artificially introducing a variable or modifying the discrete-time scheme.

\subsection{G-ODE and the dynamical system model 
of~\citet{muehlebach2019dynamical}}

Consider the generalized ODE model with $g(x)=\frac{1}{2}\|x\|^2$.
When $f$ is  convex, 
let
\begin{equation*} 
e^{\alpha(t)}=\frac{2}{t+3},\quad e^{\beta(t)}=\frac{(t+3)^2}{4L},\quad \ac(t) =\frac{2t}{(t+3)^2}.
\end{equation*}
In this case, \eqref{eqn:G-ODE-C} can be expressed as
\begin{equation}\nonumber
\ddot{X}+\frac{3}{t+3}\dot{X}+\frac{1}{L}\nabla f \bigg (X+\frac{t}{t+3}\dot{X} \bigg )=0,
\end{equation}
which is equivalent to the dynamical system model proposed 
by~\citet{muehlebach2019dynamical}; see Equation (1) with coefficients defined as (14) in \citet{muehlebach2019dynamical}.

When $f$ is $\mu$-strongly convex, 
we consider the following choice of coefficients:
\begin{equation*} 
e^{\alpha}=\frac{1}{\sqrt{\kappa}},\quad \dot{\beta}=\frac{\sqrt{\kappa}-1}{\kappa+1},\quad \ac =\frac{\sqrt{\kappa}-1}{\sqrt{\kappa}(\sqrt{\kappa}+1)}.
\end{equation*}
Then, \eqref{eqn:G-ODE-UC} is given by
\begin{equation*}
\ddot{X}+\frac{2}{\sqrt{\kappa}+1}\dot{X}+\frac{1}{L}\nabla f \bigg (X+\frac{\sqrt{\kappa}-1}{\sqrt{\kappa}+1}\dot{X} \bigg )=0,
\end{equation*}
which is the same as the ODE model proposed 
by~\citet{muehlebach2019dynamical}; see Equation (1) with coefficients defined as (2) in \citet{muehlebach2019dynamical}.

\subsection{G-ODE and the high-resolution 
ODEs 
of~\citet{shi2021understanding,chen2022revisiting}
}\label{sec:shi}

The high-resolution ODEs proposed by \citet{shi2021understanding} and \citet{chen2022revisiting}
can also be obtained using our models. We first consider the high-resolution ODE for a convex function $f$ \citep{shi2021understanding}, 
\begin{equation}\label{eqn:shi-cv}
    \ddot{X} + \frac{3}{t}\dot{X}+\sqrt{s}\nabla^2 f(X)\dot{X}+(1+\frac{3\sqrt{s}}{2t})\nabla f(X)=0
\end{equation}
with any step size $0\le s\leq 1/L$. If $t$ is sufficiently large, the following approximations hold: 
\begin{equation}\label{eq:shiapp1} 
    \sqrt{s}\nabla^2 f(X)\dot{X}+(1+\frac{3\sqrt{s}}{2t})\nabla f(X) \approx \sqrt{s}\nabla^2 f(X)\dot{X}+\nabla f(X) \approx \nabla f(X+\sqrt{s}\dot{X}).
\end{equation}
Thus, we can rewrite \eqref{eqn:shi-cv} as  
\begin{equation*}
    \ddot{X}+\frac{3}{t}\dot{X}+\nabla f(X+\sqrt{s}\dot{X})=0,
\end{equation*}
which belongs to \eqref{eqn:G-ODE-C}.
Meanwhile, the high-resolution ODE for a $\mu$-strongly convex function $f$ is written as 
\begin{equation}\label{eqn:shi-sc}
\ddot{X}+2\sqrt{\mu}\dot{X}+\sqrt{s}\nabla^2 f(X)\dot{X}+(1+\sqrt{\mu s})\nabla f(X)=0
\end{equation}
with any step size $0\le s\leq 1/L$. 
If $\mu \ll L$, we observe that 
\begin{equation}\label{eq:shiapp2}
    \sqrt{s}\nabla^2 f(X)\dot{X}+(1+\sqrt{\mu s})\nabla f(X) \approx \sqrt{s}\nabla^2 f(X)\dot{X}+  \nabla f(X) \approx \nabla f(X+\sqrt{s}\dot{X}).
\end{equation}
Thus, \eqref{eqn:shi-sc} is rewritten as 
\begin{equation*}
    \ddot{X}+2\sqrt{\mu}\dot{X}+\nabla f(X+\sqrt{s}\dot{X})=0
\end{equation*} 
in a form that fits within \eqref{eqn:G-ODE-UC}. 
The approximations we use in \eqref{eq:shiapp1} and \eqref{eq:shiapp2} for the models  
of \citet{shi2021understanding}
are empirically validated in Sections \ref{sec:exp1} and \ref{sec:6.2}, respectively. 

Furthermore, we consider the high-resolution ODE with varying step sizes \citep{chen2022revisiting},  
\begin{equation}\label{eqn:fofosks}
    \ddot{X}(t)+2\sqrt{\mu} \dot{X}(t)+\nabla f(X+\frac{\sqrt{s}}{1+2\sqrt{\mu s}}\dot{X}(t))=0.
\end{equation}
By letting $e^{\alpha}=\sqrt{\mu}$, $\beta(t)=\sqrt{\mu}t$, $g(x)=\frac{1}{2}\|x\|^2$, and $\ac(t)=\frac{\sqrt{\mu s}}{1+2\sqrt{\mu s}}$, we observe that \eqref{eqn:G-ODE-UC} includes \eqref{eqn:fofosks}.

\section{Restart Schemes for Nesterov's Methods}\label{sec:restart}
In this section, as our third main contribution, we demonstrate that restart schemes for Nesterov's methods can be designed using our generalized models. First, we  rewrite both \eqref{eqn:G-ODE-C} and \eqref{eqn:G-ODE-UC} in the unified form:
\begin{equation*}
\begin{split}
    &\ddot{X}+c_1(t)\dot{X}+c_2(t)\nabla f(Y) = 0,\\
    &Y=X+b(t)\dot{X}.
    \end{split}
\end{equation*}
If $c_1, c_2 ,b > 0$, as $f$ is convex, it follows that 
\begin{equation*}
     b(t) \langle \dot{X},\nabla f(Y)-f(X) \rangle = \langle Y-X, \nabla f(Y)-f(X) \rangle \geq 0.
\end{equation*}
We define the restarting time $T$ as 
\begin{equation*}
    T = T(x_0;f) = \sup \Bigl\{ t>0, \forall u\in(0,t), \frac{d\|\dot{X}(u)\|^2}{du}\geq0 \Bigr\}.
\end{equation*}
For $t\leq T$, $f\left(X\right)$ continues to decrease because 
\begin{equation}\nonumber
    \frac{d}{dt}f(X) = \langle \dot{X},\nabla f(X)\rangle \leq \langle \dot{X},\nabla f(Y)\rangle = -\frac{1}{2c_2}\frac{d}{dt}\|\dot{X}\|^2 - \frac{c_1}{c_2}\|\dot{X}\|^2 <0. 
\end{equation} 
This suggests that $f(X(t))$  decreases monotonically with $t$ 
whenever 
\begin{equation}\label{eqn:velpos}
\langle \ddot{X}(t), \dot{X}(t) \rangle \geq 0.
\end{equation} 
Through continuous-time analysis, we observe that 
\begin{equation*}
    \frac{d\|\dot{X}(t) \|^2}{dt} = \langle \ddot{X}(t), \dot{X}(t) \rangle \approx \frac{1}{h^3}\langle x_{k+1}-2x_{k}+x_{k-1}, x_{k}-x_{k-1} \rangle.
\end{equation*}
Therefore, a discrete-time counterpart of the condition \eqref{eqn:velpos} can be written as 
\begin{equation}\label{eq:ourcondition}
\langle x_{k+1}-2x_k+x_{k-1},x_k-x_{k-1} \rangle \geq 0 \quad \forall k,
\end{equation}
under which $f(x_k)$ decreases with $k$.  We formalize this idea as follows: 
\begin{theorem}\label{thm:restart}
Consider Nesterov's method \eqref{eqn:NAtwo} with $0< s_k\leq \frac{1}{L}$ and $0< b_k\leq 1$ for all $k$. 
If 
 \eqref{eq:ourcondition} holds, 
then  $f(x_{k+1})<f(x_k)$ for all $k$.
\end{theorem}

\begin{proof}
Given that $\nabla f$ is Lipschitz continuous and $x_{k+1}=y_k-s_k\nabla f(y_k)$, it follows that 
\begin{equation*}
    f(x_{k+1})\leq f(y_k)-s_k(1-\frac{L s_k}{2})\|\nabla f(x_k)\|^2 < f(y_k).
\end{equation*}
Furthermore, as $f$ is convex, it holds that 
\begin{align*}
    f(y_k)-f(x_k) &\leq \langle\nabla f(y_k), y_k-x_k \rangle \\
    &= -\frac{b_k}{s_k}(\langle x_{k+1}-x_k,x_k-x_{k-1}\rangle-b_k\|x_k-x_{k-1}\|^2) \\
    &\leq -\frac{b_k}{s_k}\langle x_{k+1} - 2x_k+ x_{k-1}, x_k-x_{k-1}\rangle \leq 0.
\end{align*}
This implies that $f(x_{k+1})< f(x_k)$.
\end{proof}

Note that the condition $b_k\leq1$ in Theorem \ref{thm:restart} is met by the frequently used choices $b_k=\frac{k}{k+3}$ and  $b_k=\frac{\sqrt{L}-\sqrt{\mu}}{\sqrt{L}+\sqrt{\mu}}$. 

Leveraging Theorem \ref{thm:restart} and the general observation that convergence accelerates when $k$ is small, we present a restart scheme for Nesterov's method \eqref{eqn:NAtwo} in Algorithm \ref{alg:cap}. 
It is worth comparing Algorithm \ref{alg:cap} to the original restart scheme presented 
by \citet{su2016differential}. 
Both schemes are motivated by the fact that $f(X(t))$ is a decreasing function for $t\le T$. 
However, unlike the restart scheme 
of \citet{su2016differential}
that initializes Nesterov's method whenever 
\begin{equation*}
\|x_{k+1}-x_k\| > \|x_k-x_{k-1}\| 
\end{equation*} 
is satisfied, we use the condition \eqref{eq:ourcondition}. 
Simultaneously, we let $x_{k+1} = x_k - s_k\nabla f(x_k)$ when Nesterov's method restarts, i.e., when \eqref{eq:ourcondition} is violated.
As a result, we can theoretically guarantee that $f(x_k)$ monotonically decreases in discrete time. 
Note that \citet{su2016differential} did not provide a proof of the monotonicity in the sequence of the objective value obtained using their restart scheme. In fact, there exists an empirical evidence that the restart scheme 
of \citet{su2016differential}
generally fails to achieve a monotonic decrease in the objective value, which we present in Section \ref{sec:restartexp}. 
Moreover, Algorithm \ref{alg:cap} can be applied to any Nesterov's method in the form (\ref{eqn:NAtwo}), whereas the restart scheme 
of \citet{su2016differential}
covers only Nesterov's methods with constant step sizes.
Meanwhile, Algorithm \ref{alg:cap} is more effective in the non-strongly convex case than in the strongly convex case, because the coefficients of the ODEs remain almost unchanged with time when $f$ is strongly convex. 
Numerical examples that demonstrate the effectiveness of Algorithm \ref{alg:cap} are also given in Section \ref{sec:restartexp}

Meanwhile, a restart scheme is impractical in the strongly convex case, as the ODE coefficients   remain almost unchanged with time. 

\begin{algorithm}
\caption{Restart scheme for Nesterov's method \eqref{eqn:NAtwo}}
\label{alg:cap}
\begin{algorithmic}
\State $s_k, b_k$ in Nesterov's method \eqref{eqn:NAtwo}
\State $j \gets 1$
\For{$k=1$ to $k_{max}$}
\State $y_k \gets x_k+b_k(x_k-x_{k-1})$
\State $x_{k+1} \gets y_k-s_k\nabla f(y_k)$
\If{
$\langle x_{k+1}-2x_k+x_{k-1},x_k-x_{k-1} \rangle < 0$
and $j\geq k_{min}$}
    \State $j \gets 1$
    \State $x_{k+1} \gets x_k-s_k\nabla f(x_k)$
\Else
    \State $j \gets j+1$
\EndIf
\EndFor
\end{algorithmic}
\end{algorithm}

\section{Connections to Gradient Flow: Reparametrization of Time}
\label{sec:time}

In this section, we examine our generalized ODE models through the lens of gradient flows. 
Specifically, as the fourth main contribution of this study, we demonstrate that the accelerated convergence rate for our generalized ODE models can be attributed to the reparametrization of time in the standard gradient flow when $\ac (t) = 1$. If $\ac (t) = 1$, we have $Y=Z$. Thus, \eqref{eqn:G-ODE-C} and \eqref{eqn:G-ODE-UC} can be reformulated in terms of a single variable $Z$ as follows: 
\begin{equation} \label{eqn:C-a1}
        \frac{d}{dt}\nabla g (Z) = -e^{\alpha(t)+\beta(t)}\nabla f(Z)
    \end{equation}
    and
    \begin{equation} \label{eqn:SC-a1}
       \frac{d}{dt}\nabla g(Z)=- \frac{e^{\alpha(t)}}{\mu}\nabla f(Z),
    \end{equation}
    respectively.
    We establish that $Z$ has the same convergence rate of $O(e^{-\beta(t)})$
     using the following Lyapunov function:
 \begin{equation}\nonumber
\tilde{V}(Z, t) = \left \{
\begin{array}{ll}
D_g (0, Z) + e^{\beta (t)} f(Z) &\mbox{if $\mu = 0$}\\
e^{\beta (t)} (\mu D_g (0, Z) + f(Z) ) &\mbox{if $\mu > 0$}.
\end{array}
\right.
 \end{equation}   
Moreover, $Z$ possesses a descent property that $X$ lacks. 
    
 \begin{theorem}\label{thm:conv2}
 Let $Z(t) :[0, \infty) \to \mathbb{R}^n$ be a solution to \eqref{eqn:C-a1} when $f$ is  convex and a solution to \eqref{eqn:SC-a1} when $f$ is $\mu$-uniformly convex with respect to $g$. 
 Then, the energy function $\tilde{\mathcal{E}} (t) = \tilde{V}(Z(t), t)$ is monotonically non-increasing for both convex and uniformly convex cases, and
\[
f(Z(t)) \leq e^{-\beta(t)} \tilde{\mathcal{E}} (0) \quad \forall t \geq 0.
\]
\end{theorem}   

The proof of Theorem~\ref{thm:conv2} is provided in Appendix~\ref{app:conv2}.
The ODE models~\eqref{eqn:C-a1} and \eqref{eqn:SC-a1} can be unified into a single model using a new continuously differentiable coefficient function $\tau (t): [0,\infty) \to [0,\infty)$ such that 
\begin{equation}\label{eqn:condition2}
\tau(0) = 0, \quad \dot{\tau} > 0, \quad \tau(t) \to \infty \mbox{ as $t \to \infty$}.
\end{equation}

\begin{corollary}\label{cor:Z}
Consider the following ODE:
\begin{equation}\label{eqn:QGF0}
\frac{d}{dt}\nabla g(Z)=-\dot{\tau}\nabla f(Z)
\end{equation}
with $Z(0) = x_0$.
Then, when $f$ is  convex, 
this ODE with $\tau(t)= \int_{0}^{t}e^{\alpha(t')+\beta(t')}dt'$
is equivalent to \eqref{eqn:C-a1} and achieves the following convergence rate:
\[
f(Z(t)) \leq O\bigg (\frac{1}{\tau(t)} \bigg ).
\]
Moreover, when $f$ is $\mu$-uniformly convex with respect to $g$, 
the ODE with $\tau(t)= \int_{0}^{t}\frac{e^{\alpha(t')}}{\mu}dt'$ is equivalent to \eqref{eqn:SC-a1} and has the following convergence rate:
\[
f(Z(t)) \leq O (e^{-\mu \tau (t)}).
\]
\end{corollary}

The proof of Corollary~\ref{cor:Z} is provided in Appendix~\ref{app:Z}.

\subsection{Quasi-gradient flow}
Note that the coefficient function $\tau (t)$ has its inverse, denoted by $\tau^{-1}$, for $t \geq 0$. 
Reparametrizing time via $\tau$, we define a new flow as
\begin{equation}\label{eqn:repara}
Q(t) = Z(\tau^{-1} (t)), \quad t \geq 0.
\end{equation}
The ODE~\eqref{eqn:QGF0} can then be expressed in terms of $Q$ as follows:
\begin{equation}\label{eqn:QGF} \tag{QGF}
\frac{d}{dt}\nabla g(Q)=-\nabla f(Q),
\end{equation}
which we call the \emph{quasi-gradient flow} (QGF).
Applying the forward Euler method yields the discrete-time scheme
\begin{equation*}
\frac{q_{k+1}-q_k}{h}=-[\nabla^2 g(q_k)]^{-1} \nabla f(q_k).
\end{equation*} 
This can be interpreted as the quasi-Newton method when $[\nabla^2 g(q_k)]^{-1}$ is chosen to approximate the inverse Hessian of $f$ at $q_k$. 
In particular, when $g(x) = \frac{1}{2} \| x\|^2$, 
\eqref{eqn:QGF} is compactly written as
\begin{equation}\label{eqn:GF} \tag{GF}
\frac{d}{dt} Q=-\nabla f(Q),
\end{equation}
which corresponds to the standard gradient flow (GF). Furthermore, its discrete-time version is written as the following gradient descent scheme:
\begin{equation*}
\frac{q_{k+1}-q_k}{h}=-\nabla f(q_k).
\end{equation*}

\subsection{Understanding  acceleration as changes in the speed of time}

According to Corollary~\ref{cor:Z}, \eqref{eqn:QGF} and \eqref{eqn:GF} have the convergence rates of $O  (\frac{1}{t}  )$ and $O (e^{-\mu t})$, respectively. These are in line with  well-known results in the convergence analysis for gradient flow~\citep[Proposition 1.1]{scieur2017integration}. 
Notably, when $\tau(t) > t$, the generalized ODE achieves a faster convergence rate compared to the (quasi-)gradient flow. 

We now provide an intuitive explanation for the following question: 
\vspace{-0.05in}
\begin{center}
\emph{Why does the generalized ODE  converge faster than the gradient flow?}
\end{center}
\vspace{-0.05in}
Recall that the (quasi-)gradient flow is obtained by reparametrizing time in the generalized ODE via \eqref{eqn:repara}.
The coefficient function $\tau$ essentially dictates how fast the solution $Z$ to the generalized ODE evolves along the trajectory of \eqref{eqn:QGF}. 
Thus, the accelerated convergence in the generalized ODE can be interpreted as a consequence of \emph{changing the speed of time} in the gradient flow.\footnote{\citet{wibisono2016variational} showed that the BLF is closed under time dilation. \citet{attouch2018fast} proposed an ODE model that converges faster through a  reparametrization of time, noting that its implicit discretization leads to faster algorithms. However, neither study explicitly connected their models to gradient flow nor provided an intuitive interpretation of the acceleration.}

However, this interpretation comes with a caveat. The above discussion holds under the assumption that $\ac (t) = 1$. In reality, the values of $\ac(t)$ corresponding to the existing ODE models---whose relationships with our generalized models are analyzed in Section \ref{sec:ex}---seldom satisfy this condition. Thus, we cannot directly translate  this observation to discrete-time Nesterov's method. In summary, while we reveal that altering the speed of time induces  acceleration in continuous time, this does not fully explain why Nesterov's method outperforms gradient descent in discrete time. 
Nonetheless, we believe that  emphasizing the role of time reparametrization could pave the way for new research avenues in understanding Nesterov's method.

\section{Numerical Experiments}\label{sec:exp}
In this section, we provide the results of numerical experiments to substantiate our analyses in Sections \ref{sec:ex}--\ref{sec:time}.
Specifically, we compare the trajectories generated by the special cases \eqref{eqn:ODE-C} and \eqref{eqn:ODE-SC} of our generalized models against those produced by existing ODE models, as well as the iterates of Nesterov's methods. Moreover, we examine the efficacy of our restart scheme and explore the connections to gradient flow. 
The source code of our experiments is available online.\footnote{\url{https://github.com/CORE-SNU/G-ODE-NAG}}

\subsection{ODE-C and the ODE models 
of \citet{su2016differential, shi2021understanding}
}\label{sec:exp1}

In this subsection, we compare \eqref{eqn:ODE-C} with the models 
of \citet{su2016differential,shi2021understanding}.
We adopt two Nesterov's methods, namely (NAG-C), which is expressed by \eqref{eqn:Nes-CV-sb}, and (NAG-C-C), which is an approximation of (NAG-C) representing Nesterov's method with a constant step size. For both (NAG-C) and (NAG-C-C), our model and the two models 
of \citet{su2016differential,shi2021understanding}
are written as \eqref{eqn:ODE-CV-h}, \eqref{eqn:sec}, and \eqref{eqn:shi-cv}, respectively. We fix $L=s=h=1$, aiming for $x_k=X(t=k)$. We consider the objective function $f(x_1,x_2)=0.02 x_1^2 + 0.005 x_2^2$ with an initial point $x_0 = (1, 1)$. 
The initial conditions for the ODEs are given by $X(0)=(1,1)$ and $\dot{X}(0)=(0,0)$. 

\begin{figure}
\centering 
\begin{tabular}{cc}
\subfigure[Trajectories]{\includegraphics[width=0.43\textwidth ]{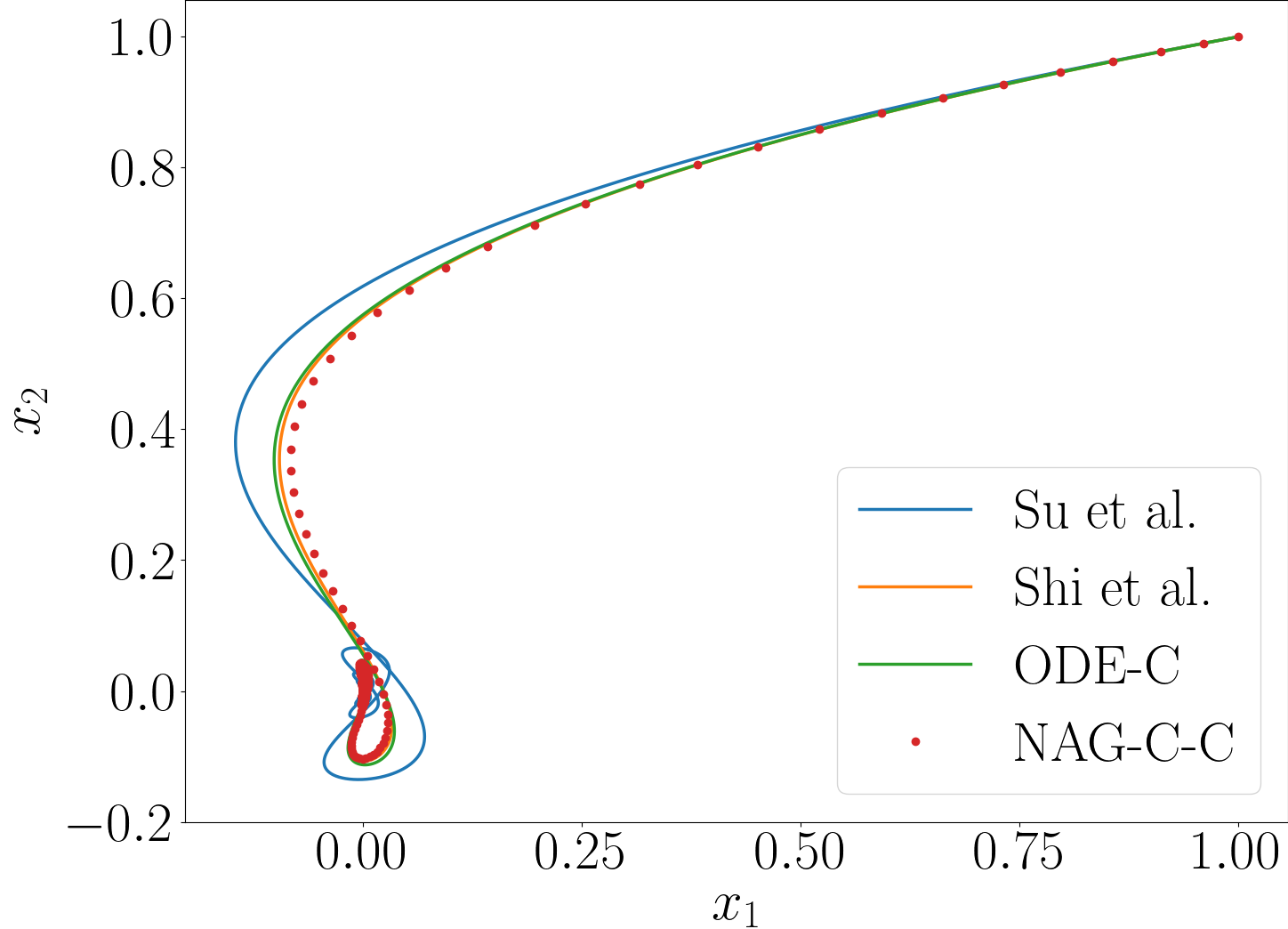}}& \subfigure[Zoomed trajectories]{\includegraphics[width=0.43\textwidth ]{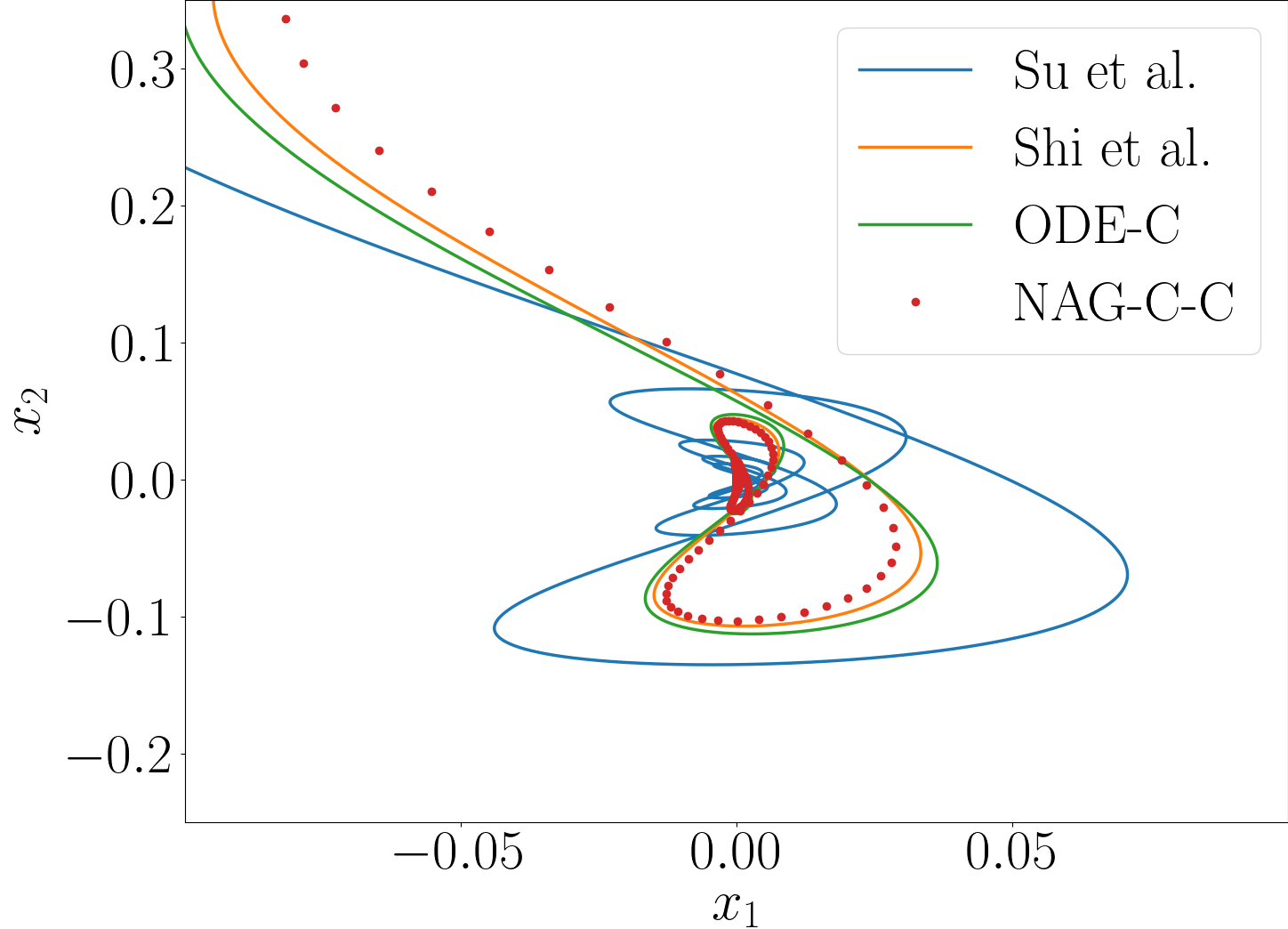}}\\ \subfigure[Errors $f - f^*$]{\includegraphics[width=0.43\textwidth ]{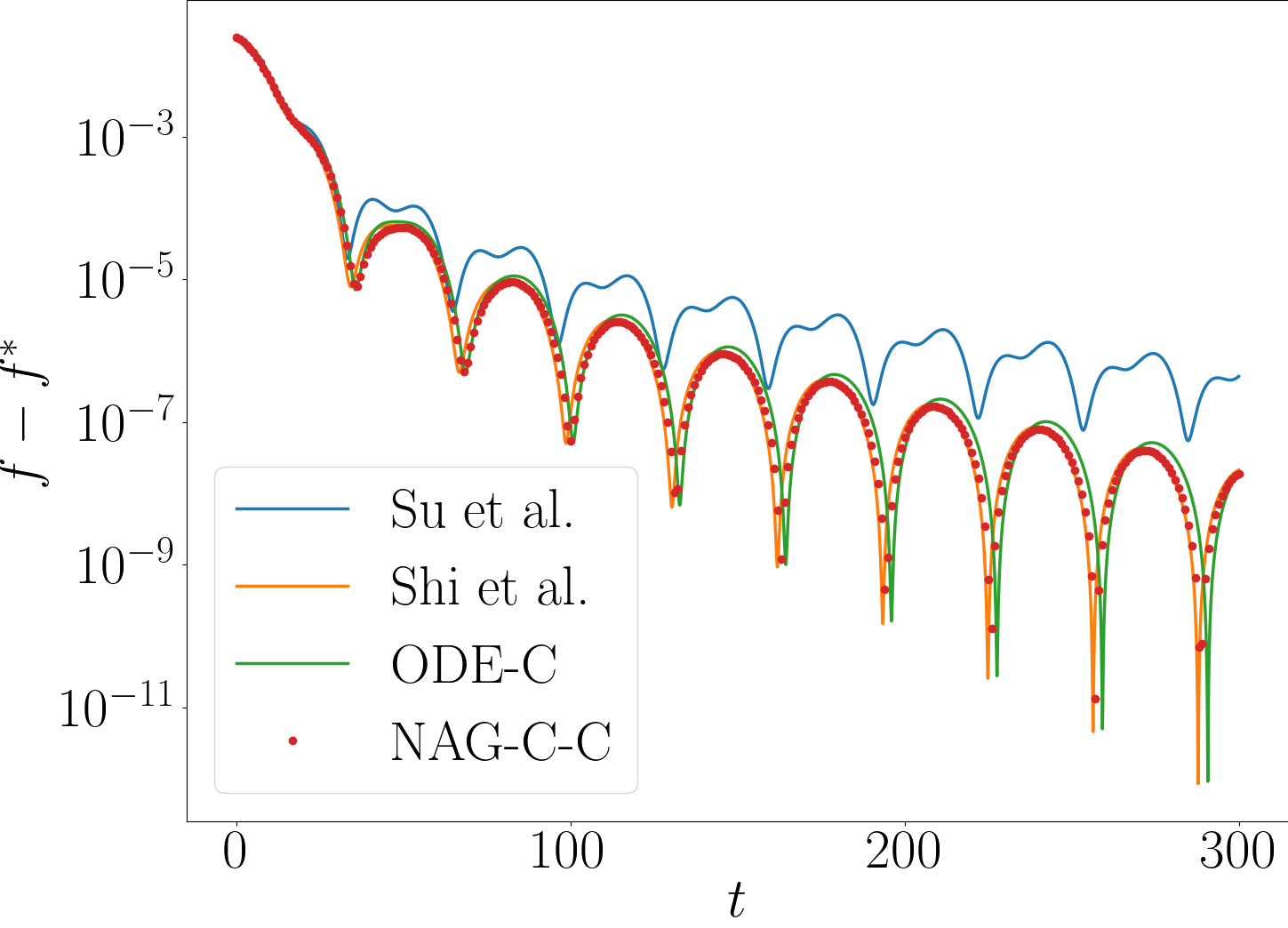}}& \subfigure[Errors $\|x_k-X(t=k)\|_2$]{\includegraphics[width=0.43\textwidth ]{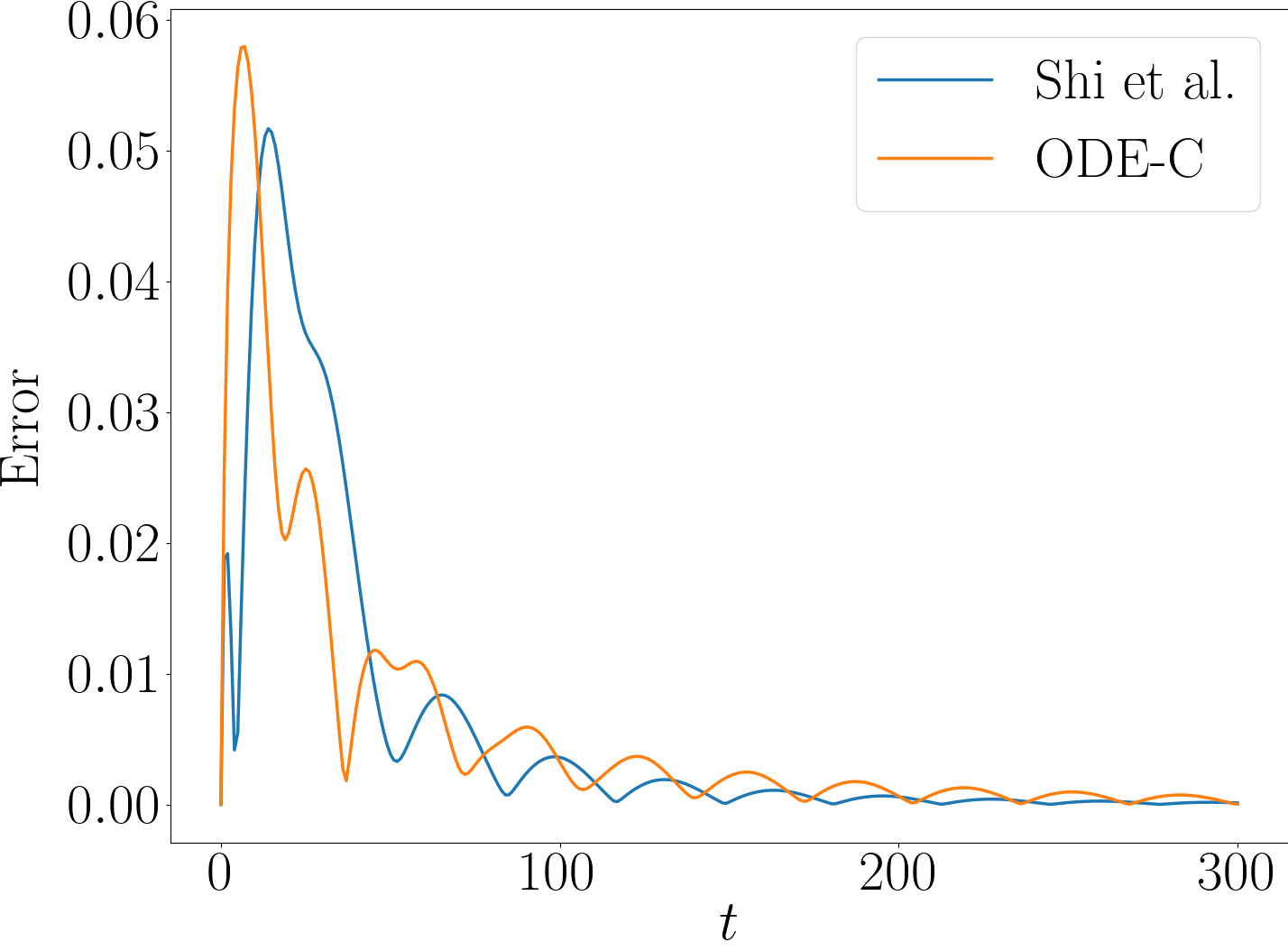}}\\
\end{tabular}
\caption{Comparison of \eqref{eqn:ODE-C} and the existing ODEs \citep{su2016differential, shi2021understanding} for (NAG-C-C).}\label{fig:exp1}
\end{figure}

\begin{figure}
\centering 
\begin{tabular}{cc}
\subfigure[Trajectories]{\includegraphics[width=0.43\textwidth ]{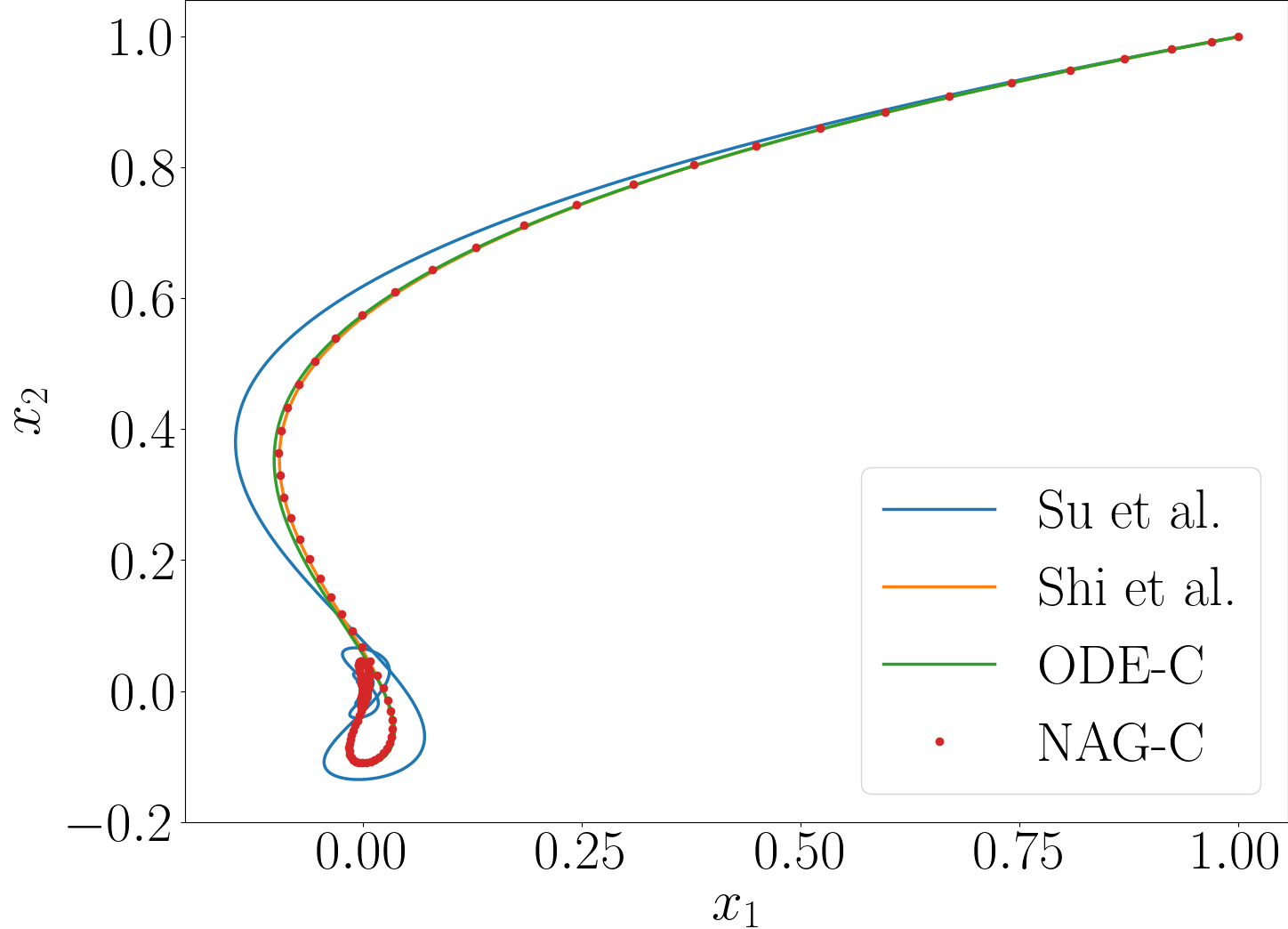}}& \subfigure[Zoomed trajectories]{\includegraphics[width=0.43\textwidth ]{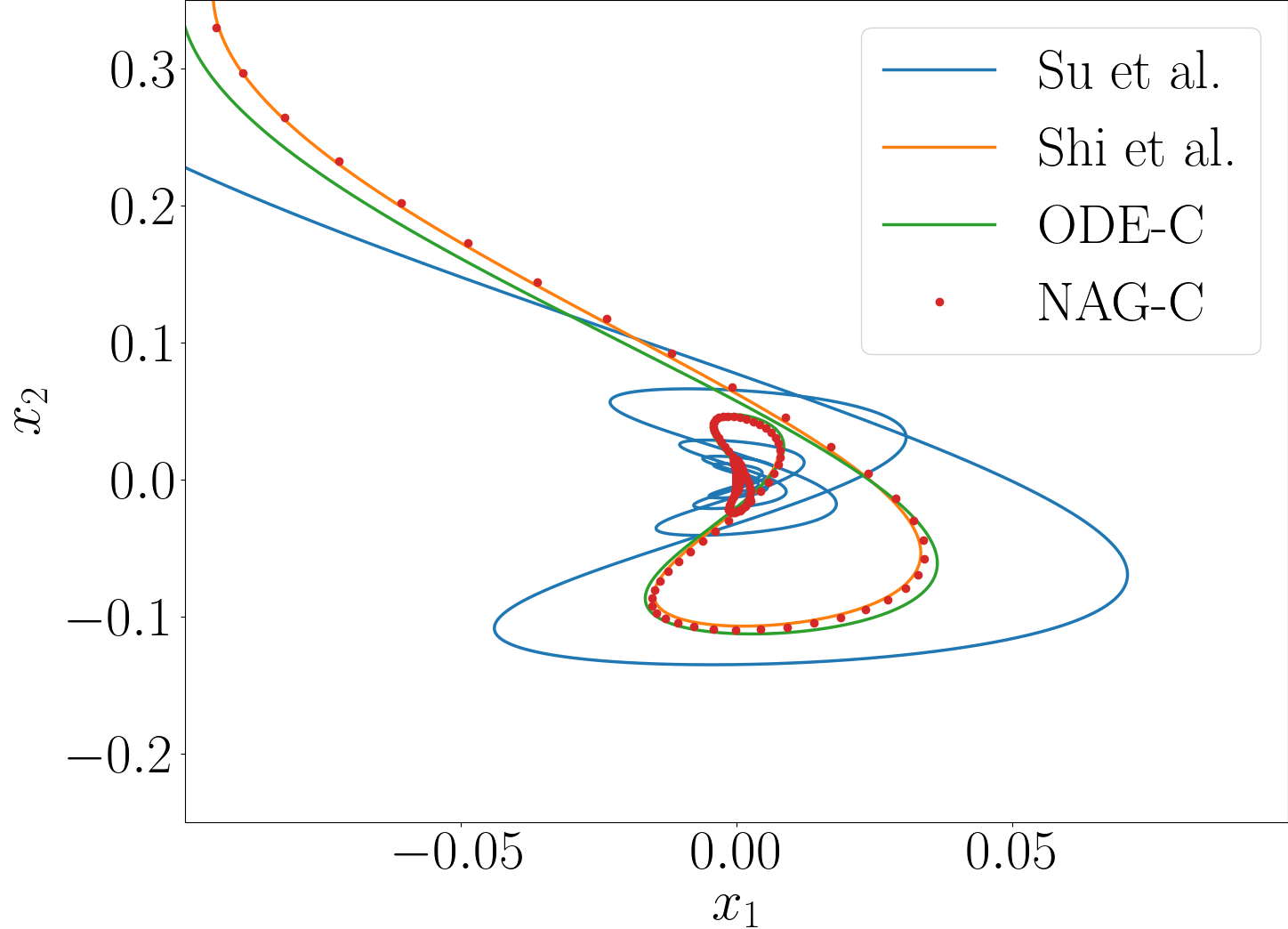}}\\ \subfigure[Errors $f - f^*$]{\includegraphics[width=0.43\textwidth ]{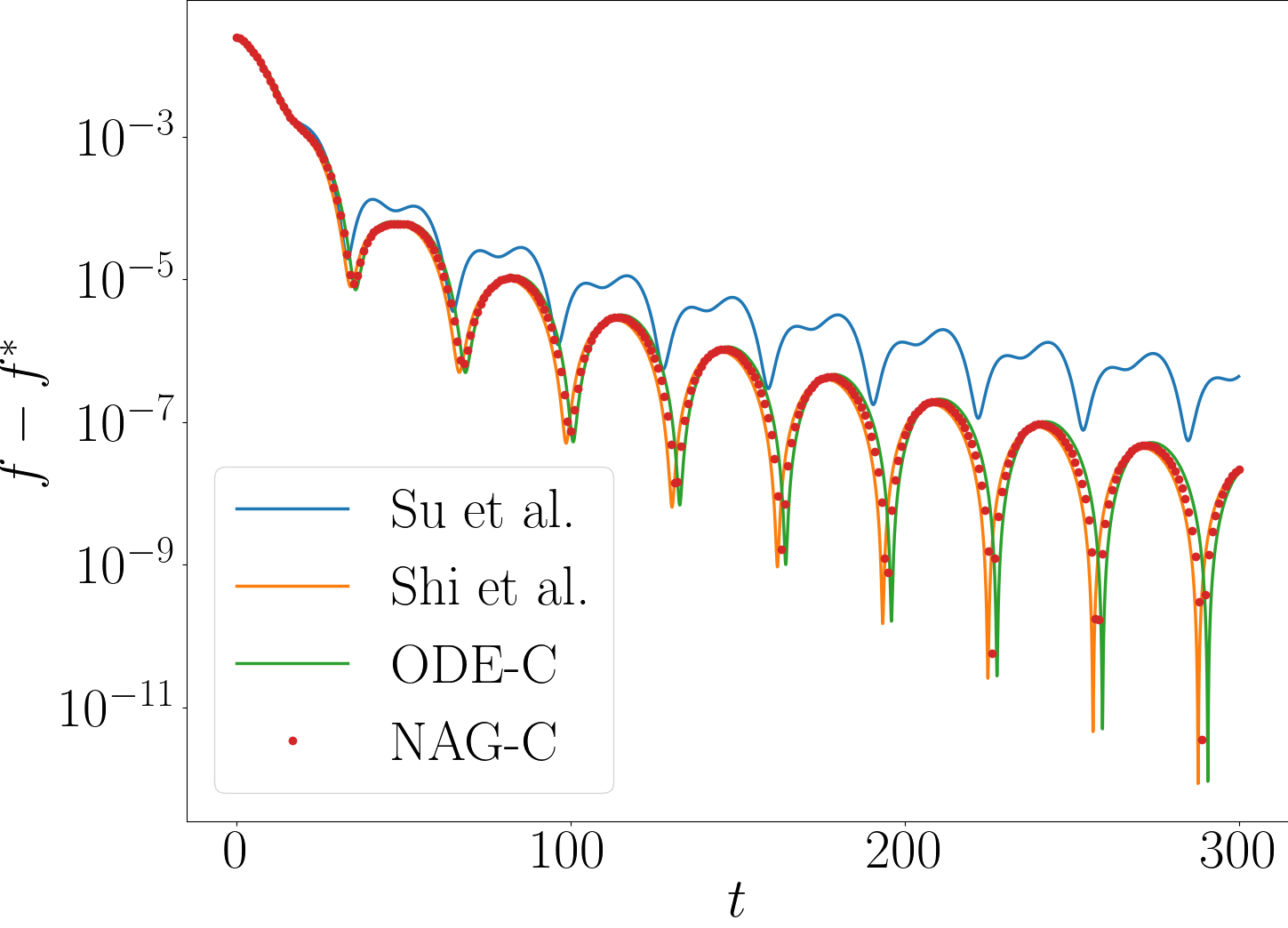}}& \subfigure[Errors $\|x_k-X(t=k)\|_2$]{\includegraphics[width=0.43\textwidth ]{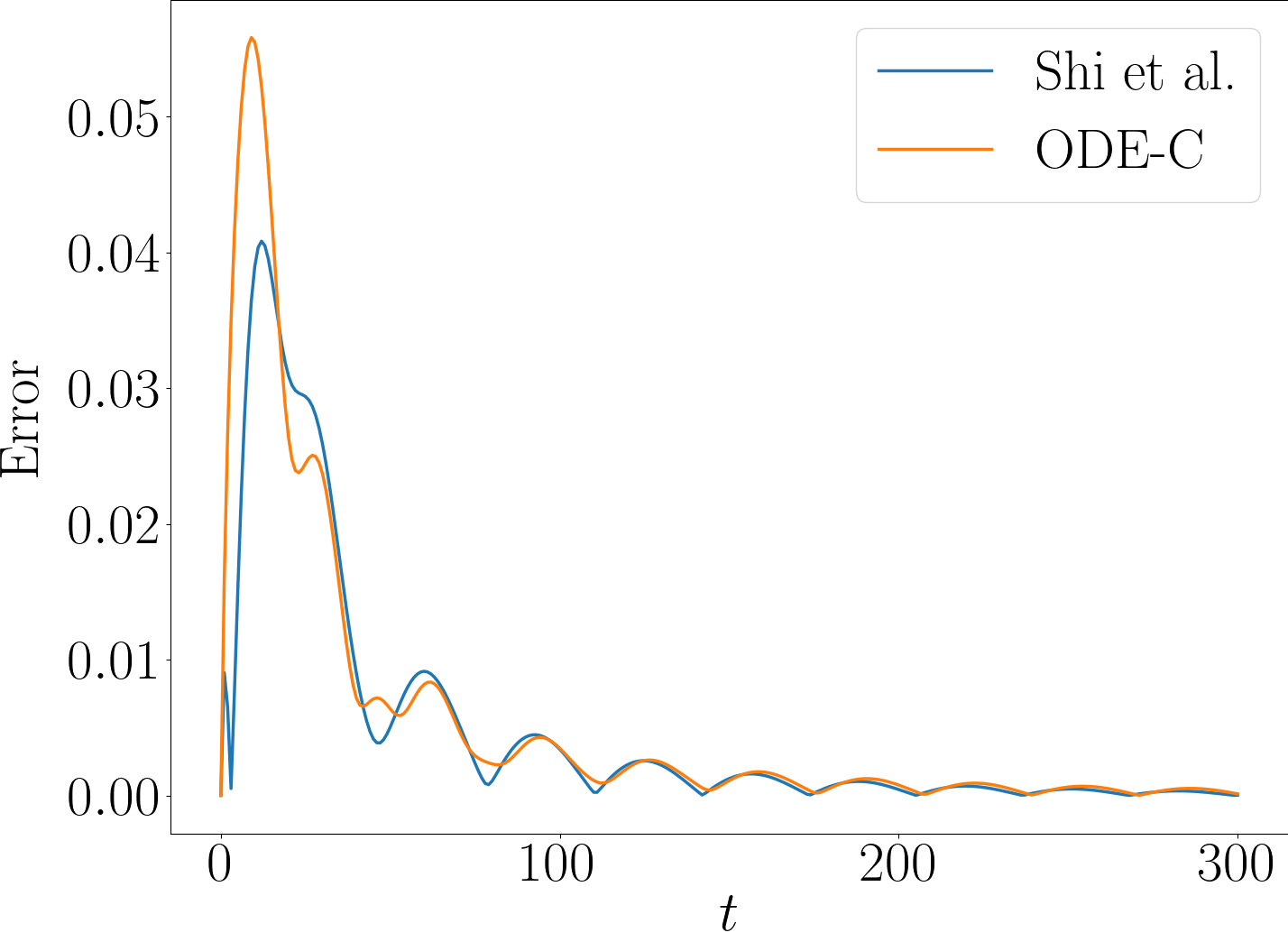}}\\
\end{tabular}
\caption{Comparison of \eqref{eqn:ODE-C} and the existing ODEs \citep{su2016differential, shi2021understanding} for (NAG-C).}\label{fig:exp2}
\end{figure}

The simulation results for (NAG-C-C) and (NAG-C) are illustrated in Figures \ref{fig:exp1} and \ref{fig:exp2}, respectively. Figure \ref{fig:exp1} shows that 
\eqref{eqn:ODE-C} describes (NAG-C-C) better than the ODE in \citet{su2016differential}. Specifically, for $100\le k\le 300$, \eqref{eqn:ODE-C} reduced the $L_2$-norm error $\Vert x_k - X(t=k)\Vert_2$ by 69.2\% on average compared to 
\citet{su2016differential}.
This improvement stems from  \eqref{eqn:ODE-C}'s ability to compute the gradient at a future time, as discussed in Section \ref{sec:gen}. In Figure \ref{fig:exp2}, \eqref{eqn:ODE-C} and the ODE 
of \citet{shi2021understanding}
have similar trajectories, which experimentally confirms the relationship studied in Section \ref{sec:shi}. 
The error between \eqref{eqn:ODE-C} and 
(NAG-C-C) shown in Figure \ref{fig:exp1} is due to the double approximation, compared to which the error between \eqref{eqn:ODE-C} and (NAG-C) shown in Figure \ref{fig:exp2} is smaller by 68.2\% in terms of the $L_2$-norm, averaged over $100\le k\le 300$. Note that Figure~\ref{fig:exp1}~(c) and Figure~\ref{fig:exp2}~(c) use a logarithmic scale, indicating that the error oscillation does not increase with time. 

\subsection{ODE-SC and the ODE models 
of \citet{wilson2021lyapunov,shi2021understanding}
}\label{sec:6.2}

In this subsection, we compare \eqref{eqn:ODE-SC} with the two ODE models proposed by \citet{wilson2021lyapunov} and \citet{shi2021understanding}. 
We consider the same objective function $f$ as in the previous subsection, treating  it as $\mu$-strongly convex. The Nesterov's methods employed are   (NAG-SC) 
with  coefficients as specified in \eqref{eqn:coeff_sc} and (NAG-SC-C), which is an approximation of (NAG-SC) with a constant step size. Our model and the two models 
of \citet{wilson2021lyapunov,shi2021understanding}
are represented by 
\eqref{eqn:ODE-SC-h}, \eqref{eqn:wilson_ode} and \eqref{eqn:shi-sc}, respectively. 
We let $\mu=0.001$, while keeping the other parameters consistent with  the previous subsection. 

\begin{figure}
\centering 
\begin{tabular}{cc}
\subfigure[Trajectories]{\includegraphics[width=0.43\textwidth ]{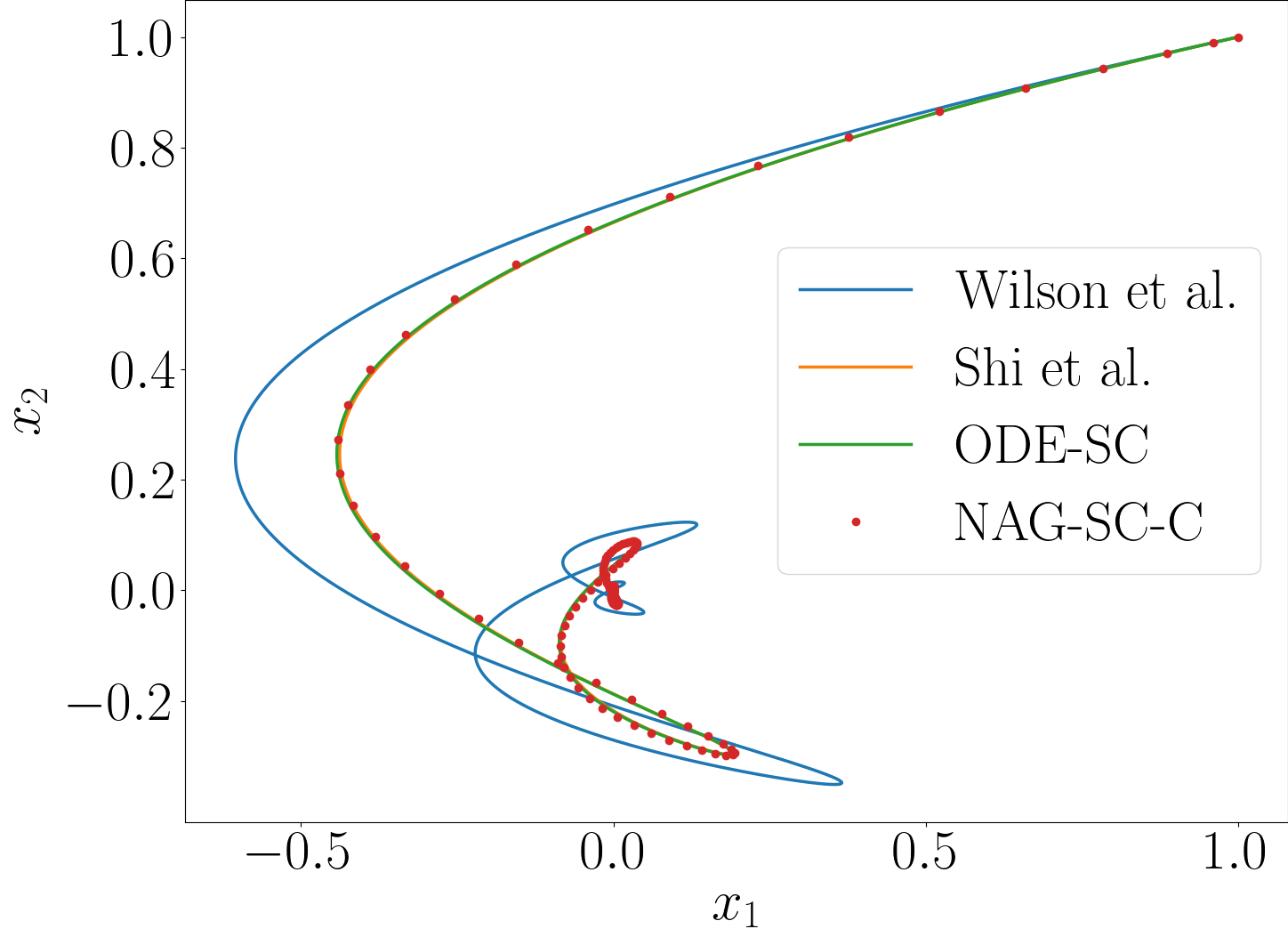}}& \subfigure[Zoomed trajectories]{\includegraphics[width=0.43\textwidth ]{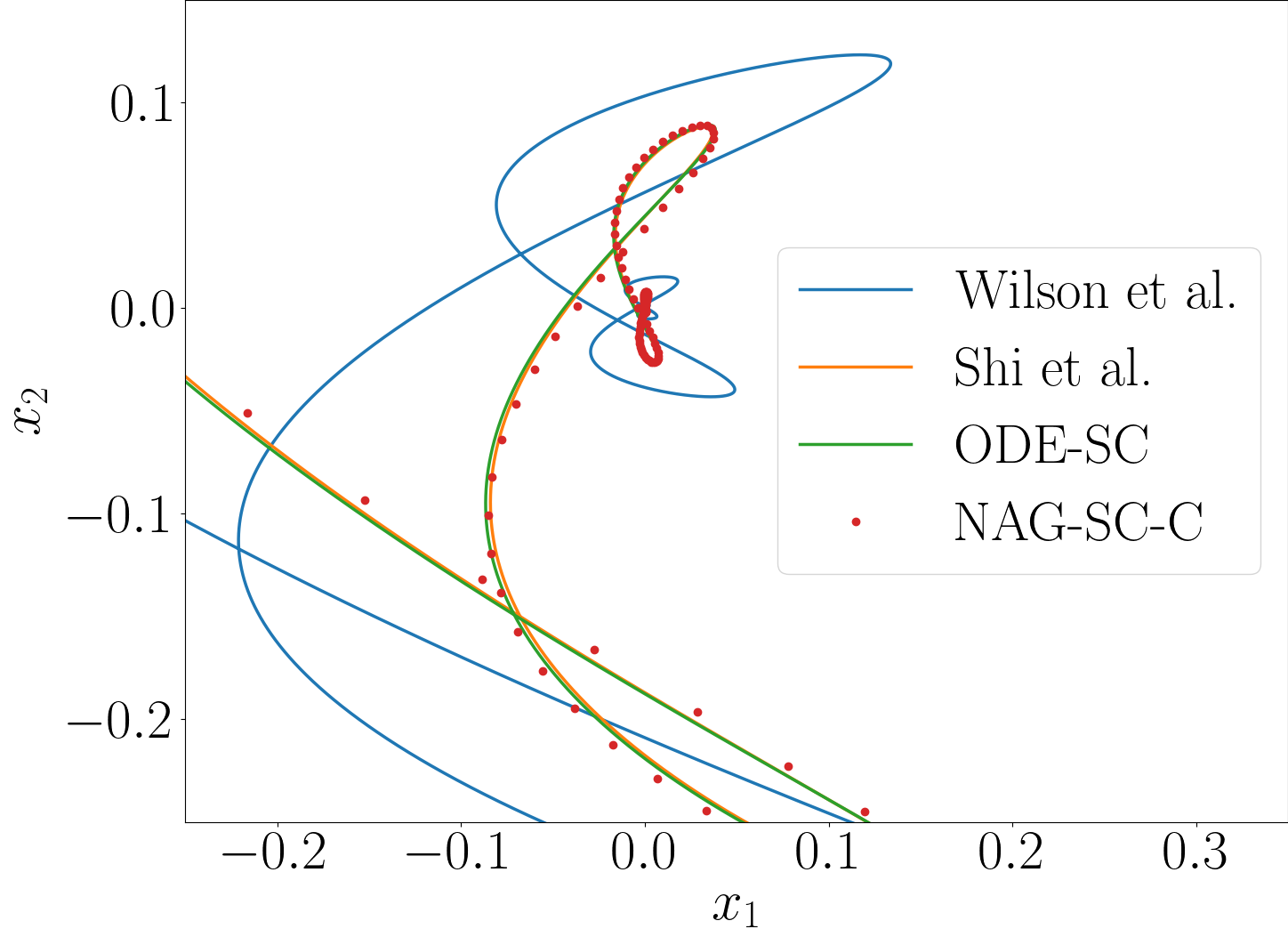}}\\ \subfigure[Errors $f - f^*$]{\includegraphics[width=0.43\textwidth ]{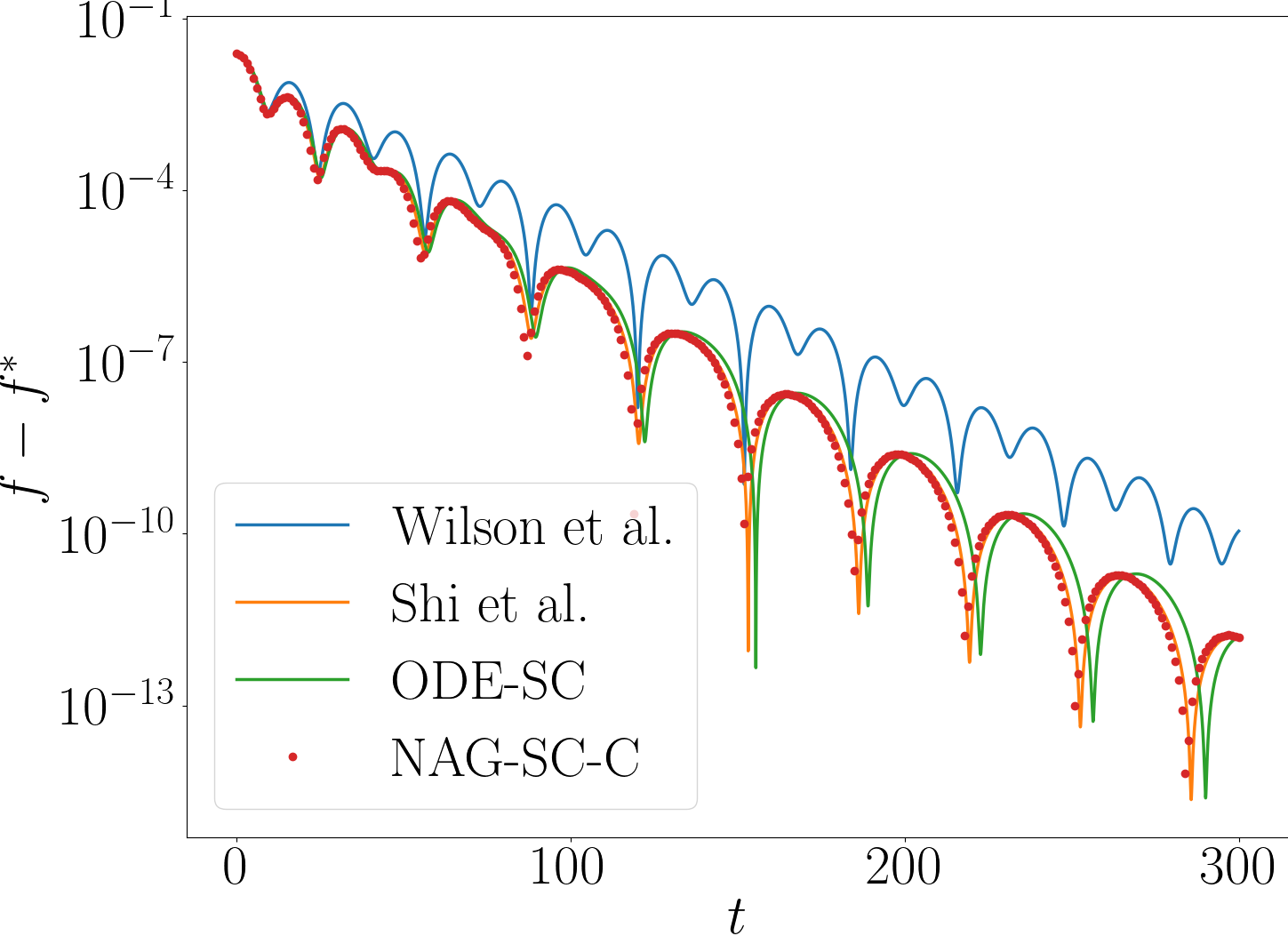}}& \subfigure[Errors $\|X(t=k)-x_k\|_2$]{\includegraphics[width=0.43\textwidth ]{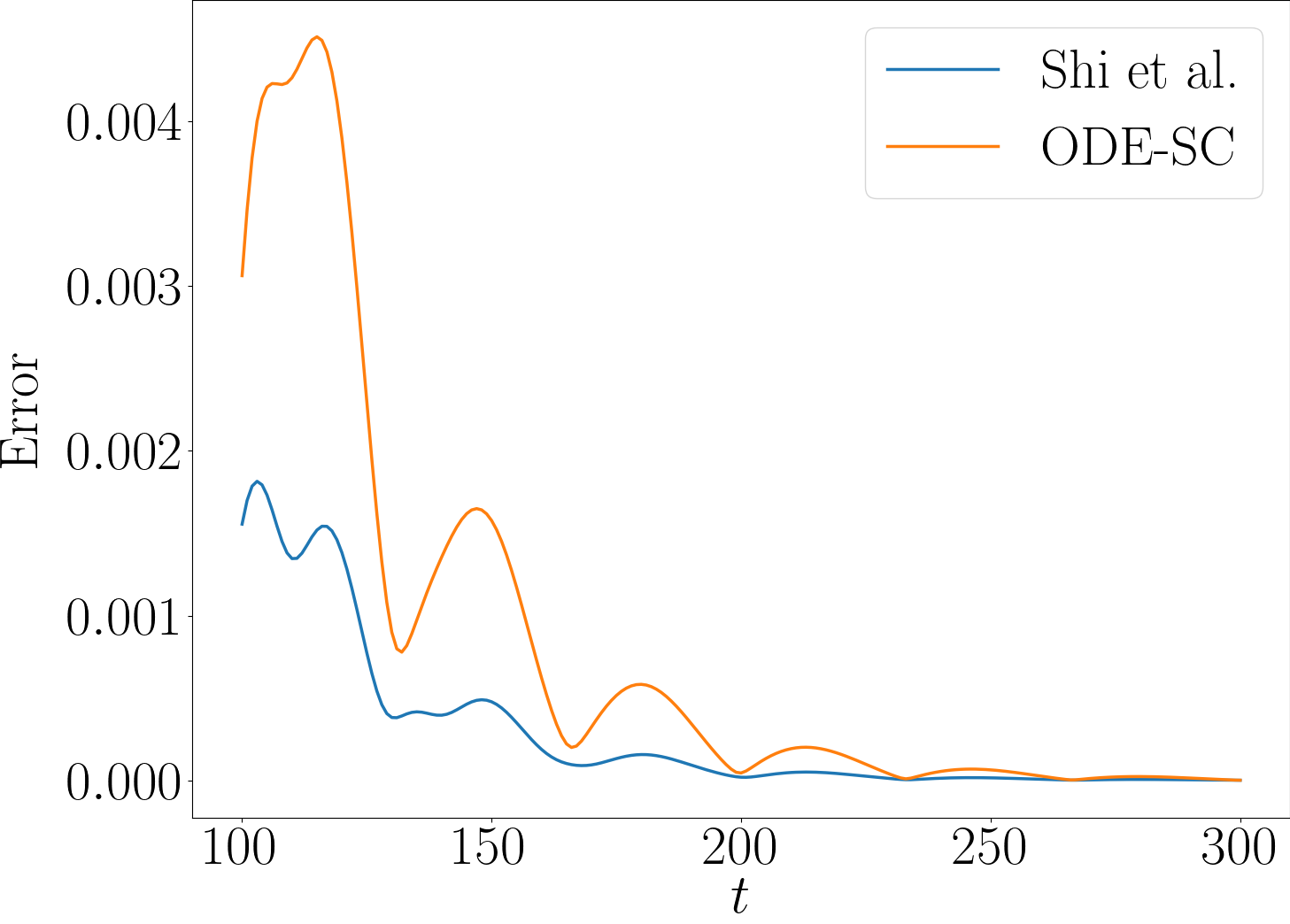}}\\
\end{tabular}
\caption{Comparison of \eqref{eqn:ODE-SC} and existing ODEs \citep{wilson2021lyapunov, shi2021understanding} for (NAG-SC-C).}\label{fig:exp3}
\end{figure}

\begin{figure}
\centering 
\begin{tabular}{cc}
\subfigure[Trajectories]{\includegraphics[width=0.43\textwidth ]{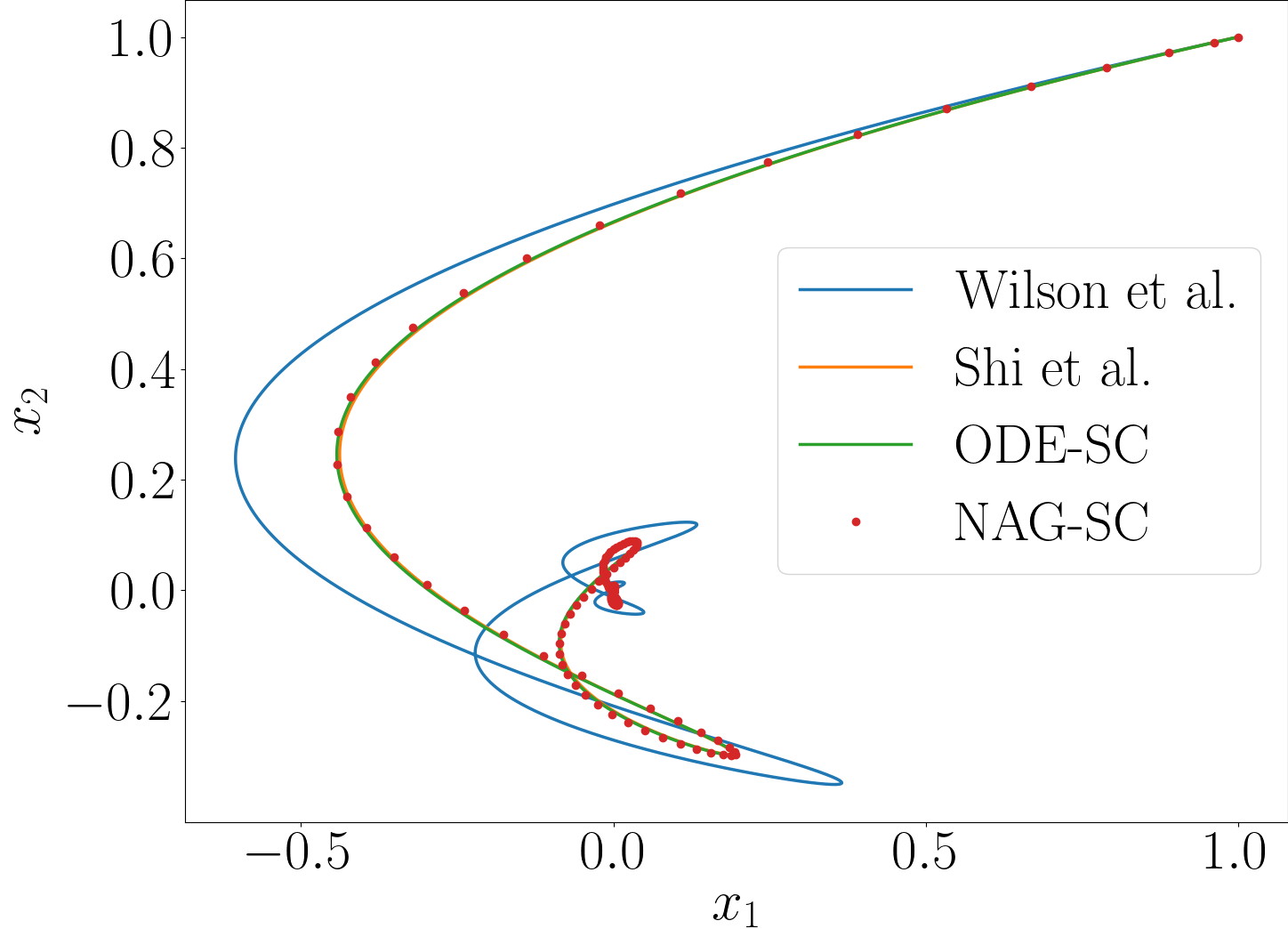}}& \subfigure[Zoomed trajectories]{\includegraphics[width=0.43\textwidth ]{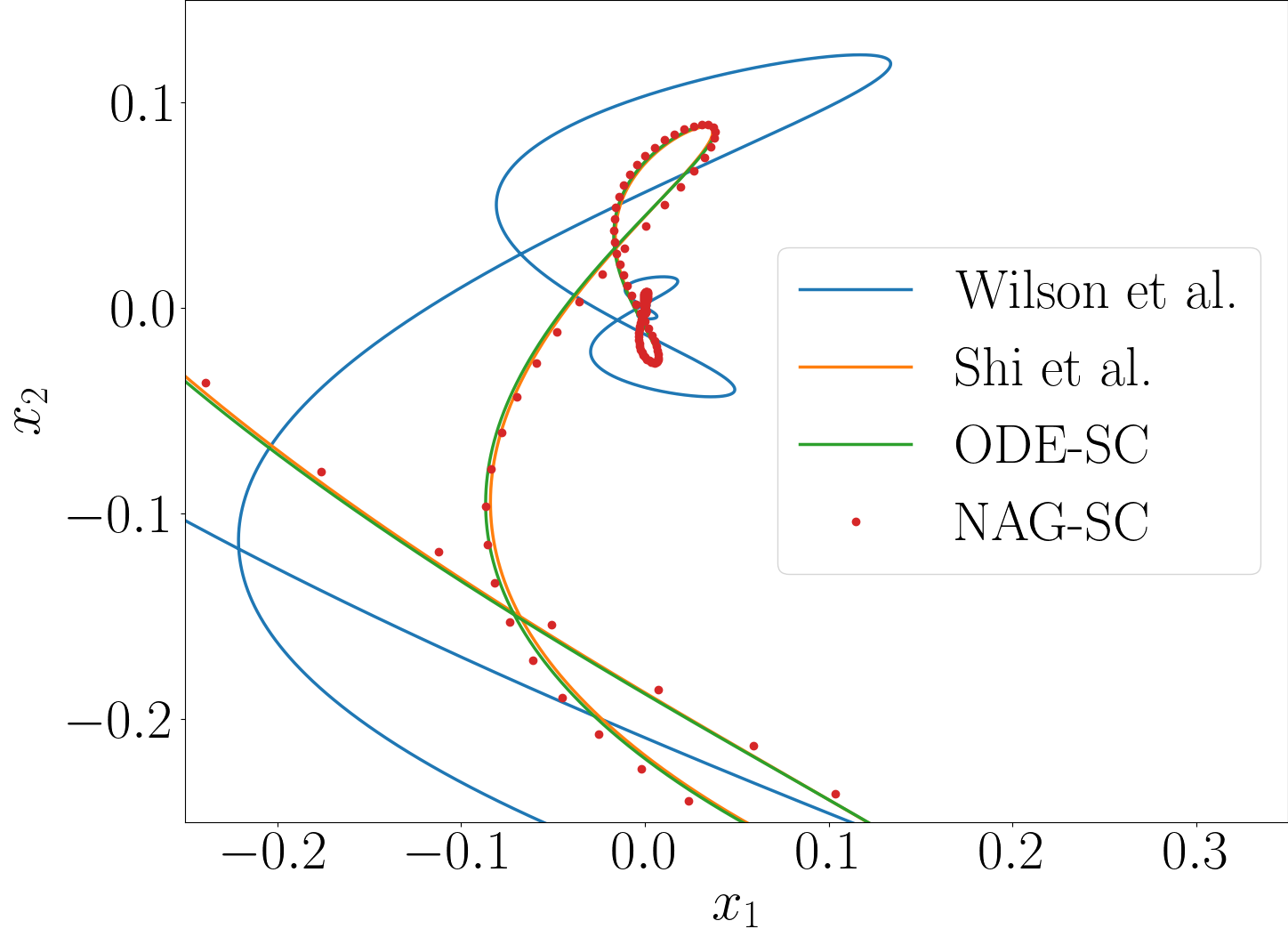}}\\ \subfigure[Errors $f - f^*$]{\includegraphics[width=0.43\textwidth ]{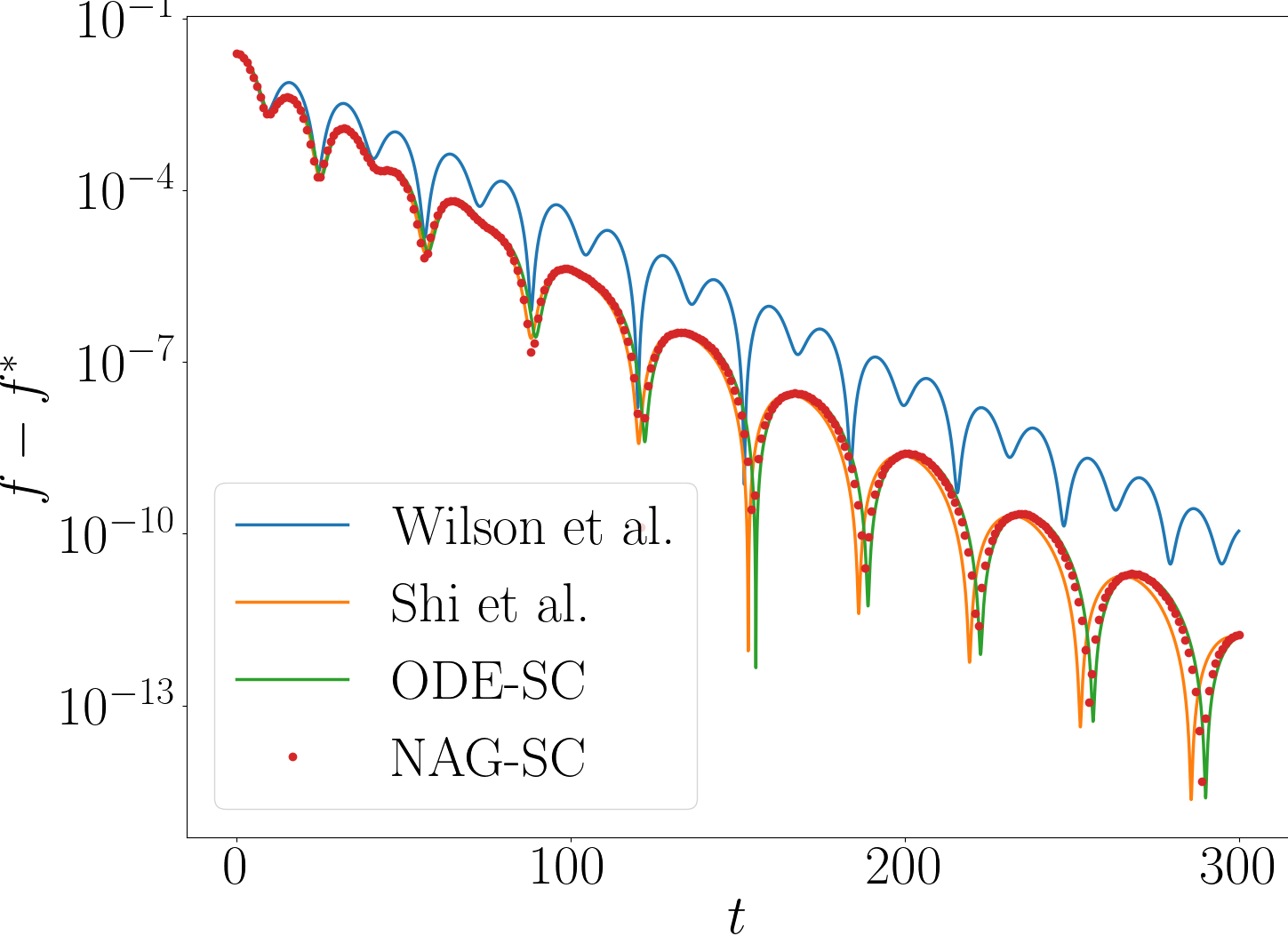}}& \subfigure[Errors $\|X(t=k)-x_k\|_2$]{\includegraphics[width=0.43\textwidth ]{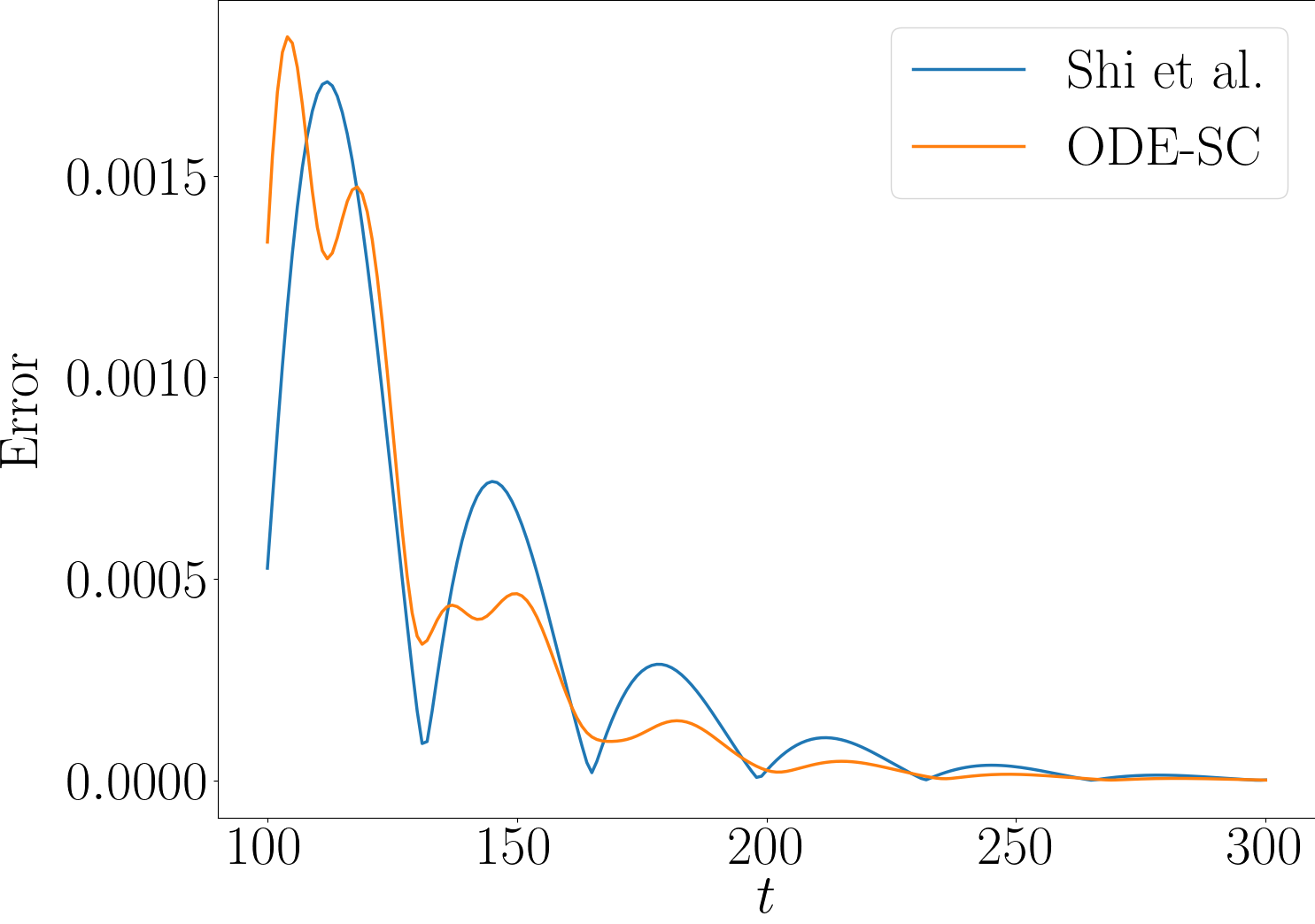}}\\
\end{tabular}
\caption{Comparison of \eqref{eqn:ODE-SC} and  existing models \citep{wilson2021lyapunov,shi2021understanding}
for (NAG-SC)}\label{fig:exp4}
\end{figure}

Figures \ref{fig:exp3} and \ref{fig:exp4} show the results for (NAG-SC-C) and (NAG-SC), respectively. 
In both cases, \eqref{eqn:ODE-SC} provides a more accurate description of  Nesterov's method 
compared to the ODE model 
of \citet{wilson2021lyapunov},
achieving a reduction in the average $L_2$-norm error 
by 82.3\% for (NAG-SC-C) and by 93.7\% for (NAG-SC) over the range $100\le k\le 300$.
This further underscores the importance of computing gradients at future points. While 
\eqref{eqn:ODE-SC} and the model 
of \citet{shi2021understanding}
yield similar trajectories, they do not converge over time due to the convergence condition $\mu\ll L$ being independent of $t$ (see Section \ref{sec:shi}). 
Analyzing the error between Nesterov's methods and the ODEs reveals that \eqref{eqn:ODE-SC} is more accurate in Figure \ref{fig:exp4} than in Figure \ref{fig:exp3}. 
Specifically, for $100\le k\le 300$, the average $L_2$-norm error of \eqref{eqn:ODE-SC} for (NAG-SC) is 65.4\% lower than that of \eqref{eqn:ODE-SC} for (NAG-SC-C).
This discrepancy arises because \eqref{eqn:ODE-SC} directly describes (NAG-SC), while additional  approximations are required in \eqref{eq:sbsb} for (NAG-SC-C).

\subsection{Changes in discretization interval $h=t/k$}\label{sec:exp3}

In this subsection, we empirically validate Proposition \ref{thm:Nes-Flow}, specifically investigating whether our ODE models and Nesterov's methods converge as $h$ decreases. We also compare our models with the BLFs \eqref{eqn:BLF-C} and \eqref{eqn:BLF-SC}, derived by replacing $\nabla f(Y)$ with $\nabla f(X)$ in \eqref{eqn:G-ODE-C} and \eqref{eqn:G-ODE-UC}, respectively.
The simulation parameters are consistent with those in previous subsections, with the  exception of $h$. 

\begin{figure*}[ht]
	\centering
	\subfigure[$h=1$]{\includegraphics[width=0.3\textwidth]{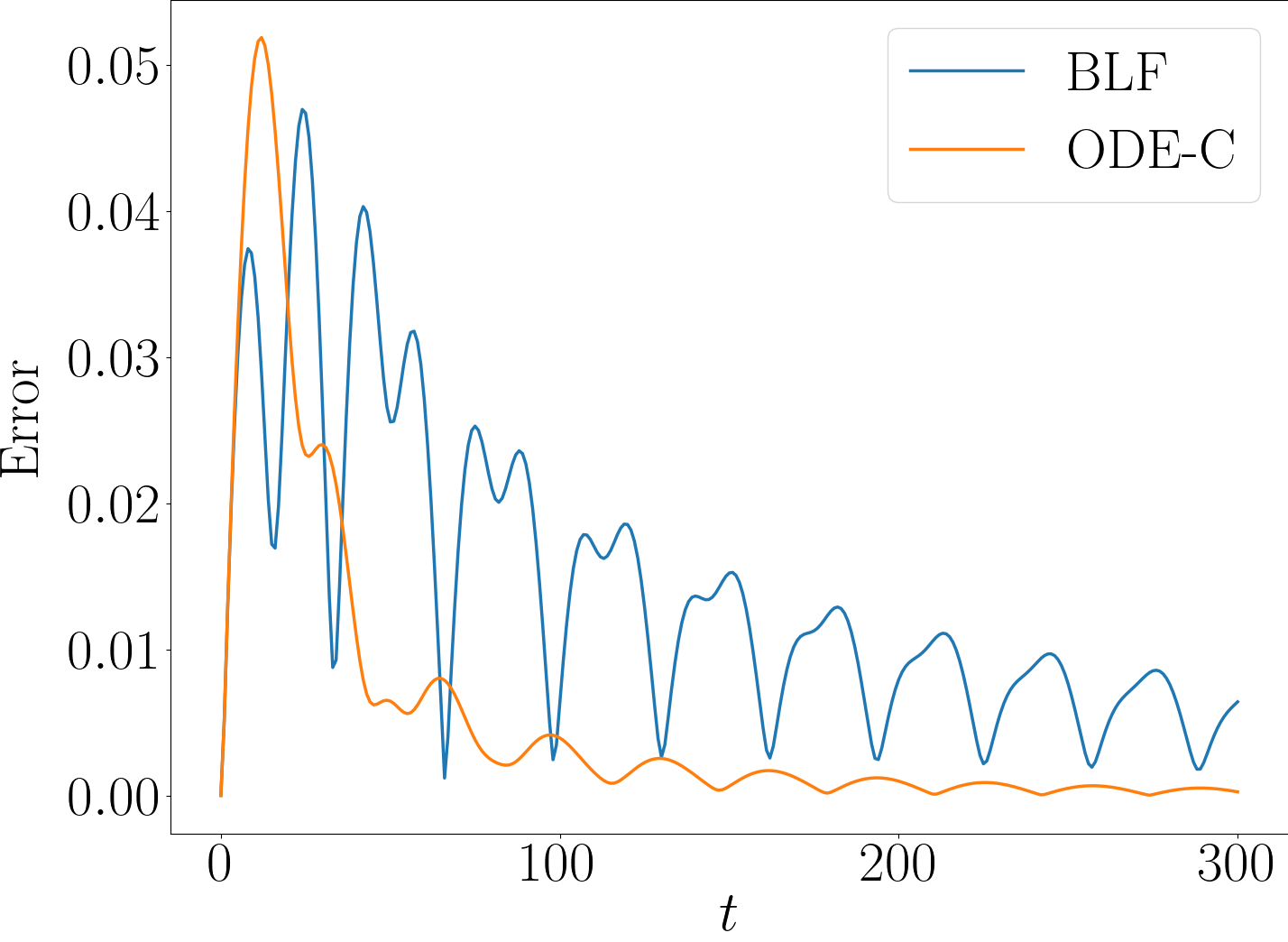}} \hspace{0.1in}
	\subfigure[$h=0.1$]{\includegraphics[width=0.3\textwidth]{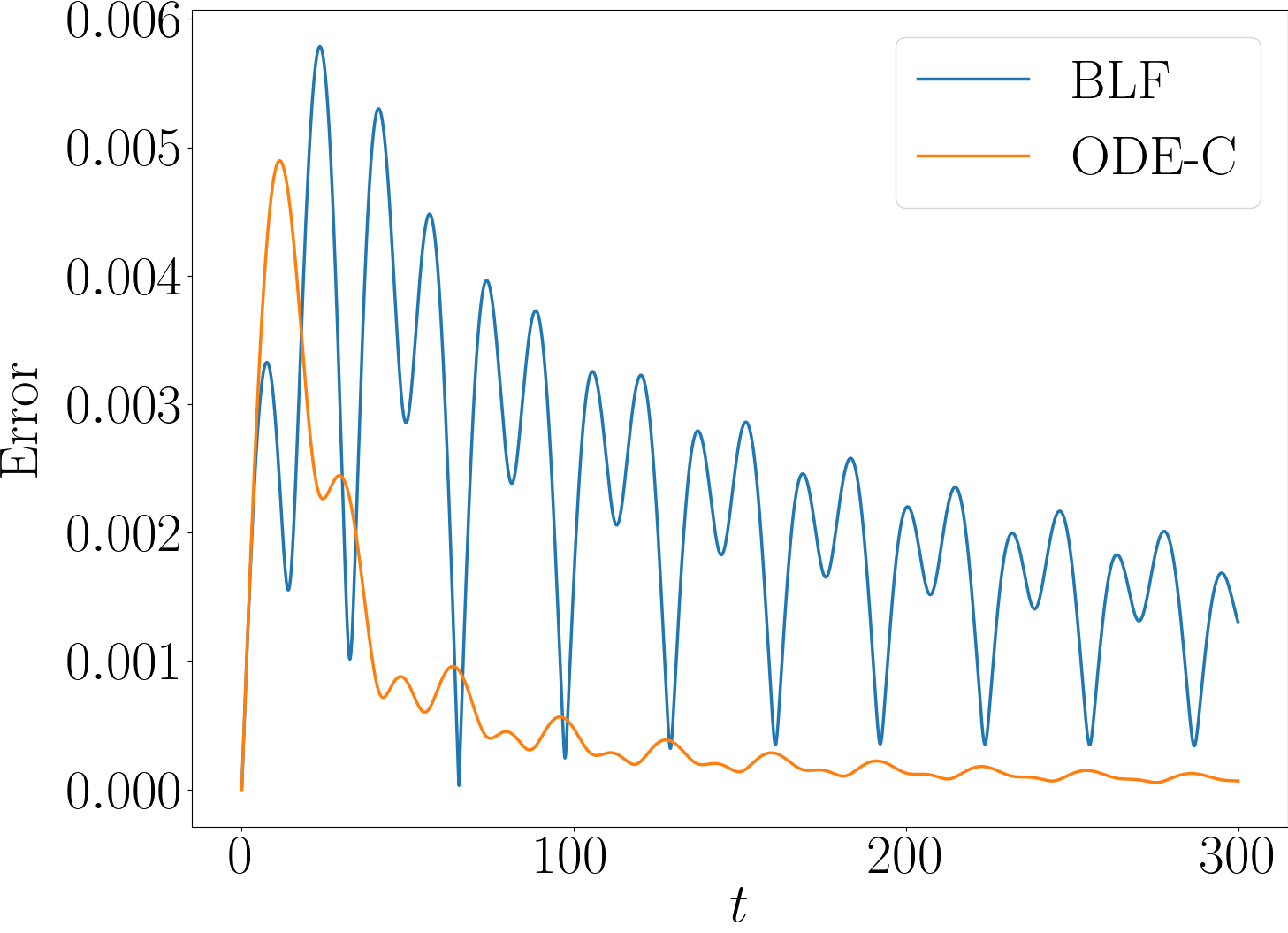}} \hspace{0.1in}
	\subfigure[$h=0.01$]{\includegraphics[width=0.3\textwidth]{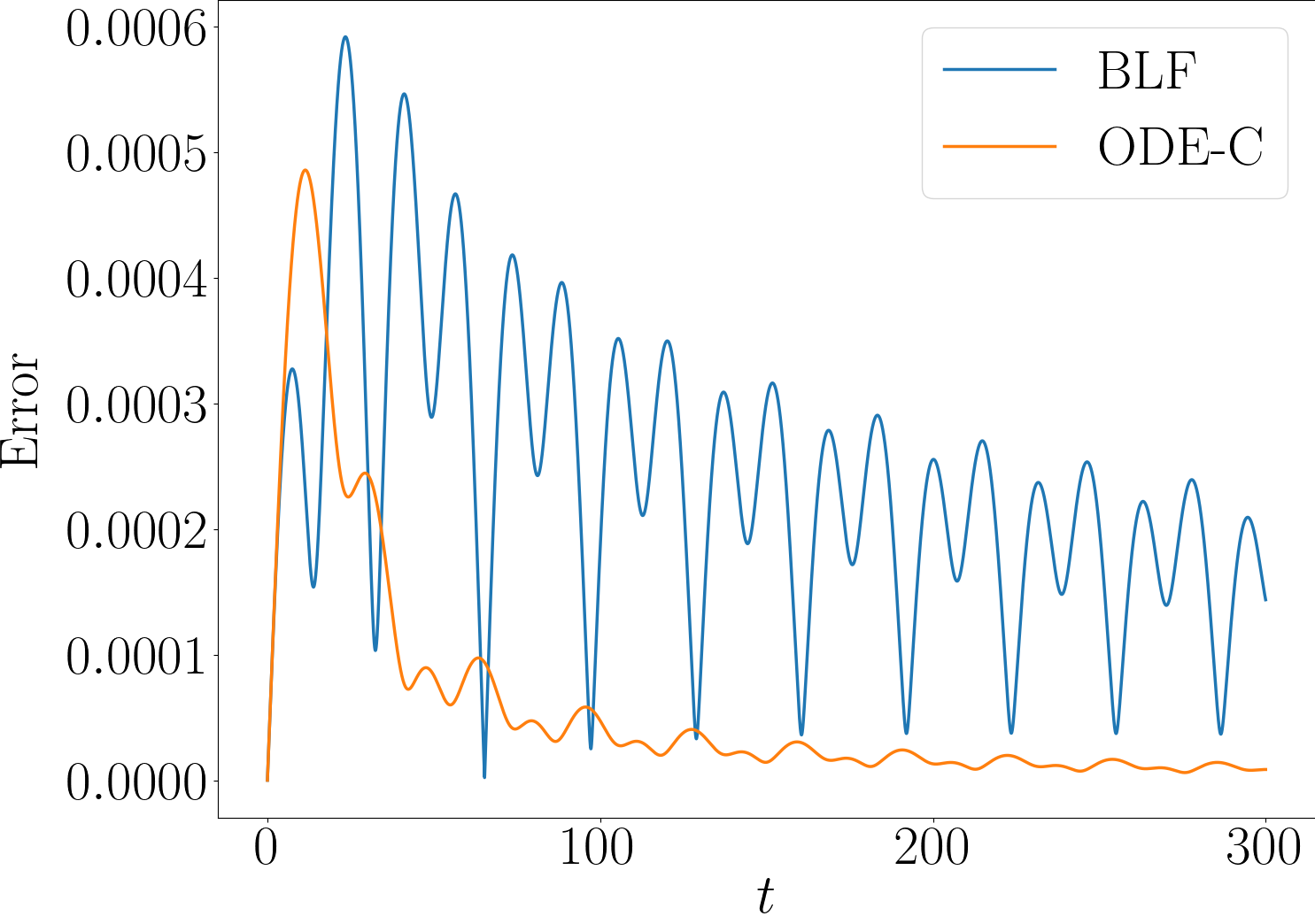}}
	\caption{Errors $\|X(t=hk)-x_k\|_2$ of 	
	our model and BLF compared to the iterates of Nesterov's method for the non-strongly convex case with $f(x_1,x_2)=0.02 x_1^2 + 0.005 x_2^2$.} 
	\label{fig:exp5}
\end{figure*}

\begin{figure*}[ht]
	\centering
	\subfigure[$h=1$]{\includegraphics[width=0.3\textwidth]{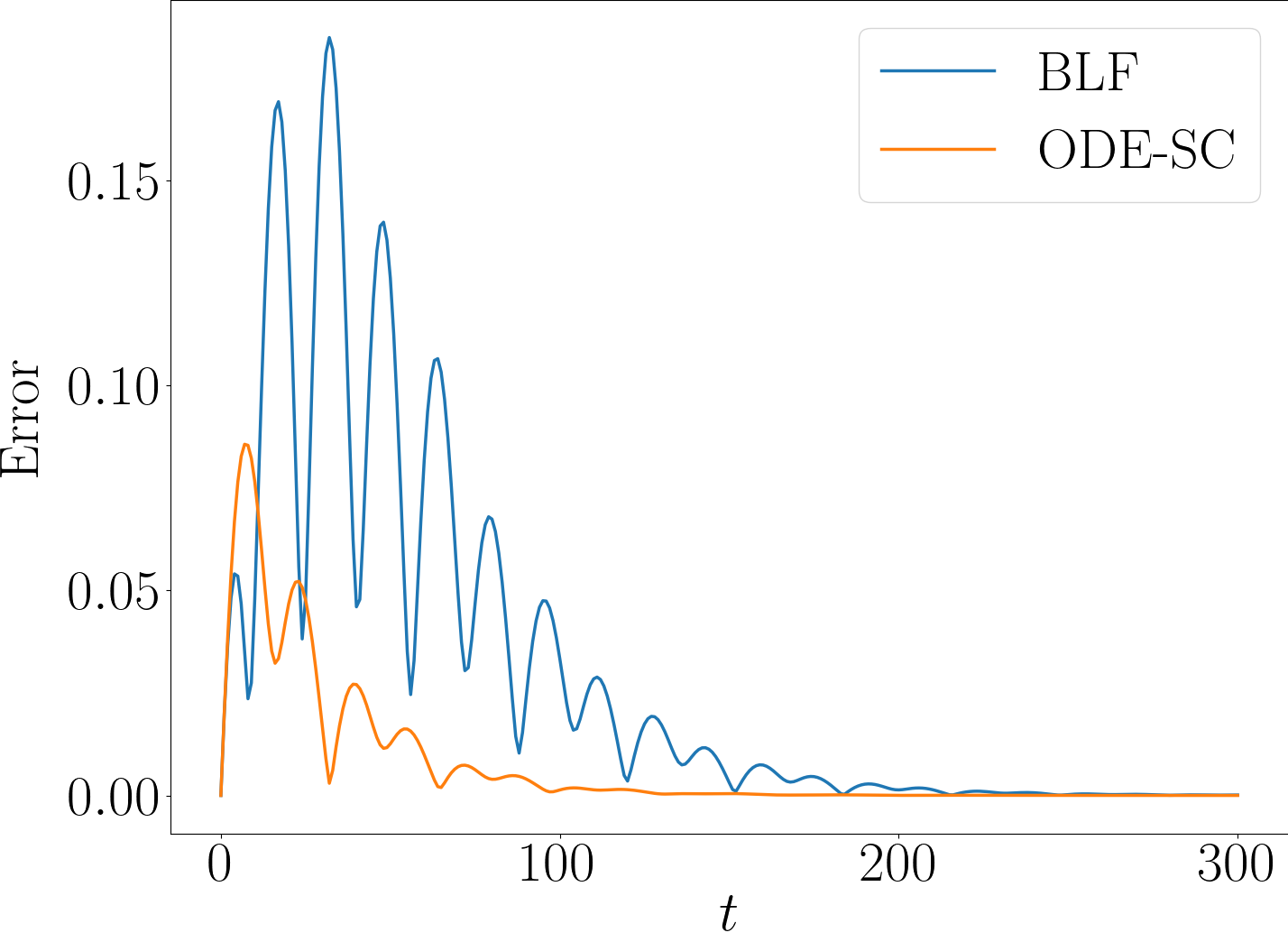}} \hspace{0.1in}
	\subfigure[$h=0.1$]{\includegraphics[width=0.3\textwidth]{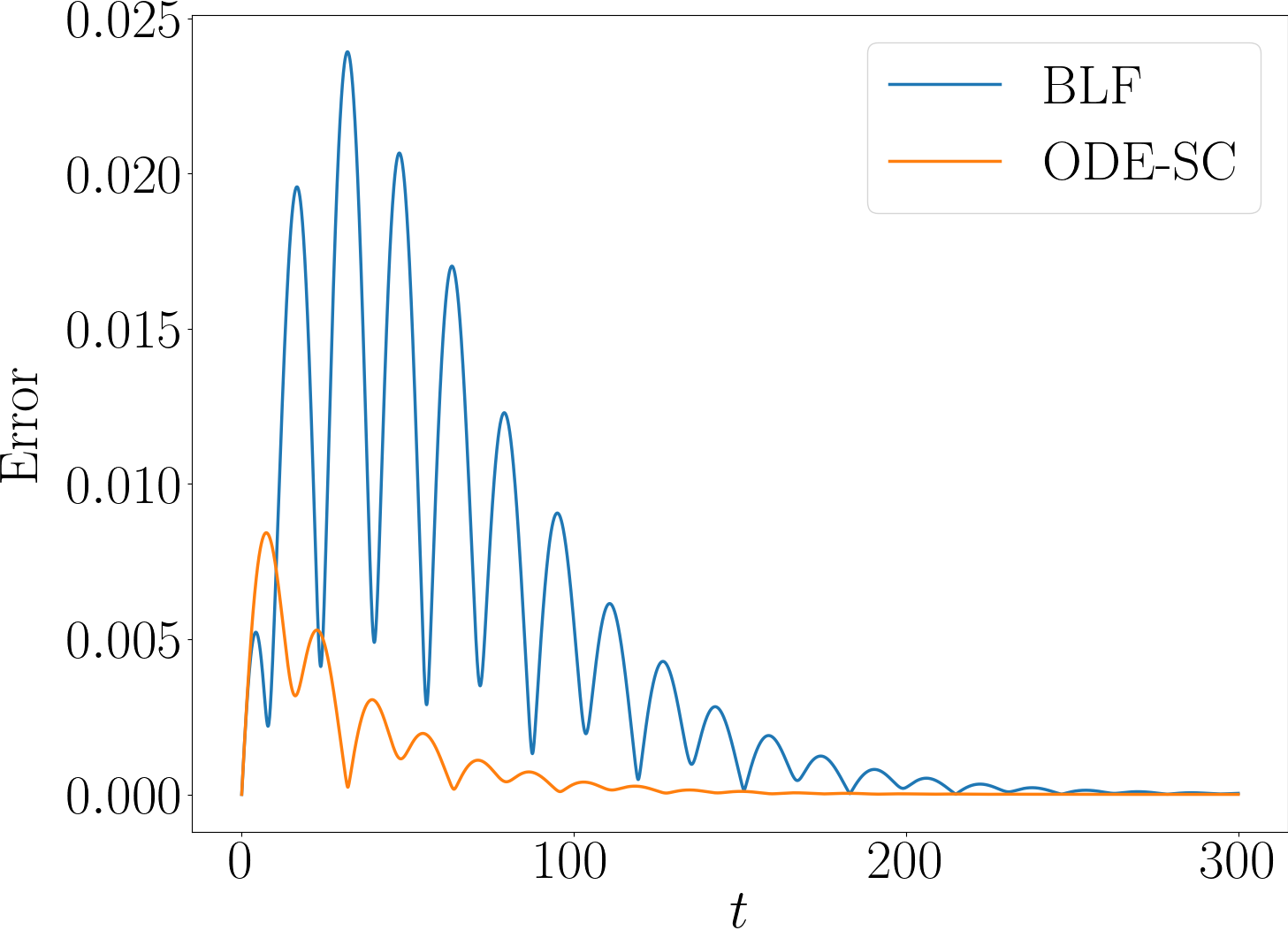}} \hspace{0.1in}
	\subfigure[$h=0.01$]{\includegraphics[width=0.3\textwidth]{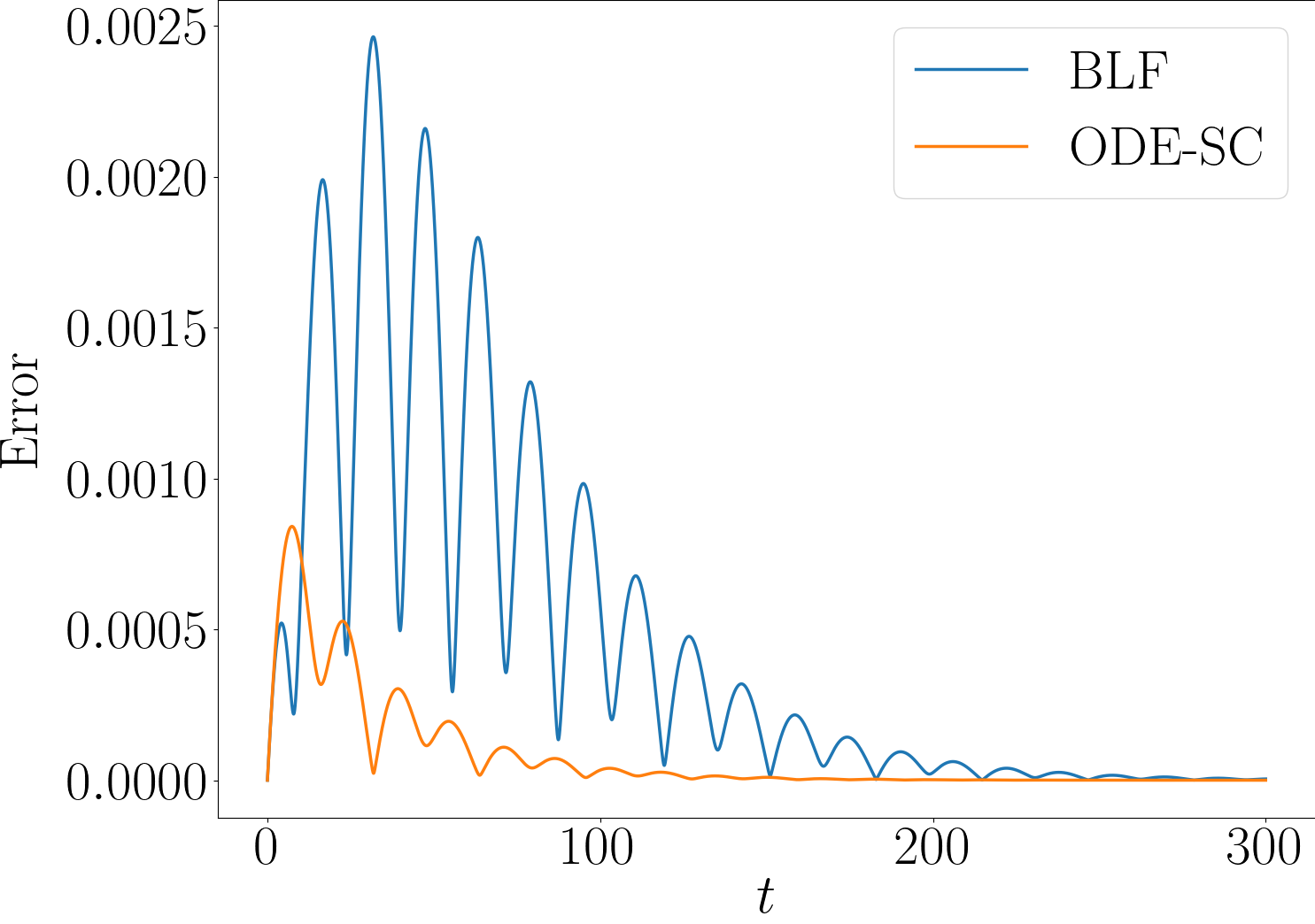}}
	\caption{Errors $\|X(t=hk)-x_k\|_2$ of our model and BLF compared to the iterates of Nesterov's method for the strongly convex case with $f(x_1,x_2)=0.02 x_1^2 + 0.005 x_2^2$.}
	\label{fig:exp6}
\end{figure*}

\begin{figure*}[ht]
	\centering
	\subfigure[$h=1$]{\includegraphics[width=0.3\textwidth]{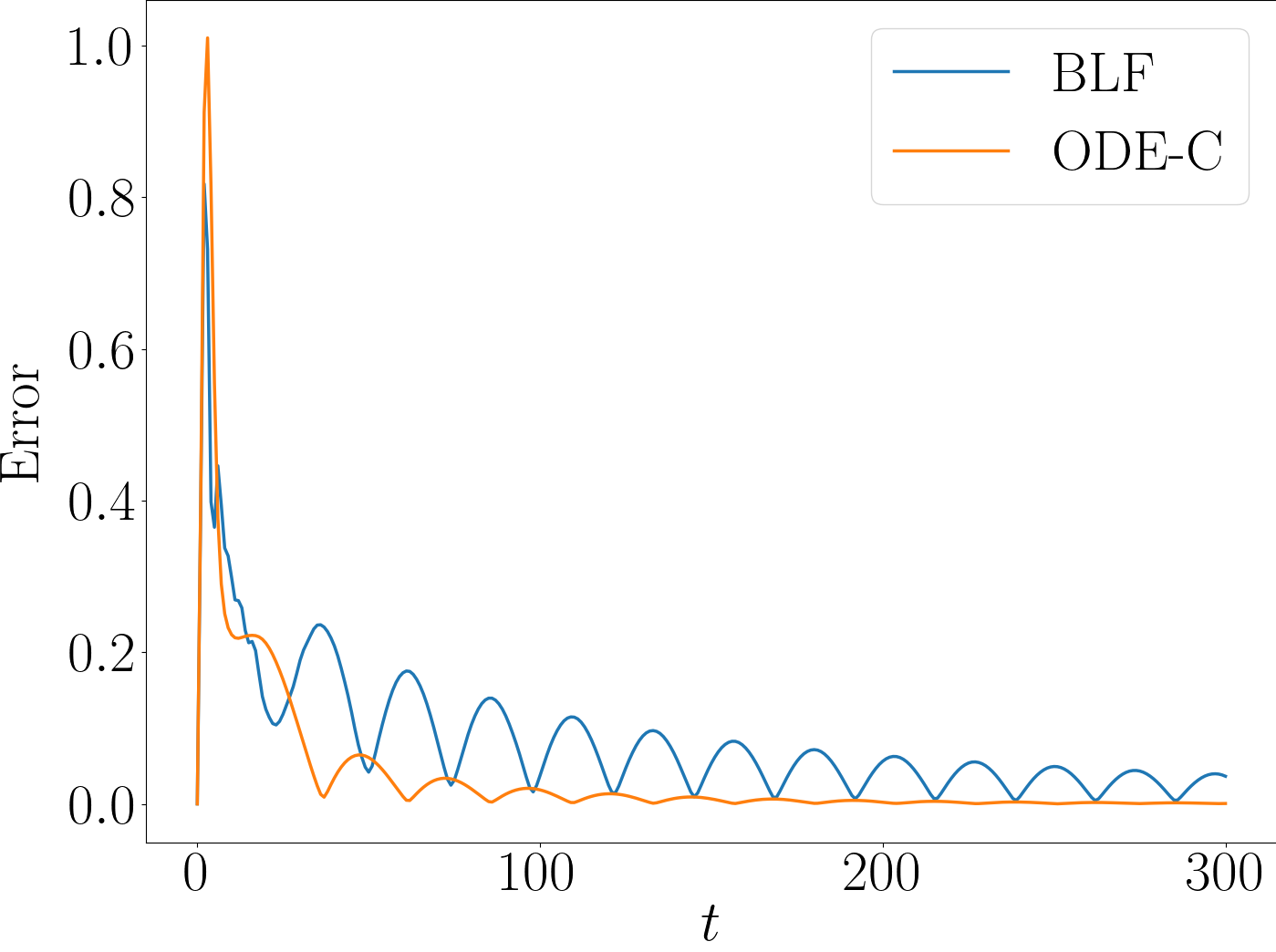}} \hspace{0.1in}
	\subfigure[$h=0.1$]{\includegraphics[width=0.3\textwidth]{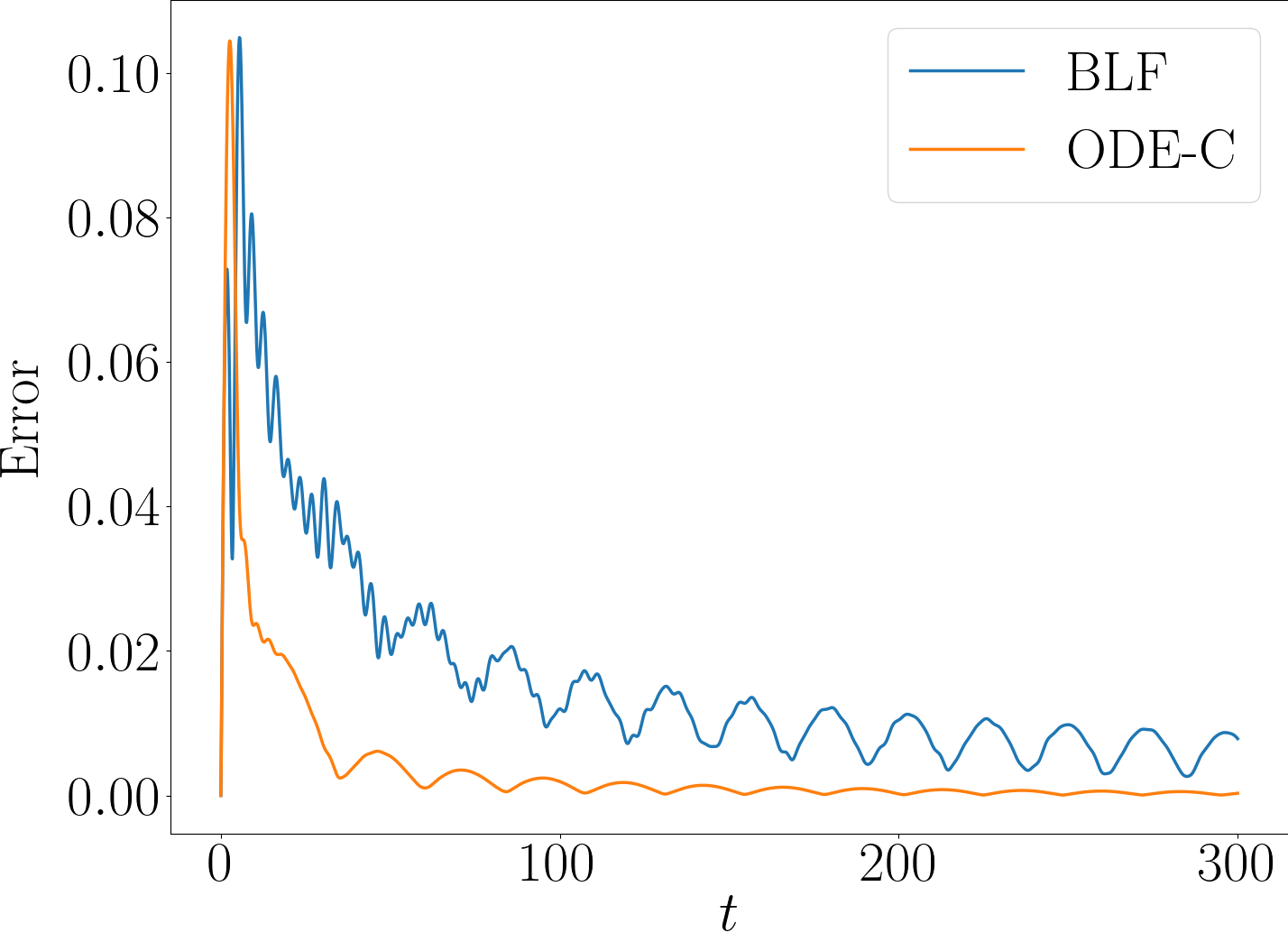}} \hspace{0.1in}
	\subfigure[$h=0.01$]{\includegraphics[width=0.3\textwidth]{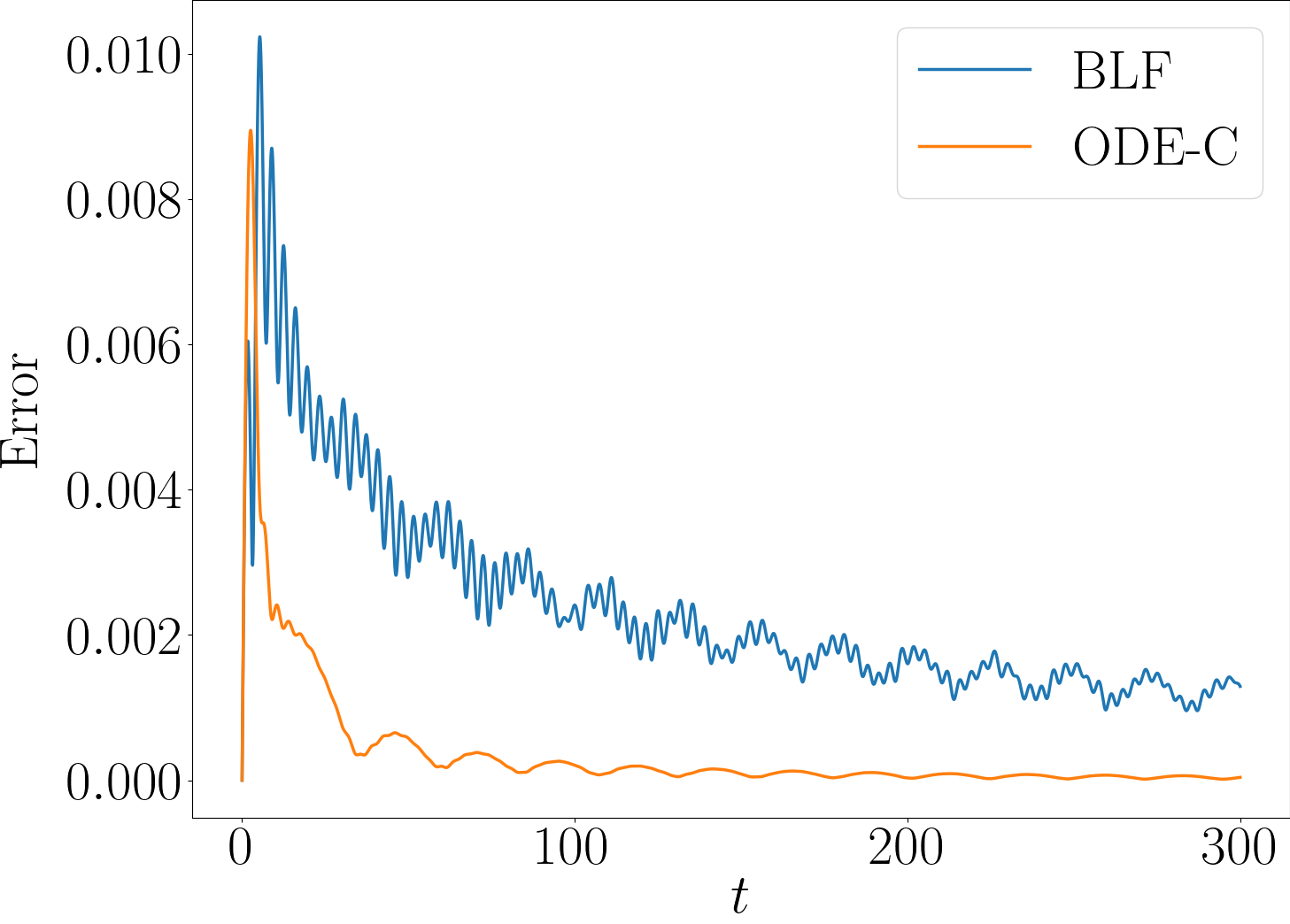}}
	\caption{Errors $\|X(t=hk)-x_k\|_2$ of our model and BLF compared to the iterates of Nesterov's method for the non-strongly convex case with $f(x)=\frac{1}{2}x^T M x$.}
	\label{fig:exp7}
\end{figure*}

\begin{figure*}[ht]
	\centering
	\subfigure[$h=1$]{\includegraphics[width=0.3\textwidth]{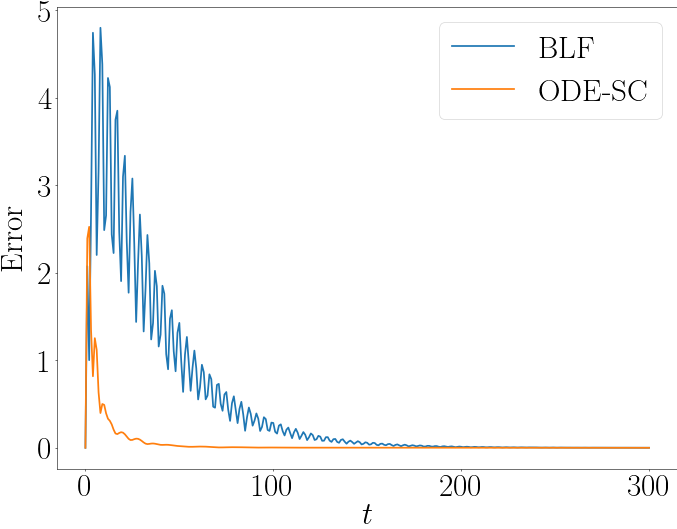}} \hspace{0.1in}
	\subfigure[$h=0.1$]{\includegraphics[width=0.3\textwidth]{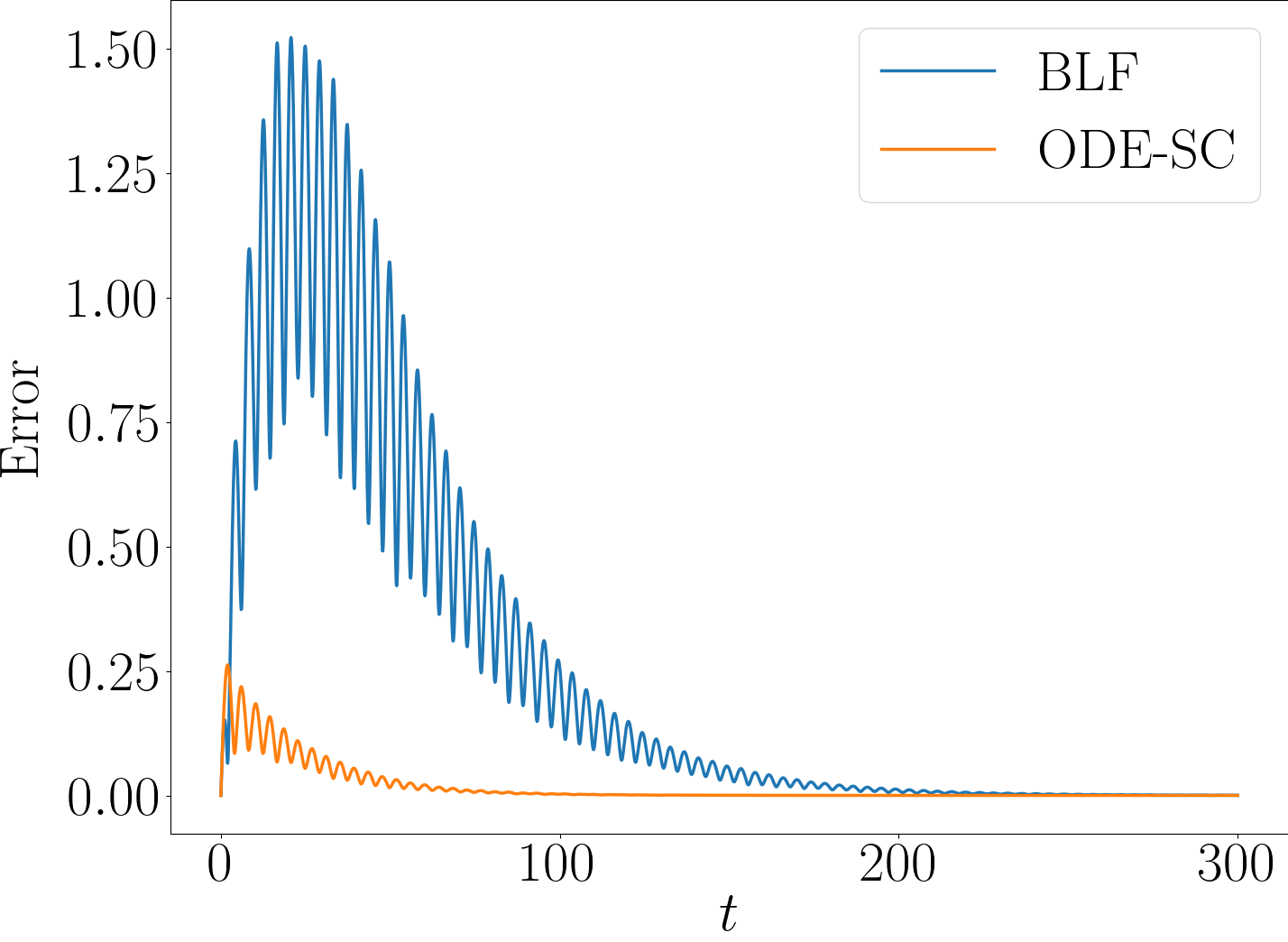}} \hspace{0.1in}
	\subfigure[$h=0.01$]{\includegraphics[width=0.3\textwidth]{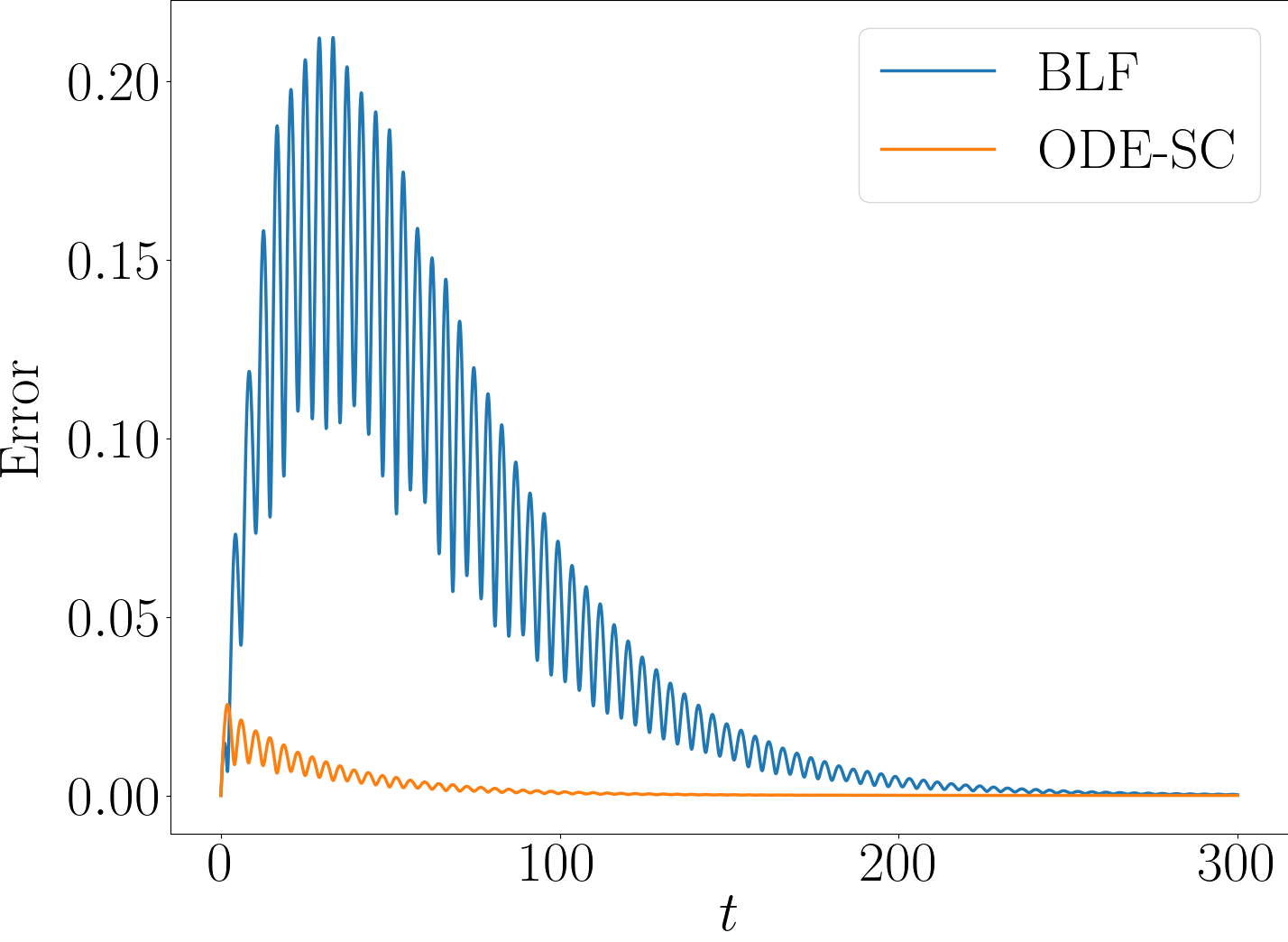}}
	\caption{Errors $\|X(t=hk)-x_k\|_2$ of our model and BLF compared to the iterates of Nesterov's method for the strongly convex case with $f(x)=\frac{1}{2}x^T M x$.}\label{fig:exp8}
\end{figure*}

Figures \ref{fig:exp5} and \ref{fig:exp6} display the results, showing deviation errors computed in terms of the $L_2$-norm. In addition, we repeat the experiments with a different objective function $f(x) = \frac{1}{2}x^T M x$, where the positive-definite matrix $M\in\mathbb{R}^{n\times n}$ has a random eigenvalue ranging between $\mu$ and $L$. We set $n=200$ and $L = 1$, with $\mu = 0$ for the convex case and $\mu = 0.001$ for the $\mu$-strongly convex case. The starting point $X(0)$ is chosen randomly  with $\dot{X}(0)=(0,0)$. The results are shown in Figures \ref{fig:exp7} and \ref{fig:exp8}. 

From the simulation results, it is evident that the deviation error decreases with $h$, in line with Proposition~\ref{thm:Nes-Flow}. Furthermore,   our ODE model's trajectory more closely approximates the iterates of Nesterov's method compared to the BLFs for all values of $h$. This improved alignment is due to  our models incorporating  gradient correction, which aligns with the fundamental principle of Nesterov's methods: calculating the gradient at the predicted future position.

\subsection{Restart schemes}\label{sec:restartexp}
\begin{figure*}[ht]
	\centering
	\subfigure[$f(x_1,x_2)=0.02 x_1^2 + 0.005 x_2^2$]{\includegraphics[width=0.4\textwidth]{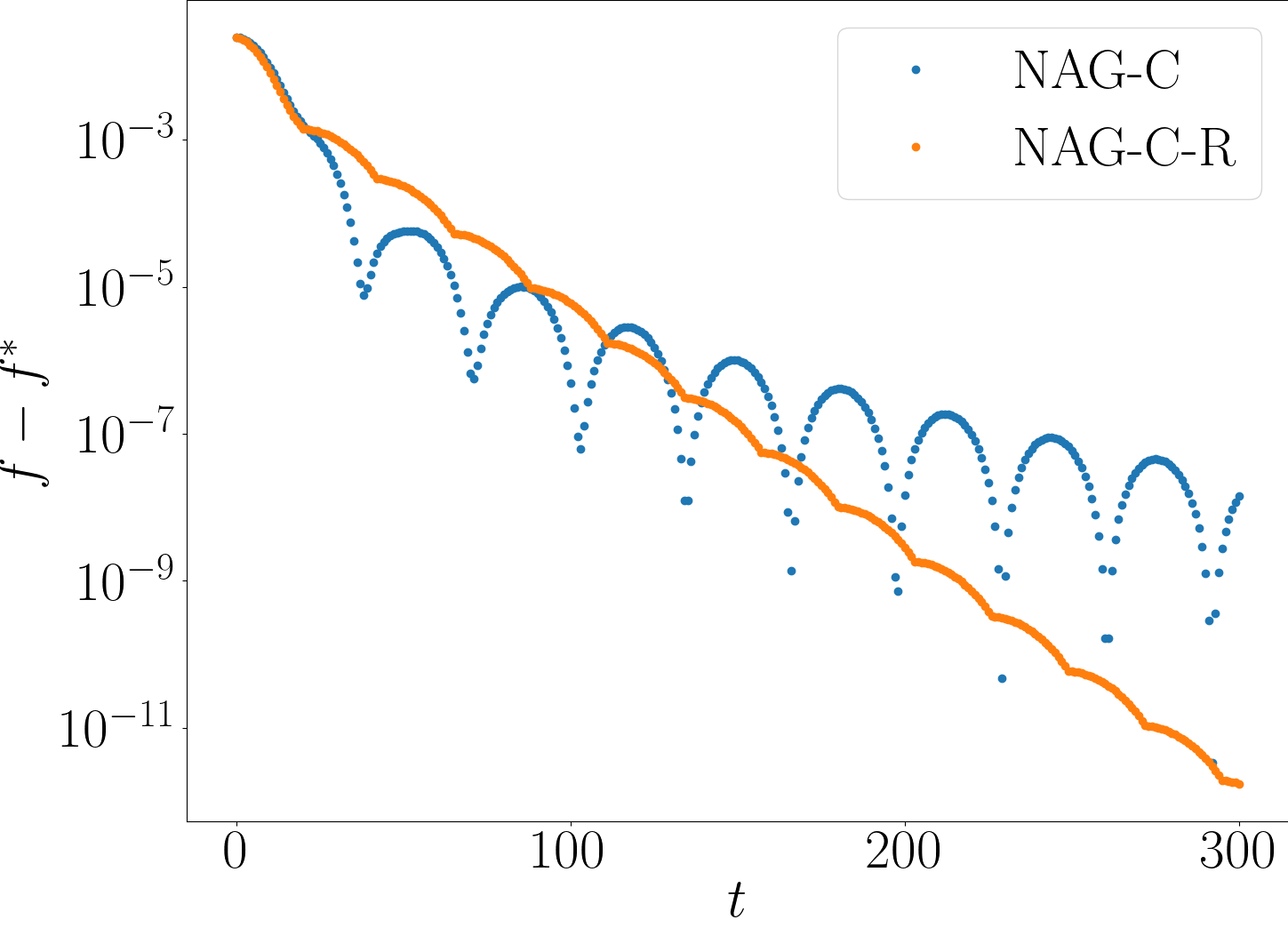}} \hspace{0.1in}
	\subfigure[$f(x)=\frac{1}{2}x^T M x$]{\includegraphics[width=0.4\textwidth]{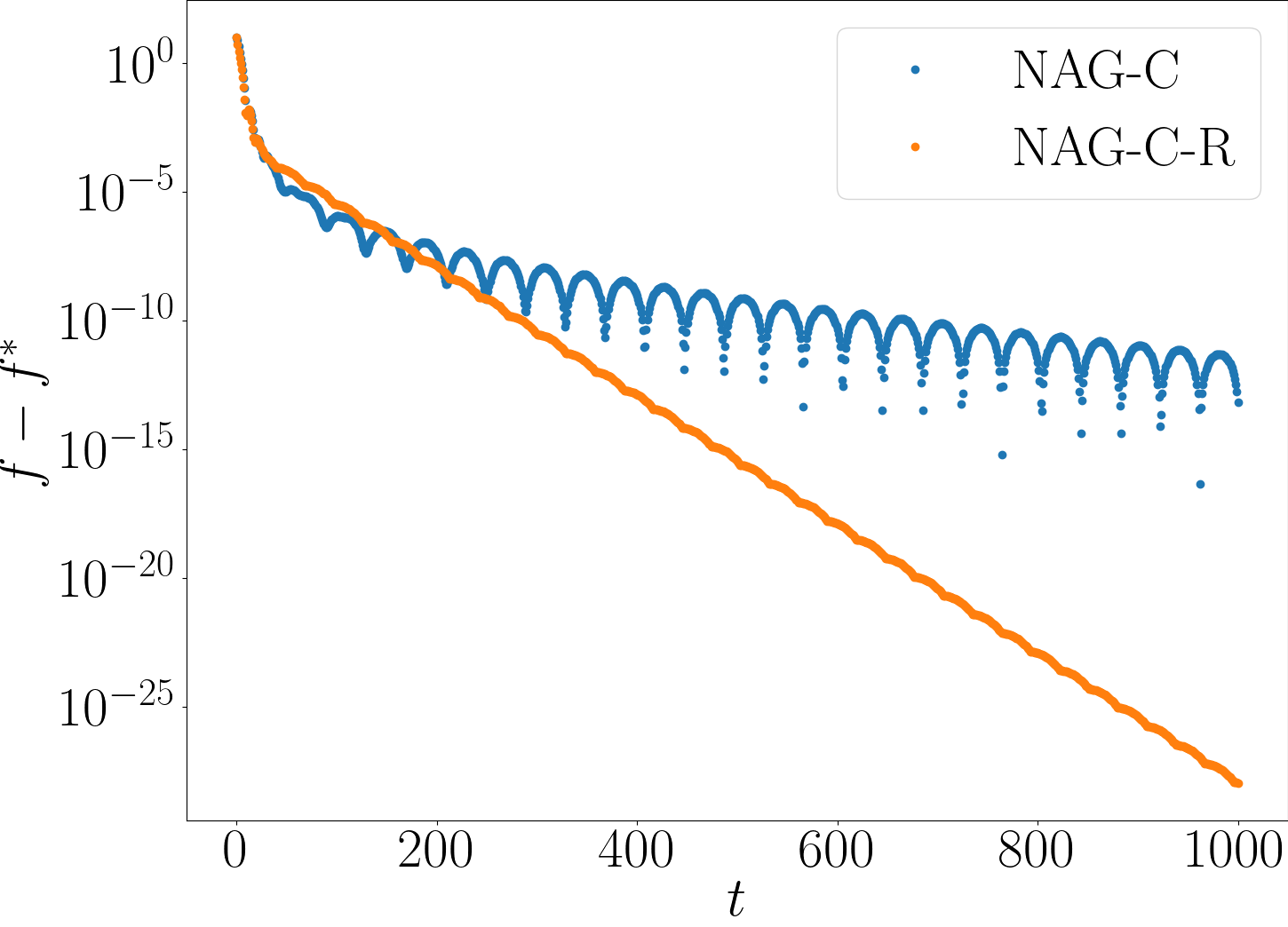}} 
	\caption{Error trajectories of (NAG-C) and (NAG-C-R).}\label{fig:exp-r}
\end{figure*}

In this subsection, we first examine the performance of 
our restart scheme, Algorithm \ref{alg:cap}, for (NAG-C) using the same settings as in Sections \ref{sec:exp1} and \ref{sec:exp3}. 
Figure \ref{fig:exp-r} illustrates the error trajectories for (NAG-C) and (NAG-C-R), with the latter incorporating our restart scheme 
with $k_{min}=20$. The results indicate that the restart scheme improves the performance of (NAG-C) by preventing fluctuations. 

\begin{figure}[t]
	\begin{center}\centerline{\includegraphics[width=0.5\columnwidth]{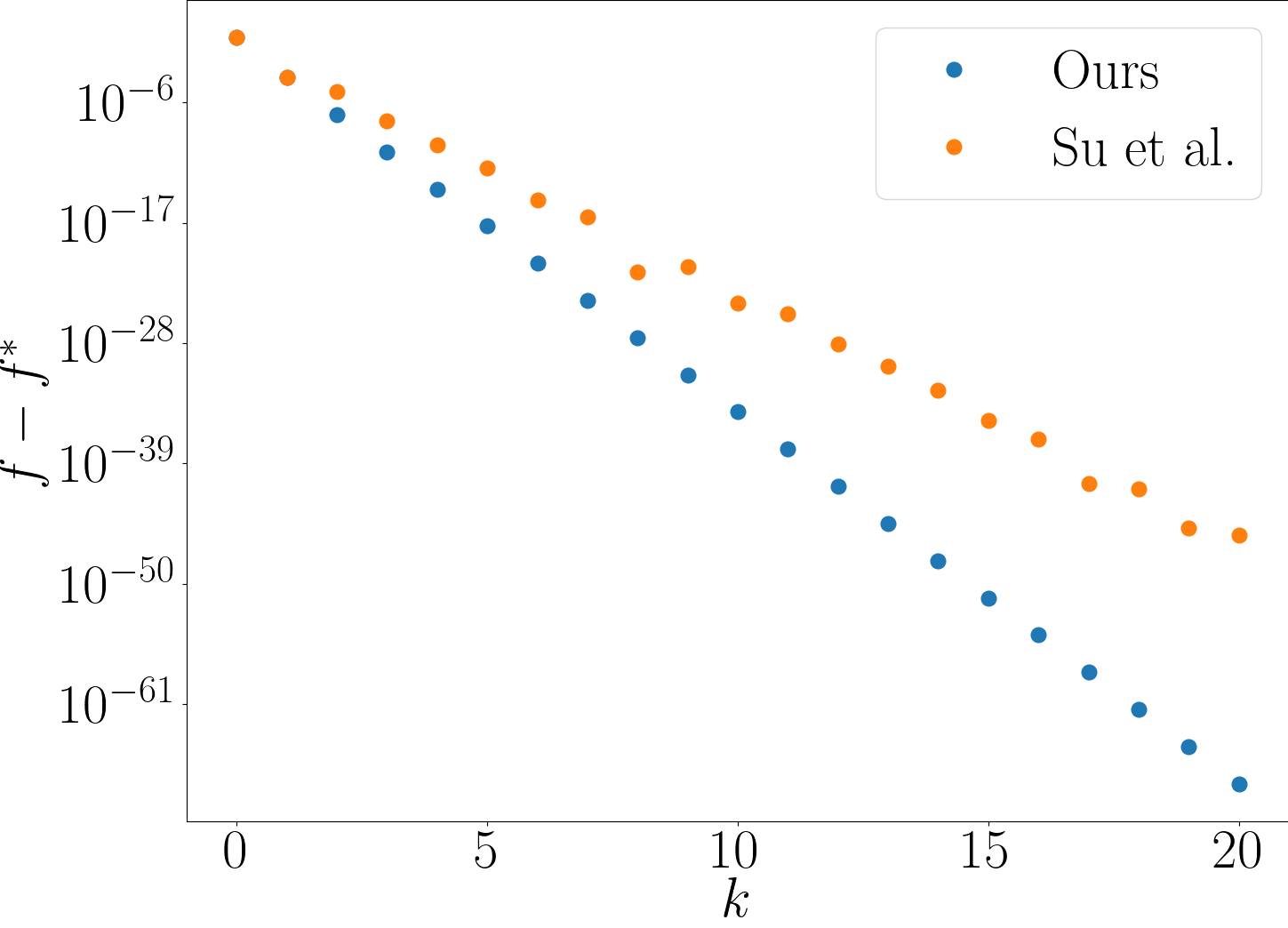}}
\caption{Error trajectories of (NAG-C-C) obtained using  our restart scheme and the restart scheme 
of \citet{su2016differential}.}
		\label{fig:counter}
	\end{center}
\end{figure}

Additionally, we compare our restart scheme to the  method proposed by~\citet{su2016differential} for (NAG-C-C). We consider the objective function $f(x_1, x_2)=0.5x_1^2+0.49x_2^2$ with the initial point $x_0 = (1, 1)$ and set $k_{min}=1$. The error trajectories of (NAG-C-C) generated using the two restart schemes are shown in Figure \ref{fig:counter}. We observe that the objective values obtained using the restart scheme 
of \citet{su2016differential}
increase at $k=9$. More specifically, the objective value at $k=9$ is approximately 2.8 times greater than that at $k=8$. On the contrary, the objective values obtained using our restart scheme  decrease monotonically, which is consistent with the theoretical analysis in Section \ref{sec:restart}.

\subsection{Connections to gradient descent}
\begin{figure*}[ht]
	\centering
	\subfigure[ODE-SC and \eqref{eqn:time-XZ}]{\includegraphics[width=0.4\textwidth]{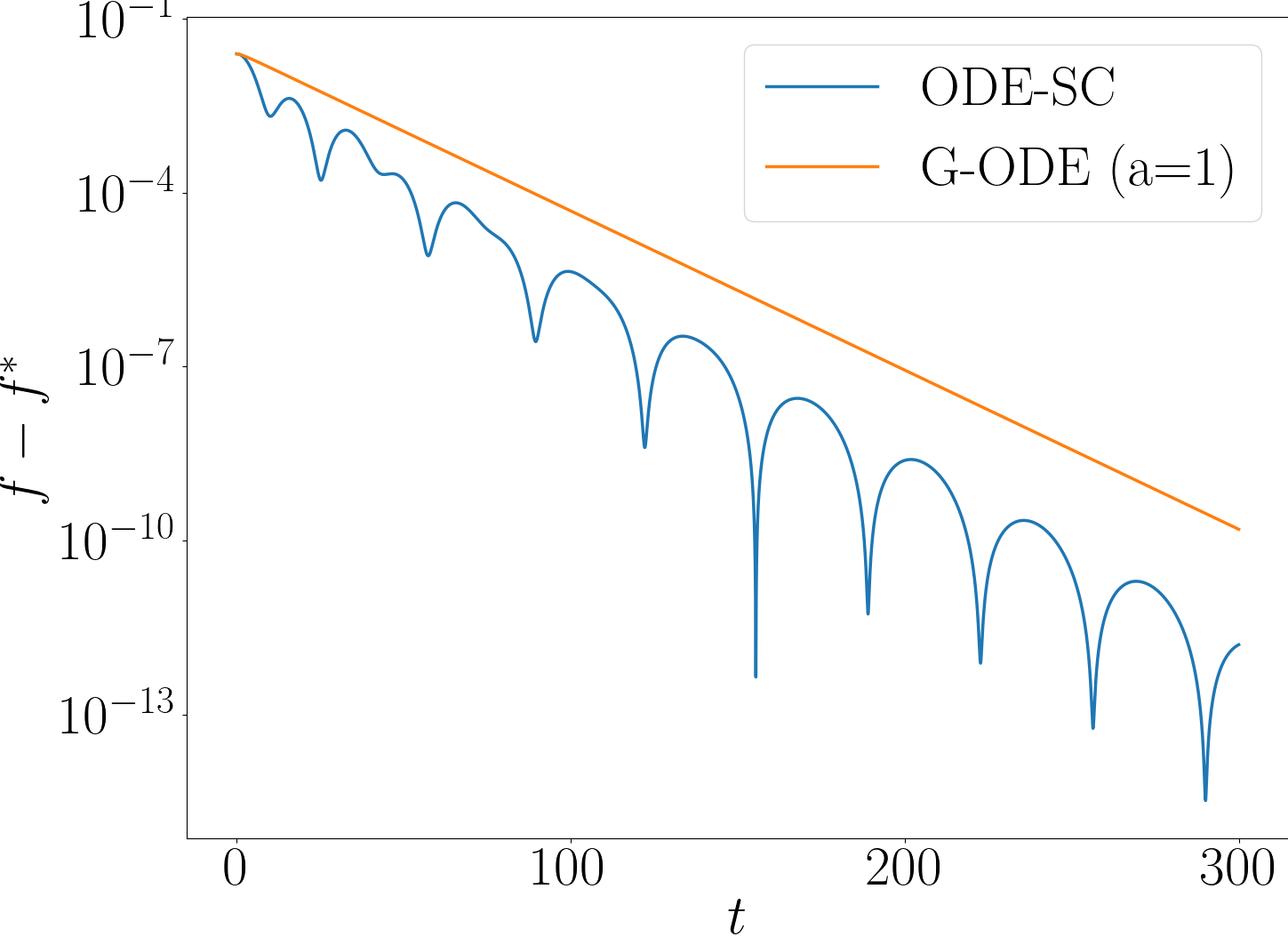}} \hspace{0.1in}
	\subfigure[Error $\|X(t=\tau^{-1}(k))-r_k\|_2$]
	{\includegraphics[width=0.4\textwidth]{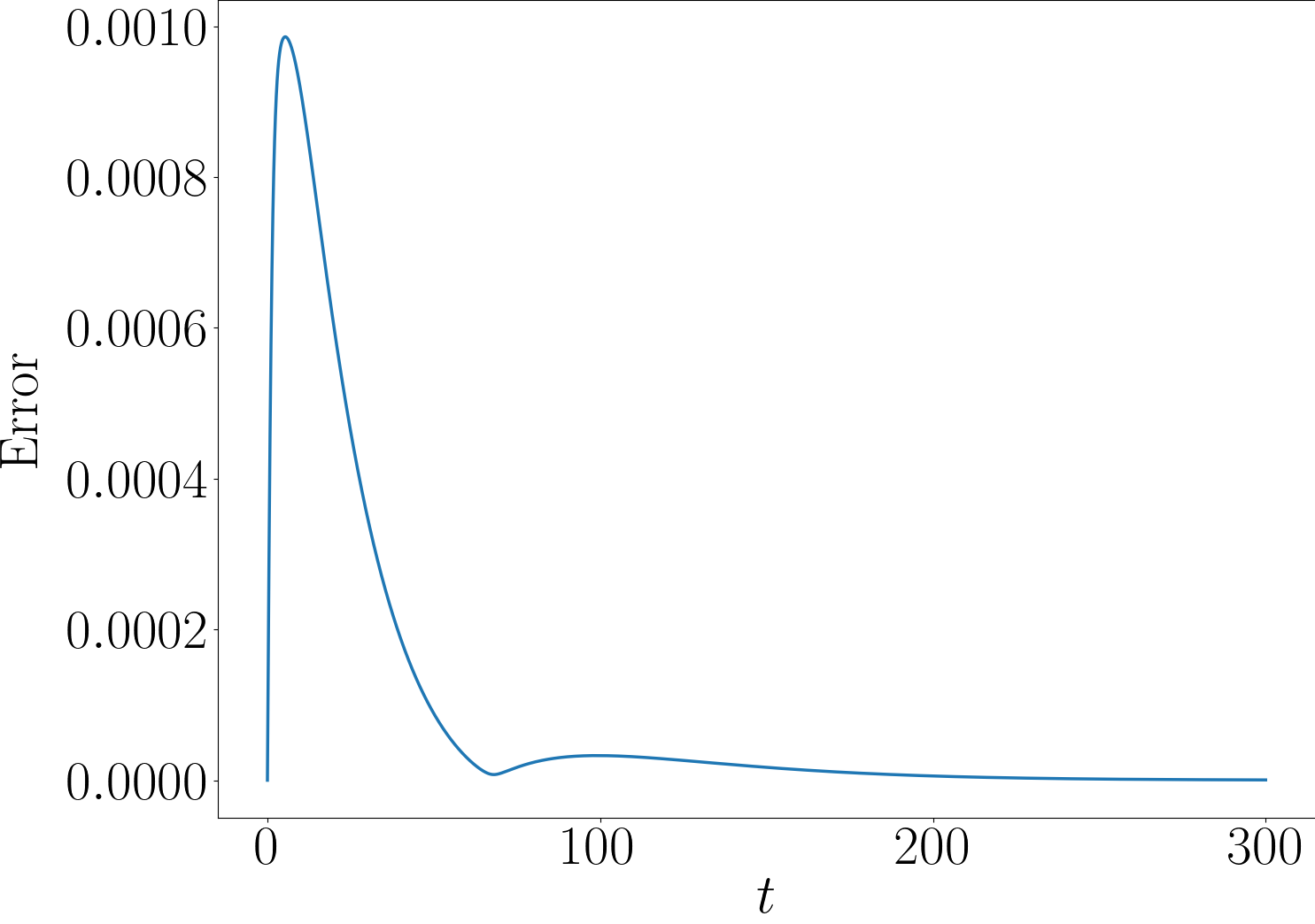}} 
	\caption{Numerical verifications of the discussions in Section \ref{sec:time}.}
\label{fig:exp-t}
\end{figure*}

In this subsection, we explore the relationship between our generalized  models and gradient flow as described in Section \ref{sec:time} using the same experiment settings as in Section \ref{sec:6.2}. First, we compare \eqref{eqn:ODE-SC} from Section \ref{sec:comp:wilson} with our generalized model, where ${\bf a}(t)=1$.
This comparison helps elucidate 
  the connection between Nesterov's methods and gradient descent. The ODE model \eqref{eqn:ODE-SC} is rewritten as follows: 
\begin{equation}\label{eqn:ODE-SC-time}
\begin{split}
    Z &= X + \sqrt{\frac{L}{\mu}}\dot{X} \\
    Y &= X + \ac(t) \frac{L}{\mu}\dot{X} \\
    \dot{Z} &= -(1-\ac(t))\dot{X} - \frac{1}{\sqrt{\mu L}}\nabla f(Y)
\end{split}
\end{equation}
where
\begin{equation*}
    \ac (t) = \frac{e^{\sqrt{\frac{\mu}{L}}h}-1}{2e^{\sqrt{\frac{\mu}{L}}h}-1}=1.
\end{equation*}
Meanwhile, the generalized model is written as 
\begin{equation}\label{eqn:time-XZ}
\begin{split}
    \dot{Z}&=-\frac{1}{\sqrt{\mu}}\nabla f(Z)\\
    Z&=X+\frac{1}{\sqrt{\mu}}\dot{X}
    \end{split}
\end{equation}
Figure \ref{fig:exp-t} (a) displays the error trajectories for \eqref{eqn:ODE-SC-time} and \eqref{eqn:time-XZ}. We observe that the objective function value $f(X)$ generated by \eqref{eqn:time-XZ} decreases monotonically  with $t$, which is expected since $f(Z)$  also  decreases monotonically   (Theorem \ref{thm:conv2}). Moreover, the objective values for \eqref{eqn:ODE-SC} and \eqref{eqn:time-XZ} tend to decrease at similar rates for any $t$.

Furthermore, we consider a discretization of \eqref{eqn:time-XZ}. With $Q(t)=Z(\tau^{-1}(t))$ and $R(t)=X(\tau^{-1}(t))$ for $\tau(t)=t/\sqrt{\mu}$ (see Corollary \ref{cor:Z}), \eqref{eqn:time-XZ} can be  rewritten as 
\begin{equation*}
\begin{split}
    \dot{Q}&=-\nabla f(Q)\\
    Q&=R+\frac{1}{\mu}\dot{R}. 
\end{split}
\end{equation*}
Discretizing this yields the form 
\begin{equation}\nonumber
\begin{split}
    q_{k+1}&=q_k-\nabla f(q_k) \\
    q_{k+1}&=r_{k+1}+\frac{1}{\mu}(r_{k+1}-r_k).
    \end{split}
\end{equation}
This form includes both gradient descent and momentum terms. Figure \ref{fig:exp-t} (b) shows the $L_2$-norm difference of $r_k$ and $R(t=k) =X(\tau^{-1}(k))$. Since the difference is sufficiently small   for any $t$, the discretization is justified. Recalling that gradient descent can be expressed as \eqref{eqn:time-XZ}, we confirm that the  time reparametrization discussions in Section \ref{sec:time} are valid. 

\section{Concluding Remarks}
We have introduced generalized ODE models that are naturally derived from Nesterov's methods and are supported by  strong theoretical results, including accelerated convergence rates, as demonstrated in Theorem~\ref{thm:conv}. Our models can serve as a unified tool to analyze and understand various existing ODE models and Nesterov's methods.
The restart scheme designed using our models ensures a monotonic decrease in objective function values for a broader class of Nesterov's methods.

Notably, the link between our models and standard gradient flow indicates that modifying the speed of time results in  acceleration in continuous time. While this intuitive explanation of acceleration is currently  valid only  in continuous time, it may offer new insights into Nesterov's methods in discrete time. 
For instance, developing a rate-matching discretization scheme for our generalized ODE models could be a promising direction for future research.


\section*{Acknowledgments}
This work was supported in part by the Information and Communications Technology Planning and Evaluation grant funded by MSIT (2022-0-00124, 2022-0-00480).


\appendix
\section{Reformulation of Nesterov's Methods}\label{app:reform}
We define the sequence $\{A_k\}$ as 
\begin{equation*}
    A_k = \frac{A_0}{\prod_{i=0}^{k-1} (1-\theta_i)}
\end{equation*}
for some $A_0 > 0$,  from which we can derive $\theta_k = (A_{k+1} - A_k)A^{-1}_{k+1}$. 
Since $0<\theta_k< 1$, $\{A_k\}$ is a monotonically increasing sequence. 
We can further deduce from~\eqref{eqn:para} that
\begin{align*}
    &\gamma_{k+1}-\mu =(\gamma_k-\mu)(1-\theta_k)= (\gamma_0-\mu)\prod_{i=0}^{k} (1-\theta_i) = \frac{A_0}{A_{k+1}}(\gamma_0-\mu)\\
    &s_k = \frac{\theta_k^2}{\mu+\frac{A_0}{A_{k+1}}(\gamma_0-\mu)} = \frac{(A_{k+1}-A_k)^2}{A_{k+1}(\mu A_{k+1} +A_0(\gamma_0-\mu))}.
\end{align*}
Finally, by letting 
\[
\gamma_0 = 
\begin{cases}
\begin{aligned}
&\frac{1}{A_0}&&\text{if $\mu=0$}\\
&\mu&&\text{if $\mu>0$},
\end{aligned}
\end{cases}
\]
we set $s_k$ and $a_k$ as in \eqref{eqn:s}.

It is noteworthy that \eqref{eqn:NA2} satisfies  an $O(\frac{1}{A_k})$ convergence rate for $0 < s_k \leq \frac{1}{L}$.
The growth of $A_k$, which is critical to the $O(\frac{1}{A_k})$ convergence rate, can be estimated using the condition $0 < s_k \leq \frac{1}{L}$.
When $f$ is convex, assuming that $A_k$ is polynomial, i.e., $A_k = C k^p$, we have
\begin{equation*}
    s_k = \frac{(A_{k+1}-A_k)^2}{ A_{k+1}} = C \frac{(pk^{p-1} + \cdots )^2}{(k+1)^p} \leq \frac{1}{L},
\end{equation*}
which implies that $p \leq 2$ and $A_k = O (k^2)$.
This is consistent with the fact that the convergence rate of Nesterov's accelerated gradient method is $O(1/k^2)$ for  convex objective functions~\citep{nesterov2003introductory}.
When $f$ is $\mu$-strongly convex, the inequality 
\[
s_k = \frac{(A_{k+1}-A_k)^2}{\mu A_{k+1}^2}\leq\frac{1}{L}
\]
should be satisfied, which implies that $A_k =O((1- \sqrt{\frac{\mu}{L}})^{-k})$. Again, this is consistent with the fact that Nesterov's accelerated gradient method achieves an $O((1- \sqrt{\frac{\mu}{L}})^{k})$ convergence rate for strongly convex objective functions~\citep{nesterov2003introductory}.

\section{Proofs}

\subsection{Proof of Theorem~\ref{thm:conv}}\label{app:conv}

Suppose first that $f$ is convex, i.e., $\mu = 0$.
We first claim that $\dot{\mathcal{E}} \leq 0$. 
To show this, we first observe that
\begin{align} \nonumber
    \langle \ac e^{-\alpha}\dot{X}, \nabla f(Y) - \nabla f(X) \rangle=\langle Y-X, \nabla f(Y) - \nabla f(X) \rangle \geq 0,
\end{align}
which implies that
\begin{equation}\label{ieq:XY-CVX}
\langle \dot{X}, \nabla f(Y) \rangle \geq \langle \dot{X}, \nabla f(X) \rangle.
\end{equation}
Therefore, we have
\begin{align} \nonumber
e^{\alpha}\langle Z, \nabla f(Y)\rangle &= e^{\alpha}\langle Y, \nabla f(Y) \rangle + e^{\alpha}\langle Z-Y, \nabla f(Y) \rangle \\ \nonumber
&\geq e^{\alpha} f(Y) +(1-\ac)\langle \dot{X}, \nabla f(Y) \rangle\\ \nonumber
&\geq e^{\alpha} f(X)+e^\alpha \langle Y-X, \nabla f(X)\rangle+(1-\ac)\langle  \dot{X}, \nabla f(Y) \rangle\\ \nonumber
&\geq e^{\alpha} f(X)+ \langle  \dot{X}, \nabla f(X) \rangle\\ \nonumber
&\geq \dot{\beta} f(X) +\frac{d}{dt}f(X),
\end{align}
where the last inequality follows from the  condition~\eqref{eqn:condition}. 
Using the above inequality, we obtain that
\begin{align} \nonumber
\frac{d}{dt}  \mathcal{E} (t)&=\frac{d}{dt}    D_g(0,Z) + \frac{d}{dt}(e^{\beta}f(X)) =\left\langle \frac{d}{dt}\nabla g(Z), Z \right \rangle+ \frac{d}{dt}(e^{\beta}f(X))\\ \nonumber
& = -e^{\alpha+\beta}\langle \nabla f(Y), Z\rangle+ \frac{d}{dt}(e^{\beta}f(X))\\  \nonumber
    &\leq -\dot{\beta}e^{\beta}f(X)-e^{\beta}\frac{d}{dt}f(X)+ \frac{d}{dt}(e^{\beta}f(X))\\ \nonumber
    &=0.
\end{align} 
 Thus, we conclude that
\begin{equation} \nonumber
f(X(t))\leq e^{-\beta(t)} \mathcal{E}(t) \leq e^{-\beta(t)} \mathcal{E}(0).
\end{equation}

Suppose now that $f$ is $\mu$-uniformly convex with respect to $g$.
Again, we show that $\dot{\mathcal{E}} \leq 0$. 
It follows from the $\mu$-uniform convexity of $f$ that 
\begin{equation} \nonumber
\begin{split}
\langle \nabla f(Y), Y \rangle &\geq f(Y) + \mu D_g (0,Y)\\
f(Y) - f(X) 
&\geq \left \langle \nabla f(X), Y-X \right \rangle+ \mu D_g (Y, X)
= \ac \langle \nabla f(X), \dot{X} \rangle + \mu D_g (Y, X).
\end{split}
\end{equation}
Therefore, we have
\begin{align} \nonumber
    \langle e^{\alpha}\nabla f(Y), Z \rangle &= e^{\alpha}\langle \nabla f(Y), Y \rangle + e^{\alpha}\langle \nabla f(Y), Z-Y \rangle \\ \nonumber
    &\geq e^{\alpha}f(Y) + \mu e^{\alpha}D_g (0,Y) + (1-\ac )\langle \nabla f(Y), \dot{X}\rangle\\ \nonumber
    &\geq e^{\alpha}f(X)+ \ac \langle \nabla f(X), \dot{X}\rangle + \mu e^{\alpha}D_g (Y,X) + \mu e^{\alpha}D_g (0,Y) + (1-\ac)\langle \nabla f(Y), \dot{X}\rangle. \nonumber
\end{align}
Using the inequality \eqref{ieq:XY-CVX}, we have
\begin{align} \nonumber
    \langle e^{\alpha}\nabla f(Y), Z \rangle     &\geq e^{\alpha}f(X)+ \ac\langle \nabla f(X), \dot{X}\rangle + \mu e^{\alpha}D_g(Y,X) + \mu e^{\alpha}D_g(0,Y) + (1-\ac)\langle \nabla f(X), \dot{X}\rangle\\ \nonumber
    &= e^{\alpha}f(X)+ \frac{d}{dt}f(X) + \mu e^{\alpha}D_g(Y,X)+\mu e^{\alpha}D_g(0,Y) \\  \nonumber
    & \geq \dot{\beta}f(X) + \frac{d}{dt}f(X)+\mu \dot{\beta}(D_g(0,Y)+ D_g(Y,X)),
\end{align}
where the last inequality follows from the   condition~\eqref{eqn:condition}. 
It follows from the above inequality and the three point identity
\[
\langle \nabla g(Z)-\nabla g(Y), Z \rangle  = D_g(0,Z) + D_g(Z,Y) - D_g(0,Y)
\]
that
\begin{align} \nonumber
    \frac{d}{dt} D_g(0,Z)&= \left \langle\frac{d}{dt}\nabla g(Z),Z \right \rangle \\ \nonumber
    &= -\dot{\beta}\langle \nabla g(Z)-\nabla g(Y),Z\rangle-\frac{e^{\alpha}}{\mu}\langle\nabla f(Y),Z \rangle \\ \nonumber
    &\leq-\dot{\beta}(D_g(0,Z) + D_g(Z,Y) - D_g(0,Y)) - \dot{\beta}\frac{f(X)}{\mu} \\ \nonumber
    &\quad- \frac{d}{dt}\frac{f(X)}{\mu} - \dot{\beta}(D_g(0,Y)+ D_g(Y,X))\\ \nonumber
    &\leq -\dot{\beta}D_g(0,Z) - \dot{\beta}\frac{f(X)}{\mu} - \frac{d}{dt}\frac{f(X)}{\mu},
\end{align}
where the last inequality follows from the non-negativity of Bregman divergence. 
Multiplying both sides by $e^{\beta}$ yields
\begin{equation} \nonumber
\frac{d}{dt}(e^{\beta}D_g(0,Z)) \leq -\frac{d}{dt}\bigg (e^{\beta}\frac{f(X)}{\mu} \bigg ).
\end{equation}
Finally, the time derivative of $\mathcal{E}$ is bounded as follows:
\begin{equation}\nonumber
\begin{split}
\frac{d}{dt} \mathcal{E} (t) 
=\mu \frac{d}{dt} (e^{\beta}   D_g(0,Z)) +  \frac{d}{dt} ( e^{\beta}f(X) )  \leq 0.
\end{split}
\end{equation}
Thus, we conclude that
\begin{equation} \nonumber
f(X(t))\leq e^{-\beta(t)} \mathcal{E}(t) \leq e^{-\beta(t)} \mathcal{E}(0)
\end{equation}
as desired.

\subsection{Proof of Proposition~\ref{thm:Nes-Flow}}\label{app:Nes-Flow}

To prove 
Proposition~\ref{thm:Nes-Flow}, we use the following classical results about nonlinear systems.

\begin{lemma}\citep[Theorem 3.2]{khalil2002nonlinear}\label{lem:unique}
Consider the following ODE 
\[
\dot{\bold{x}}  = F(t, \bold{x})
\]
with $\bold{x} (t_0) = \bold{x}_0 \in \mathbb{R}^n$.
Suppose that $F(t, \bold{x})$ is piecewise continuous in $t$ for each $\bold{x} \in \mathbb{R}^n$ and satisfy
\[
\| F(t, \bold{x}) - F(t, \bold{y}) \| \leq L_1\| \bold{x} -\bold{y} \| \quad \forall \bold{x},\bold{y} \in \mathbb{R}^n, \forall t \in [t_0, t_1]
\]
for some positive constant $L_1$. 
Then, the ODE has a unique solution over $[t_0, t_1]$. 
\end{lemma}

\begin{lemma}\citep[Theorem 3.5]{khalil2002nonlinear} \label{lem:per}
Suppose that the function $F(t, \bold{x})$ be piecewise continuous in $t$ and Lipschitz continuous in $\bold{x}$ on $[t_0, t_1] \times \mathcal{X}$ with a Lipschitz constant $L_1$, where $\mathcal{X} \subset \mathbb{R}^n$ is an open connected set. 
We also assume that  $G(t, \bold{x})$ is uniformly bounded on $[t_0, t_1] \times \mathcal{X}$, i.e.,
\[
\| G(t, \bold{x}) \| \leq L_2 \quad \forall (t, \bold{x}) \in [t_0, t_1] \times \mathcal{X} 
\]
for some positive constant $L_2$. 
Let $\bold{y}(t)$ and $\bold{z}(t)$ be the solutions of
\begin{subequations}
\begin{align}
\dot{\bold{y}} &= F(t, \bold{y}), \quad \bold{y} (t_0) = \bold{y}_0 \label{eqn:unp}\\
\dot{\bold{z}} &= F(t, \bold{z}) + G(t, \bold{z}), \quad \bold{z}(t_0) = \bold{z}_0 \label{eqn:per}
\end{align}
\end{subequations}
such that $\bold{y}(t), \bold{z}(t) \in \mathcal{X}$ for all $t \in [t_0, t_1]$.
Then, we have
\[
\| \bold{y}(t) - \bold{z}(t) \| \leq \| \bold{y}_0 - \bold{z}_0 \| e^{ L_1 (t-t_0)} + \frac{L_2}{L_1} (e^{L_1(t-t_0)} - 1) \quad \forall t \in [t_0, t_1].
\]
\end{lemma}

We now prove Proposition~\ref{thm:Nes-Flow}.
To show that \eqref{eqn:ODE-C} and \eqref{eqn:ODE-SC} have a unique solution, we use Lemma~\ref{lem:unique}.
Eliminating $Y$, \eqref{eqn:ODE-C} and \eqref{eqn:ODE-SC} can be expressed as
\begin{equation}\label{eqn:unp-C}
\begin{split}
 \dot{Z} &= -\dot{\A}(t)\nabla f(X+\ac(t)(Z-X)) \\
\dot{X} &= \frac{\dot{\A}(t)}{\A(t)}(Z-X)
\end{split}
\end{equation}
and
\begin{equation}\label{eqn:unp-SC}
\begin{split}
 \dot{Z} &= -\frac{\dot{\A}(t)}{\A(t)}\bigg ( Z-X-\ac(t)(Z-X) + \frac{1}{\mu}\nabla f(X+\ac(t)(Z-X)) \bigg ) \\
            \dot{X} &= \frac{\dot{\A}(t)}{\A(t)}(Z-X),
\end{split}
\end{equation}
respectively. 
Recall that $\nabla f$ is  Lipschitz continuous, $\A$ is continuously differentiable and $\ac$ is continuous. 
Thus, the condition in Lemma~\ref{lem:unique} is satisfied if $\ac, \dot{\A}$ and $\frac{\dot{\A}}{\A}$ are bounded on $[0,T]$.
Clearly, $\ac$ and $\dot{\A}$ are bounded on $[0,T]$ because they are continuous. 
Regarding $\frac{\dot{\A}}{\A}$, we recall that $\A(t)$ is positive and always increasing. This implies $\frac{\dot{\A}}{\A}$ is also bounded on $[0,T]$.
Therefore, according to Lemma~\ref{lem:unique}, 
\eqref{eqn:ODE-C} and \eqref{eqn:ODE-SC} have a unique solution.

To compare the solution of our ODE model 
and the iterates of Nesterov's method, 
we take  continuously differentiable curves 
$X_h (t), Y_h(t)$ and $Z_h(t)$ such that
\[
X_h (hk) = x_k, \quad Y_h (hk) = y_k, \quad Z_h (hk) = z_k.
\]

To avoid confusion, $O_h$ describes the situation where $h\rightarrow0$ in big $O$ notation. $o_h$ is also the same in small $O$ notation. Using the Tayler expansion, the scheme~\eqref{eqn:C-h} for the convex case can be rewritten as
\begin{equation*}
    \begin{aligned}
            Y_h&=X_h + \ac(t)(Z_h-X_h)\\
            \dot{Z}_h+o_h(1) &= -(\dot{\A}(t)+o_h(1))\nabla f(Y_h) \\
            \dot{X}_h+o_h(1) &= \frac{\dot{\A}(t) + o_h(1)}{ \A(t)} (Z_h-X_h + O_h(h)),
    \end{aligned}
\end{equation*}
and the scheme~\eqref{eqn:SC-h} for the strongly convex case can be expressed as
\begin{align*}
            Y_h&=X_h + \ac(t)(Z_h-X_h)\\
            \dot{Z}_h+o_h(1) &= -\frac{\dot{\A}(t)+o_h(1)}{\A(t)+O_h(h)}\bigg (Z_h-Y_h-\frac{1}{\mu}\nabla f(Y_h) \bigg ) \\
            \dot{X}_h+o_h(1) &= \frac{\dot{\A}(t)}{ \A(t)} (Z_h-X_h + O_h(h))
    \end{align*}
as $h \to 0$. 
Eliminating $Y_h$ yields the following nonlinear systems: for a  convex $f$
\begin{equation}\label{eqn:per-C}
    \begin{split}
            \dot{Z}_h &= -\dot{\A}(t)\nabla f(X_h + \ac(t)(Z_h-X_h)) + o_h(1) \\
            \dot{X}_h &= \frac{\dot{\A}(t)}{\A(t)} (Z_h-X_h)+o_h(1),
    \end{split}
\end{equation}
and for a $\mu$-strongly convex $f$
\begin{equation}\label{eqn:per-SC}
    \begin{split}
            \dot{Z}_h &= -\frac{\dot{\A}(t)}{\A(t)}\bigg (Z_h-X_h-\ac(t)(Z_h-X_h) + \frac{1}{\mu}\nabla f(X_h+\ac(t)(Z_h-X_h)) \bigg ) + o_h(1) \\
            \dot{X}_h &= \frac{\dot{\A}(t)}{\A(t)}(Z_h-X_h) + o_h(1).
\end{split}
\end{equation}
These two systems are perturbed versions of \eqref{eqn:unp-C} and \eqref{eqn:unp-SC}.        
Thus, we can use Lemma~\ref{lem:per} to obtain the following bound on the difference between $X$ and its perturbed counterpart $X_h$:    
\begin{equation*}
    \|X_h(t)-X(t)\|\leq \frac{L_2}{L_1}(e^{L_1 t}-1) \quad \forall t \in [0,T],
\end{equation*}
where $L_1$ denotes the Lipschitz constant of $F$ when \eqref{eqn:unp-C} or \eqref{eqn:unp-SC} is written in the form of \eqref{eqn:unp},
and $L_2$ denotes a constant that bounds $\| G\|$ when \eqref{eqn:per-C} or \eqref{eqn:per-SC} is written in the form of \eqref{eqn:per}.
We also observe that the Lipschitz constant $L_1$ is independent of $h$ since $\A(t)$ is independent of $h$ and $0\leq \ac (t) \leq 1$. 
Furthermore, $L_2 = o_h(1)$ by definition. 
 Thus, we conclude that   
\begin{equation*}
    \lim_{h\rightarrow0} \|X_h(t)-X(t)\| = 0.
\end{equation*}
Considering the instances $t = hk$ for all $0\leq k \leq \left \lfloor \frac{T}{h} \right \rfloor$, the result follows.

\subsection{Proof of Theorem~\ref{thm:conv2}}\label{app:conv2}

We first consider the case of convex $f$, that is, $\mu = 0$. 
We first observe that $f(Z(t))$ is monotonically non-increasing because
\begin{equation}\label{eqn:RQGF-decrease}
\begin{split}
    \frac{d}{dt}f(Z(t))&=\langle \nabla f(Z), \dot{Z}\rangle = -e^{-\alpha-\beta}\langle \nabla^2 g(Z)\dot{Z}, \dot{Z}\rangle \leq 0.
\end{split}
\end{equation}
We now claim that $\dot{\tilde{\mathcal{E}}} \leq 0$. 
To see this, we first consider the time derivative of the Bregman divergence $D_g(0,Z)$:
\begin{align*}
    \frac{d}{dt} D_g(0,Z(t))&= \left \langle Z, \frac{d}{dt}\nabla g(Z) \right \rangle = -e^{\alpha + \beta}\langle Z, \nabla f(Z) \rangle \leq -e^{\alpha + \beta} f(Z) 
    \leq -\dot{\beta}e^{\beta}f(Z), 
\end{align*}
where the last inequality holds due to the  condition~\eqref{eqn:condition}. 
Therefore, 
\begin{align*}
     \frac{d}{dt} \tilde{\mathcal{E}}(t)&\leq -\dot{\beta}e^{\beta}f(Z)+  \frac{d}{dt} (e^{\beta}f(Z)) = e^{\beta} \frac{d}{dt}f(Z)\leq 0,
\end{align*}
where the last inequality follows from \eqref{eqn:RQGF-decrease}.
Therefore, we conclude that
\begin{equation} \nonumber
f(Z(t)) \leq e^{-\beta(t)}\tilde{\mathcal{E}}(t) \leq e^{-\beta(t)}\tilde{\mathcal{E}}(0),
\end{equation}
and the statement holds for $\mu = 0$. 

We now consider the case of  $\mu > 0$. 
As before, $f(Z(t))$ is monotonically non-increasing because
\begin{equation}\label{eqn:decrease}
\begin{split}
    \frac{d}{dt}f(Z)&=\langle \nabla f(Z), \dot{Z}\rangle = -\mu e^{-\alpha}\langle \nabla^2 g(Z)\dot{Z}, \dot{Z}\rangle \leq 0.
    \end{split}
\end{equation}
To show that $\dot{\tilde{\mathcal{E}}} \leq 0$, 
we first bound the time derivative of  $D_g(0,Z)$ as follows:
\begin{align} \nonumber
    \frac{d}{dt} D_g(0,Z(t))&= \left \langle Z, \frac{d}{dt}\nabla g(Z) \right \rangle = -\frac{e^{\alpha}}{\mu}\langle Z, \nabla f(Z) \rangle \leq -\frac{\dot{\beta}}{\mu}f(Z) - \dot{\beta} D_g(0,Z), 
\end{align}
where the last inequality follows from the  condition~\eqref{eqn:condition} and the inequality 
$\langle \nabla f(Z), Z \rangle \geq f(Z) + \mu D_g (0,Z)$ which holds due to the $\mu$-uniform convexity of $f$. 
This gives
\begin{align*}
    \frac{d}{dt} \tilde{\mathcal{E}}(t)&
     \leq -\frac{\dot{\beta}e^{\beta}}{\mu}f(Z) + \frac{d}{dt} \left ( 
     \frac{e^{\beta}}{\mu}f(Z)
     \right ) = \frac{e^{\beta}}{\mu}\frac{d}{dt}f(Z) \leq 0, \nonumber
\end{align*}
where the last inequality follows from \eqref{eqn:decrease}.
Thus, we finally obtain that
\begin{equation} \nonumber
f(Z(t)) \leq \mu e^{-\beta(t)}\tilde{\mathcal{E}}(t) \leq \mu e^{-\beta(t)}\tilde{\mathcal{E}}(0),
\end{equation}
and the result follows.

\subsection{Proof of Corollary~\ref{cor:Z}}\label{app:Z}

Suppose first that $f$ is  convex, i.e., $\mu = 0$. 
Given $\alpha$ and $\beta$ satisfying the condition \eqref{eqn:condition}, let $\tau: [0, \infty) \to [0, \infty)$ be defined as
\[
 \tau(t)= \int_{0}^{t}e^{\alpha(t')+\beta(t')}dt'\geq e^{\beta(t)}-e^{\beta(0)} \geq 0.
 \]
 Then, we have
 \[
 \tau (0) = 0, \quad \dot{\tau} = e^{\alpha + \beta} >0,
 \]
 and $\tau (t) \to \infty$ as $t \to \infty$. Thus, $\tau$ satisfies the condition~\eqref{eqn:condition2}.
On the other hand, given  $\tau$ satisfying the condition \eqref{eqn:condition2}, we choose $\alpha$ and $\beta$ as
\[
{\alpha}=\ln \frac{\dot{\tau}}{\tau+C}, \quad \beta= \ln (\tau+C),
\]
where $C$ is an arbitrary positive constant. 
Then, the condition~\eqref{eqn:condition} holds. 
Therefore, \eqref{eqn:QGF0} with $\tau(t)= \int_{0}^{t}e^{\alpha(t')+\beta(t')}dt'$ is equivalent to the ODE~\eqref{eqn:C-a1}.
It follows from Theorem~\ref{thm:conv} that 
\[
f(Z(t)) \leq O(e^{-\beta(t)}) \leq O\bigg( \frac{1}{\tau (t)} \bigg ).
\]

Suppose now that $\mu > 0$. 
Given $\alpha$ and $\beta$ satisfying the condition \eqref{eqn:condition}, let $\tau: [0, \infty) \to [0, \infty)$ be chosen as
\begin{equation*}
    \tau(t)= \int_{0}^{t}\frac{e^{\alpha(t')}}{\mu}dt'\geq \frac{1}{\mu}(\beta(t)-\beta(0)) \geq 0.
\end{equation*}
Then, 
 \[
 \tau (0) = 0, \quad \dot{\tau}=\frac{e^{\alpha}}{\mu}>0,
 \]
 and $\tau (t) \to \infty$ as $t \to \infty$. 
Thus, the condition~\eqref{eqn:condition2} holds. 
On the other hand, given  $\tau$ satisfying the condition \eqref{eqn:condition2}, we select $\alpha$ and $\beta$ as
\begin{equation*} 
    {\alpha}=\ln (\mu\dot{\tau}), \quad \beta= \mu \tau.
\end{equation*}
Then, $\alpha$ and $\beta$ satisfy the condition~\eqref{eqn:condition}. 
Thus, \eqref{eqn:QGF0} with $\tau(t)= \int_{0}^{t}\frac{e^{\alpha(t')}}{\mu}dt'$ is equivalent to the ODE~\eqref{eqn:SC-a1}.
According to Theorem~\ref{thm:conv}, we obtain the following convergence rate for \eqref{eqn:QGF0}:
\[
f(Z(t)) \leq O(e^{-\beta(t)}) = O( e^{-\mu \tau} )
\]
as desired.

\bibliography{reference}

\begin{thebibliography}{58}
\providecommand{\natexlab}[1]{#1}
\providecommand{\url}[1]{\texttt{#1}}
\expandafter\ifx\csname urlstyle\endcsname\relax
  \providecommand{\doi}[1]{doi: #1}\else
  \providecommand{\doi}{doi: \begingroup \urlstyle{rm}\Url}\fi

\bibitem[Adly and Attouch(2024)]{adly2024accelerated}
Samir Adly and Hedy Attouch.
\newblock Accelerated optimization through time-scale analysis of inertial
  dynamics with asymptotic vanishing and {Hessian-driven} dampings.
\newblock \emph{Optimization}, pages 1--38, 2024.

\bibitem[Alimisis et~al.(2020)Alimisis, Orvieto, Becigneul, and
  Lucchi]{alimisis2020continuous}
Foivos Alimisis, Antonio Orvieto, Gary Becigneul, and Aurelien Lucchi.
\newblock A continuous-time perspective for modeling acceleration in riemannian
  optimization.
\newblock In \emph{Proceedings of the Twenty Third International Conference on
  Artificial Intelligence and Statistics}, pages 1297--1307. PMLR, August 2020.

\bibitem[Allen-Zhu and Orecchia(2014)]{allen2014linear}
Zeyuan Allen-Zhu and Lorenzo Orecchia.
\newblock Linear coupling: An ultimate unification of gradient and mirror
  descent.
\newblock \emph{arXiv preprint arXiv:1407.1537}, 2014.

\bibitem[Attouch et~al.(2016)Attouch, Peypouquet, and Redont]{attouch2016fast}
Hedy Attouch, Juan Peypouquet, and Patrick Redont.
\newblock Fast convex optimization via inertial dynamics with {Hessian} driven
  damping.
\newblock \emph{Journal of Differential Equations}, 261\penalty0 (10):\penalty0
  5734--5783, 2016.

\bibitem[Attouch et~al.(2018)Attouch, Chbani, Peypouquet, and
  Redont]{attouch2018fast}
Hedy Attouch, Zaki Chbani, Juan Peypouquet, and Patrick Redont.
\newblock Fast convergence of inertial dynamics and algorithms with asymptotic
  vanishing viscosity.
\newblock \emph{Mathematical Programming}, 168\penalty0 (1):\penalty0 123--175,
  2018.

\bibitem[Attouch et~al.(2019)Attouch, Chbani, and Riahi]{attouch2019fast}
Hedy Attouch, Zaki Chbani, and Hassan Riahi.
\newblock Fast proximal methods via time scaling of damped inertial dynamics.
\newblock \emph{SIAM Journal on Optimization}, 29\penalty0 (3):\penalty0
  2227--2256, 2019.

\bibitem[Aujol et~al.(2022)Aujol, Dossal, Labarri{\`e}re, and
  Rondepierre]{aujol2022fista}
Jean-Fran{\c{c}}ois Aujol, Charles~H Dossal, Hippolyte Labarri{\`e}re, and Aude
  Rondepierre.
\newblock {FISTA} restart using an automatic estimation of the growth
  parameter.
\newblock 2022.
\newblock URL \url{https://hal.science/hal-03153525v4}.

\bibitem[Baes(2009)]{baes2009estimate}
Michel Baes.
\newblock Estimate sequence methods: extensions and approximations.
\newblock Internal report, Institute for Operations Research, ETH, Z{\"u}rich,
  Switzerland, 2009.

\bibitem[Beck and Teboulle(2009)]{beck2009fast}
Amir Beck and Marc Teboulle.
\newblock A fast iterative shrinkage-thresholding algorithm for linear inverse
  problems.
\newblock \emph{SIAM Journal on Imaging Sciences}, 2\penalty0 (1):\penalty0
  183--202, 2009.

\bibitem[Berthier et~al.(2021)Berthier, Bach, Flammarion, Gaillard, and
  Taylor]{berthier2021continuized}
Rapha{\"e}l Berthier, Francis Bach, Nicolas Flammarion, Pierre Gaillard, and
  Adrien Taylor.
\newblock A continuized view on {Nesterov} acceleration.
\newblock \emph{arXiv preprint arXiv:2102.06035}, 2021.

\bibitem[Betancourt et~al.(2018)Betancourt, Jordan, and
  Wilson]{betancourt2018symplectic}
Michael Betancourt, Michael~I Jordan, and Ashia~C Wilson.
\newblock On symplectic optimization.
\newblock \emph{arXiv preprint arXiv:1802.03653}, 2018.

\bibitem[Bottou et~al.(2018)Bottou, Curtis, and Nocedal]{Bottou2018}
L.~Bottou, F.~E. Curtis, and J.~Nocedal.
\newblock Optimization methods for large-scale machine learning.
\newblock \emph{SIAM Review}, 60\penalty0 (2):\penalty0 223--311, 2018.

\bibitem[Bubeck et~al.(2015)Bubeck, Lee, and Singh]{bubeck2015geometric}
S{\'e}bastien Bubeck, Yin~Tat Lee, and Mohit Singh.
\newblock A geometric alternative to {Nesterov's} accelerated gradient descent.
\newblock \emph{arXiv preprint arXiv:1506.08187}, 2015.

\bibitem[Chen and Luo(2021)]{chen2021unified}
Long Chen and Hao Luo.
\newblock A unified convergence analysis of first order convex optimization
  methods via strong {Lyapunov} functions.
\newblock \emph{arXiv preprint arXiv:2108.00132}, 4\penalty0 (9):\penalty0 10,
  2021.

\bibitem[Chen et~al.(2022{\natexlab{a}})Chen, Shi, and Yuan]{chen2022gradient}
Shuo Chen, Bin Shi, and Ya-xiang Yuan.
\newblock Gradient norm minimization of {Nesterov} acceleration: $ o (1/k^{3})
  $.
\newblock \emph{arXiv preprint arXiv:2209.08862}, 2022{\natexlab{a}}.

\bibitem[Chen et~al.(2022{\natexlab{b}})Chen, Shi, and
  Yuan]{chen2022revisiting}
Shuo Chen, Bin Shi, and Ya-xiang Yuan.
\newblock Revisiting the acceleration phenomenon via high-resolution
  differential equations.
\newblock \emph{arXiv preprint arXiv:2212.05700}, 2022{\natexlab{b}}.

\bibitem[d'Aspremont et~al.(2021)d'Aspremont, Scieur, Taylor,
  et~al.]{d2021acceleration}
Alexandre d'Aspremont, Damien Scieur, Adrien Taylor, et~al.
\newblock Acceleration methods.
\newblock \emph{Foundations and Trends{\textregistered} in Optimization},
  5\penalty0 (1-2):\penalty0 1--245, 2021.

\bibitem[Diakonikolas and Orecchia(2019)]{diakonikolas2019approximate}
Jelena Diakonikolas and Lorenzo Orecchia.
\newblock The approximate duality gap technique: A unified theory of
  first-order methods.
\newblock \emph{SIAM Journal on Optimization}, 29\penalty0 (1):\penalty0
  660--689, 2019.

\bibitem[Dozat(2016)]{dozat2016incorporating}
Timothy Dozat.
\newblock Incorporating {Nesterov} momentum into {Adam}.
\newblock In \emph{4th International Conference of Learning Representation
  Workshop Track}, pages 1--4, May 2016.

\bibitem[Drori and Teboulle(2014)]{drori2014performance}
Yoel Drori and Marc Teboulle.
\newblock Performance of first-order methods for smooth convex minimization: A
  novel approach.
\newblock \emph{Mathematical Programming}, 145\penalty0 (1):\penalty0 451--482,
  2014.

\bibitem[Fazlyab et~al.(2018)Fazlyab, Ribeiro, Morari, and
  Preciado]{fazlyab2018analysis}
Mahyar Fazlyab, Alejandro Ribeiro, Manfred Morari, and Victor~M Preciado.
\newblock Analysis of optimization algorithms via integral quadratic
  constraints: Nonstrongly convex problems.
\newblock \emph{SIAM Journal on Optimization}, 28\penalty0 (3):\penalty0
  2654--2689, 2018.

\bibitem[Fran{\c{c}}a et~al.(2021)Fran{\c{c}}a, Jordan, and
  Vidal]{francca2021dissipative}
Guilherme Fran{\c{c}}a, Michael~I Jordan, and Ren{\'e} Vidal.
\newblock On dissipative symplectic integration with applications to
  gradient-based optimization.
\newblock \emph{Journal of Statistical Mechanics: Theory and Experiment},
  2021\penalty0 (4):\penalty0 043402, 2021.

\bibitem[Hu and Lessard(2017)]{hu2017dissipativity}
Bin Hu and Laurent Lessard.
\newblock Dissipativity theory for {Nesterov's} accelerated method.
\newblock In \emph{International Conference on Machine Learning}, pages
  1549--1557. PMLR, August 2017.

\bibitem[Khalil(2001)]{khalil2002nonlinear}
Hassan~K Khalil.
\newblock \emph{Nonlinear Systems}.
\newblock Pearson, 3rd edition, 2001.

\bibitem[Kim and Yang(2022)]{Kim2022}
Jungbin Kim and Insoon Yang.
\newblock Accelerated gradient methods for geodesically convex optimization:
  Tractable algorithms and convergence analysis.
\newblock In \emph{International Conference on Machine Learning}, pages
  11255--11282. PMLR, July 2022.

\bibitem[Kim and Yang(2023{\natexlab{a}})]{kim2023unifying}
Jungbin Kim and Insoon Yang.
\newblock Unifying {Nesterov's} accelerated gradient methods for convex and
  strongly convex objective functions.
\newblock In \emph{International Conference on Machine Learning}, pages
  16897--16954. PMLR, July 2023{\natexlab{a}}.

\bibitem[Kim and Yang(2023{\natexlab{b}})]{kim2024convergence}
Jungbin Kim and Insoon Yang.
\newblock Convergence analysis of {ODE} models for accelerated first-order
  methods via positive semidefinite kernels.
\newblock In \emph{Advances in Neural Information Processing Systems},
  volume~36, pages 63023--63035, December 2023{\natexlab{b}}.

\bibitem[Krichene et~al.(2015)Krichene, Bayen, and
  Bartlett]{krichene2015accelerated}
Walid Krichene, Alexandre Bayen, and Peter~L Bartlett.
\newblock Accelerated mirror descent in continuous and discrete time.
\newblock In \emph{Advances in Neural Information Processing Systems},
  volume~28, December 2015.

\bibitem[Lessard et~al.(2016)Lessard, Recht, and Packard]{lessard2016analysis}
Laurent Lessard, Benjamin Recht, and Andrew Packard.
\newblock Analysis and design of optimization algorithms via integral quadratic
  constraints.
\newblock \emph{SIAM Journal on Optimization}, 26\penalty0 (1):\penalty0
  57--95, 2016.

\bibitem[Li et~al.(2023)Li, Yuan, Gidel, Gu, and Jordan]{li2023nesterov}
Chris~Junchi Li, Huizhuo Yuan, Gauthier Gidel, Quanquan Gu, and Michael Jordan.
\newblock Nesterov meets optimism: {Rate-optimal} separable minimax
  optimization.
\newblock In \emph{International Conference on Machine Learning}, pages
  20351--20383. PMLR, July 2023.

\bibitem[Maskan et~al.(2023)Maskan, Zygalakis, and
  Yurtsever]{maskan2024variational}
Hoomaan Maskan, Konstantinos Zygalakis, and Alp Yurtsever.
\newblock A variational perspective on high-resolution {ODEs}.
\newblock In \emph{Advances in Neural Information Processing Systems},
  volume~36, December 2023.

\bibitem[Maul{\'e}n and Peypouquet(2023)]{maulen2023speed}
Juan~Jos{\'e} Maul{\'e}n and Juan Peypouquet.
\newblock A speed restart scheme for a dynamics with {Hessian}-driven damping.
\newblock \emph{Journal of Optimization Theory and Applications}, 199\penalty0
  (2):\penalty0 831--855, 2023.

\bibitem[Moucer et~al.(2023)Moucer, Taylor, and Bach]{moucer2023systematic}
C{\'e}line Moucer, Adrien Taylor, and Francis Bach.
\newblock A systematic approach to {Lyapunov} analyses of continuous-time
  models in convex optimization.
\newblock \emph{SIAM Journal on Optimization}, 33\penalty0 (3):\penalty0
  1558--1586, 2023.

\bibitem[Muehlebach and Jordan(2019)]{muehlebach2019dynamical}
Michael Muehlebach and Michael Jordan.
\newblock A dynamical systems perspective on {Nesterov} acceleration.
\newblock In \emph{International Conference on Machine Learning}, pages
  4656--4662. PMLR, June 2019.

\bibitem[Muehlebach and Jordan(2021)]{muehlebach2021optimization}
Michael Muehlebach and Michael~I Jordan.
\newblock Optimization with momentum: Dynamical, control-theoretic, and
  symplectic perspectives.
\newblock \emph{Journal of Machine Learning Research}, 22\penalty0
  (73):\penalty0 1--50, 2021.

\bibitem[Muehlebach and Jordan(2023)]{muehlebach2023accelerated}
Michael Muehlebach and Michael~I Jordan.
\newblock Accelerated first-order optimization under nonlinear constraints.
\newblock \emph{arXiv preprint arXiv:2302.00316}, 2023.

\bibitem[Nemirovskij and Yudin(1983)]{nemirovskij1983problem}
Arkadij~Semenovi{\v{c}} Nemirovskij and David~Borisovich Yudin.
\newblock \emph{Problem complexity and method efficiency in optimization}.
\newblock Wiley-Interscience, 1983.

\bibitem[Nesterov(2005)]{nesterov2005smooth}
Yu~Nesterov.
\newblock Smooth minimization of non-smooth functions.
\newblock \emph{Mathematical Programming}, 103\penalty0 (1):\penalty0 127--152,
  2005.

\bibitem[Nesterov(2008)]{nesterov2008accelerating}
Yu~Nesterov.
\newblock Accelerating the cubic regularization of {Newton's} method on convex
  problems.
\newblock \emph{Mathematical Programming}, 112\penalty0 (1):\penalty0 159--181,
  2008.

\bibitem[Nesterov(2003)]{nesterov2003introductory}
Yurii Nesterov.
\newblock \emph{Introductory lectures on convex optimization: A basic course},
  volume~87 of \emph{Applied Optimization}.
\newblock Springer Science \& Business Media, 2003.

\bibitem[Nesterov(2018)]{nesterov2018lectures}
Yurii Nesterov.
\newblock \emph{Lectures on Convex Optimization}, volume 137 of \emph{Springer
  Optimization and Its Applications}.
\newblock Springer, 2018.

\bibitem[Nesterov(1983)]{nesterov1983method}
Yurii~E Nesterov.
\newblock A method for solving the convex programming problem with convergence
  rate {$O (1/k^2)$}.
\newblock In \emph{Dokl. akad. nauk Sssr}, volume 269, pages 543--547, 1983.

\bibitem[O'donoghue and Candes(2015)]{o2015adaptive}
Brendan O'donoghue and Emmanuel Candes.
\newblock Adaptive restart for accelerated gradient schemes.
\newblock \emph{Foundations of computational mathematics}, 15\penalty0
  (3):\penalty0 715--732, 2015.

\bibitem[Sanz~Serna and Zygalakis(2021)]{sanz2021connections}
Jes{\'u}s~Mar{\'\i}a Sanz~Serna and Konstantinos~C Zygalakis.
\newblock The connections between {Lyapunov} functions for some optimization
  algorithms and differential equations.
\newblock \emph{SIAM Journal on Numerical Analysis}, 59\penalty0 (3):\penalty0
  1542--1565, 2021.

\bibitem[Scieur et~al.(2017)Scieur, Roulet, Bach, and
  d'Aspremont]{scieur2017integration}
Damien Scieur, Vincent Roulet, Francis Bach, and Alexandre d'Aspremont.
\newblock Integration methods and optimization algorithms.
\newblock In \emph{Advances in Neural Information Processing Systems},
  volume~30, December 2017.

\bibitem[Shi et~al.(2019)Shi, Du, Su, and Jordan]{shi2019acceleration}
Bin Shi, Simon~S Du, Weijie Su, and Michael~I Jordan.
\newblock Acceleration via symplectic discretization of high-resolution
  differential equations.
\newblock In \emph{Advances in Neural Information Processing Systems},
  volume~32, December 2019.

\bibitem[Shi et~al.(2022)Shi, Du, Jordan, and Su]{shi2021understanding}
Bin Shi, Simon~S Du, Michael~I Jordan, and Weijie~J Su.
\newblock Understanding the acceleration phenomenon via high-resolution
  differential equations.
\newblock \emph{Mathematical Programming}, 195:\penalty0 79--148, 2022.

\bibitem[Su et~al.(2016)Su, Boyd, and Candes]{su2016differential}
Weijie Su, Stephen Boyd, and Emmanuel~J Candes.
\newblock A differential equation for modeling {N}esterov's accelerated
  gradient method: Theory and insights.
\newblock \emph{Journal of Machine Learning Research}, 17:\penalty0 1--43,
  2016.

\bibitem[Suh et~al.(2022)Suh, Roh, and Ryu]{suh2022continuous}
Jaewook~J Suh, Gyumin Roh, and Ernest~K Ryu.
\newblock Continuous-time analysis of accelerated gradient methods via
  conservation laws in dilated coordinate systems.
\newblock In \emph{International Conference on Machine Learning}, pages
  20640--20667. PMLR, July 2022.

\bibitem[Taylor et~al.(2016)Taylor, Hendrickx, and Glineur]{taylor2017smooth}
Adrien~B. Taylor, Julien~M. Hendrickx, and Fran{\c{c}}ois Glineur.
\newblock Smooth strongly convex interpolation and exact worst-case performance
  of first-order methods.
\newblock \emph{Mathematical Programming}, 161\penalty0 (1-2):\penalty0
  307--345, 2016.

\bibitem[Toyoda et~al.(2024)Toyoda, Nishioka, and Tanaka]{toyoda2024unified}
Mitsuru Toyoda, Akatsuki Nishioka, and Mirai Tanaka.
\newblock A unified {Euler--Lagrange} system for analyzing continuous-time
  accelerated gradient methods.
\newblock \emph{arXiv preprint arXiv:2404.03383}, 2024.

\bibitem[Tseng(2008)]{tseng2008accelerated}
Paul Tseng.
\newblock On accelerated proximal gradient methods for convex-concave
  optimization.
\newblock Technical report, University of Washington, 2008.

\bibitem[Ushiyama et~al.(2023)Ushiyama, Sato, and Matsuo]{ushiyama2023unified}
Kansei Ushiyama, Shun Sato, and Takayasu Matsuo.
\newblock A unified discretization framework for differential equation approach
  with {Lyapunov} arguments for convex optimization.
\newblock In \emph{Advances in Neural Information Processing Systems},
  volume~36, pages 26092--26120, December 2023.

\bibitem[Wibisono et~al.(2016)Wibisono, Wilson, and
  Jordan]{wibisono2016variational}
Andre Wibisono, Ashia~C Wilson, and Michael~I Jordan.
\newblock A variational perspective on accelerated methods in optimization.
\newblock \emph{Proceedings of the National Academy of Sciences}, 113\penalty0
  (47):\penalty0 E7351--E7358, 2016.

\bibitem[Wilson et~al.(2021)Wilson, Recht, and Jordan]{wilson2021lyapunov}
Ashia~C Wilson, Ben Recht, and Michael~I Jordan.
\newblock A {Lyapunov} analysis of accelerated methods in optimization.
\newblock \emph{Journal of Machine Learning Research}, 22\penalty0
  (113):\penalty0 1--34, 2021.

\bibitem[Xie et~al.(2024)Xie, Zhou, Li, Lin, and Yan]{xie2024adan}
Xingyu Xie, Pan Zhou, Huan Li, Zhouchen Lin, and Shuicheng Yan.
\newblock Adan: Adaptive {Nesterov} momentum algorithm for faster optimizing
  deep models.
\newblock \emph{IEEE Transactions on Pattern Analysis and Machine
  Intelligence}, 2024.

\bibitem[Yang et~al.(2022)Yang, Bao, Yuan, Tran, and Zomaya]{yang2022federated}
Zhengjie Yang, Wei Bao, Dong Yuan, Nguyen~H Tran, and Albert~Y Zomaya.
\newblock Federated learning with {Nesterov} accelerated gradient.
\newblock \emph{IEEE Transactions on Parallel and Distributed Systems},
  33\penalty0 (12):\penalty0 4863--4873, 2022.

\bibitem[Zhang et~al.(2018)Zhang, Mokhtari, Sra, and
  Jadbabaie]{zhang2018direct}
Jingzhao Zhang, Aryan Mokhtari, Suvrit Sra, and Ali Jadbabaie.
\newblock Direct {Runge--Kutta} discretization achieves acceleration.
\newblock In \emph{Advances in Neural Information Processing Systems},
  volume~31, December 2018.

\end{thebibliography}
\bibliographystyle{plainnat}

\end{document}